\documentclass[preprint,12pt,authoryear]{elsarticle}
\usepackage{times}

\usepackage{blindtext}
\usepackage{geometry}
\usepackage{amssymb}
\usepackage{amsmath} % for \eqref
\usepackage{amssymb}
\usepackage{booktabs}  % To thicken table lines

\usepackage{graphics}
\usepackage{subfig}
\newtheorem{assum}{Assumption}
\newtheorem{thm}{Theorem}

\newtheorem{remark}{Remark}
\newtheorem{lem}{Lemma}
\usepackage{algorithm2e}

\usepackage{hyperref}
\hypersetup{breaklinks=true, colorlinks, bookmarks=true, citecolor=Black, urlcolor=Violet,linkcolor=Black}

\def\fr{\frac}
\def\pa{\partial}

\def\lam{\lambda}

\def\ba{\begin{align}}
\def\ep{\varepsilon}

\journal{Annual Reviews in Control}

\begin{document}

\begin{frontmatter}



\title{State Estimation of the Stefan PDE:
\\
%\Large 
A Tutorial on Design and Applications to Polar Ice and Batteries
}


%\author[label1,label2]{Shumon Koga, Miroslav Krstic}

\author{Shumon Koga, Miroslav Krstic}

%\affiliation[1]{organization={Department of Electrical and Computer Engineering, University of California San Diego},%Department and Organization
%            addressline={9500 Gilman Drive}, 
%            city={La Jolla},
%            postcode={92093-0411}, 
%            state={CA},
%            country={USA}(e-mail: skoga@ucsd.edu)}
%            
%            
% \affiliation[2]{organization={Department of Mechanical and Aerospace Engineering, University of California San Diego},%Department and Organization
%            addressline={9500 Gilman Drive}, 
%            city={La Jolla},
%            postcode={92093-0411}, 
%            state={CA},
%            country={USA}(e-mail: krstic@ucsd.edu)}

\begin{abstract}
The Stefan PDE system is a representative model for thermal phase change phenomena, such as melting and solidification, arising in numerous science and engineering processes. The mathematical description is given by a Partial Differential Equation (PDE) of the temperature distribution defined on a spatial interval with a moving boundary, where the boundary represents the liquid-solid interface and its dynamics are governed by an Ordinary Differential Equation (ODE). The PDE-ODE coupling at the boundary is nonlinear and creates a significant challenge for state estimation with provable convergence and robustness. 

This tutorial article presents a state estimation method based on PDE backstepping for the Stefan system, using measurements only at the moving boundary.  PDE backstepping observer design generates an observer gain by employing a Volterra transformation of the observer error state into a desirable target system, solving a Goursat-form PDE for the transformation's kernel, and performing a Lyapunov analysis of the target observer error system. 

The observer is applied to models of problems motivated by climate change and the need for renewable energy storage: a model of polar ice dynamics and a model of charging and discharging in lithium-ion batteries. 
The numerical results for polar ice demonstrate a robust performance of the designed estimator with respect to the unmodeled salinity effect in sea ice. The results  for an electrochemical PDE model of a lithium-ion battery with a phase transition material show the elimination of more than 15 \% error in State-of-Charge estimate within 5 minutes even in the presence of sensor noise.    
\end{abstract}


% \begin{highlights}
% \end{highlights}

\begin{keyword}


Stefan system, state estimation, distributed parameter systems, backstepping, nonlinear observer, sea ice, lithium-ion batteries. 
\end{keyword}

\end{frontmatter}



%%%%%%%%%%%%%%%%%%%%%%%%%%%%%%%%%%%%%%%%%%%%%%%%%%%%%%%%%%%%%%%%%%%%%%%%%%%%%%%%%%%%%%%%%%%%%%%%%%%

 \section{Introduction}

\subsection{Phase transitions}

Ice melts into water in a hot environment. Conversely, water freezes into ice in a cold environment. These liquid-solid change phenomena are  called "{\em phase} changes." In addition to various quotidian settings, they also arise, for example, in sea ice in the polar regions \cite{maykut71}, solidification of molten metal in casting \cite{Meng03}, extrusion in polymer 3D-printing \cite{valkanaers2013}, laser sintering for metal 3D-printing \cite{Beaman95}, cryosurgery for cancer treatment \cite{Rabin1998}, thermal energy storage in buildings \cite{Sharma09}, etc. 
 
 The phase change can be regarded as a growth of the domain of one phase, accompanied with a reduction of the domain of the other phase. 
 Apart from the thermal phase change, a growth of a domain, and the consequent movement of the domain's boundary, can be seen in some non-thermal chemical, biological, and social dynamics such as tumor growth in a patient's body \cite{Spangler16}, axonal elongation for neurons' signal transmission \cite{Diehl14,Demir21neuron}, lithiated region in electrodes of lithium-ion batteries \cite{thomas2002mathematical}, crystallization process \cite{braatz2002advanced,nagy2012advances,fujiwara2005first}, domain walls in ferroelectric nanomaterials \cite{mcgilly2015}, spreading of invasive species to a new environment \cite{Du2010speading}, etc.  
 \subsection{Stefan model}

A representative mathematical dynamic model that contains moving boundaries is the ``Stefan system", which is a physical model of thermal phase changes. The associated problem of analyzing and finding the solutions to the Stefan model is referred to as the ``Stefan problem\index{Stefan problem}". Since the phase changes are caused by the temperature dynamics, the physical model involves the temperature of each phase which is distributed in space and changes in time. Hence, the mathematical formulation incorporates partial differential equations (PDEs) defined on a time-varying spatial domain, whereas the dynamics of the position of the moving boundary is governed by an ordinary differential equation (ODE) whose input depends on the PDE state. This configuration gives rise to nonlinear coupling of the PDE and ODE dynamics. As a result, though seemingly consisting of just a linear PDE and a linear scalar ODE, the Stefan problem is mathematically peculiar and not amendable to conventional analysis for PDEs and ODEs.  

The Stefan problem is named after the Slovenian-Austrian physicist Josef (Jo\v{z}ef) Stefan, who was one of the most distinguished and influential physicists of the 19th century, celebrated for his numerous contributions to thermodynamics and heat transfer from the experimental perspective. Perhaps Stefan's name is more recognized for the Stefan-Bolzman law, which revealed that a material with temperature $T$ in absolute unit emits a radiative heat transfer which is proportional to $T^4$, through Stefan's experimental work and his student Ludwig Boltzman's work on the theoretical foundation.  

After the publication of the thermal radiation law, Stefan started to focus of the thickness evolution of polar ice caps motivated by observed data of ice growth and air temperature acquired by British and German explorers during their expeditions. A long time before that, the phase change model by moving boundaries had been studied by Joseph Black in 1762. Franz Neumann developed the solution in his lectures around 1860. However, Neumann's result had not been published until Weber's paper in 1901. Stefan developed his analysis of the solution of ice growth and studied the correspondence with the empirical data, which was published in 1891 \cite{stefan91}. Since then, the model has been known as the ``Stefan problem", and has been studied widely by researchers from the middle 1900s \cite{crepeau_07}.

The mathematical and numerical analysis of the Stefan problem has been widely covered in the literature. The existence, uniqueness, monotone dependence, stability, and the differentiability of the one-phase Stefan problem have been studied by several references, see \cite{Friedman59,Cannon-book,Fasano77,Kolodner56,Rubinstein-book}. The existence and uniqueness of the classical solution of the two-phase Stefan problem were proven in \cite{Cannon71temp,Cannon71flux} with the temperature boundary conditions and the flux boundary conditions, respectively. Several numerical methods to solve the Stefan problem were investigated by finite difference and finite element methods, see \cite{kutluay97,Bonnerot77} for instance. The comparison of the numerical methods was studied in \cite{Javierre06}. 

\subsection{Control of Stefan system}

The work on  control of the Stefan system begun in the last two decades of the 1900s. For instance, feedback control of the Stefan system is developed in \cite{Hoffman82} by an ``on-off" switching design, and an inverse Stefan problem is studied in \cite{Zabaras88} by using an integral method and in \cite{Pawlow87, Pawlow90,Zabaras92,Kang95} by an optimization approach. An optimal feedback control is developed in \cite{Barbu96} by convex analysis, and a feedback control for a thermostats modeled by a hyperbolic Stefan problem is given in \cite{Colli99}. Motion planning by boundary control has been solved in \cite{dunbar2003boundary,dunbar2003motion,Petit10}, overcoming the challenges of the nonlinearity in the Stefan model. Optimal control for the Stefan problem has been developed in \cite{Hinze07,Hinze07flow} using a graph-based approach and in \cite{Bernauer11,Alessandri18} using a level-set formulation. Boundary geomertic control has been proposed in \cite{maidi2014}, and trajectory tracking control for the Vertical Gradient Freeze crystal growth process has been developed in \cite{Ecklebe19,Ecklebe20}. Numerous contributions on the feedback control of the Stefan problem with applications to continuous casting of steels have been achieved by Bentsman and coauthors, see \cite{petrus10,petrus12,petrus14,petrus17,Zhelin18,Zhelin19,Zhelin20}.  

\subsection{Backstepping designs for Stefan system}

Recently, boundary feedback controllers for the Stefan system have been designed via a ``backstepping transformation" which has been used for many other classes of infinite-dimensional systems. The initial article  \cite{Shumon19journal} introduces the designs of a state feedback control law, an observer, and an output feedback law for the one-phase Stefan system by proposing a nonlinear backstepping transformation for the moving-boundary Stefan PDE, and proves the exponential stability of the closed-loop system without imposing  {\em a priori} assumption that the temperature state respects the phase constraints, but by proving that such constraints are actually maintained under proposed feedback. Extensions have been provided in the following articles: \cite{koga_2019delay} develops a control design with time-delay in the actuator and proved delay-robustness, %\cite{koga2019iss} investigated an input-to-state stability of the control of Stefan problem with respect to an unknown heat loss, 
\cite{koga2019twophase} develops a control design for the two-phase Stefan problem, and \cite{koga2021towards} shows  stability of the closed-loop system under  sampled-data control. The backstepping controller and estimator for the Stefan system have been successfully applied to the screw extruder-based polymer 3-D printing \cite{koga2019polymer}, the laser sintering-based metal 3-D printing \cite{koga20laser}, polar ice in the Arctic \cite{koga2019arctic}, lithium-ion batteries \cite{koga2017battery}, and energy storage by paraffin with providing the experimental validation \cite{koga20experiment}. The comprehensive materials can be found in our book \cite{KKbook2021}.  

The state feedback control presented in aforementioned literature, reviewed in \cite{koga2021control}, requires the entire profile of the temperature in the liquid phase as a given information. Some imaging-based sensors such as thermographic camera (a.k.a. infrared camera or IR camera) enables to capture the temperature profile, however, they include relatively high noise as compared to single point thermal sensors, such as thermocouples. Thus, estimating the entire temperature profile given a boundary measurement is a significant task for the implementation of the control algorithm. 

\subsection{State estimators for Stefan model}

The problem of estimating variables of interest given some measured value is widely known as "state estimation". One of the most popular state estimation methods is the ``Kalman filter\index{Kalman filter}" \cite{Kalman60} which is an optimal filter for linear dynamical systems with white Gaussian noise in the model and measurements. Another well known concept is the ``Luenberger observer\index{Luenberger observer}" \cite{Luenberger71}, which stabilizes the estimation error at zero in linear deterministic systems. In finite dimensional systems, the observer gain is designed by means of pole placement or a linear matrix inequality \cite{BoydLMI}. 

%{\color{red}
While the control problem for the Stefan system has been developed as mentioned above, a state estimation for the Stefan system has been focused relatively few. In \cite{petrus17}, an online recalibration method has been proposed for the state estimation of the Stefan problem, where the unknown parameters in the model are calibrated with the discrete-time measurements. In   \cite{pernsteiner2021state}, an state estimation for the latent heat storage system modeled by the Stefan system is developed by the Extended Kalman Filter (EKF) for the reduced-order model of the Stefan system through truncation.   
%}

In this tutorial article, we review state estimation of the Stefan system for the purpose of estimating unknown variables through available measurements. The design is employed by a bacsktepping observer, which is one type of Luenberger observer for the Stefan PDE system, where the observer gain is derived by a PDE backstepping method. After several design methods are introduced, we focus on their applications to polar ice and lithium-ion batteries.  

%%%%%%%%%%%%%%%%%%%%%%%%%%%%%%%%%%%%%%%%%%%%%%%%%%%%%%%%%%%%%%%%%%%%%%%%%%%%%%%%%%%%%%%%%%%%%%%%%%%%%%

\section{Stefan System}

\begin{figure}[t]
\begin{center}
%\centering
\includegraphics[width=0.7\linewidth]{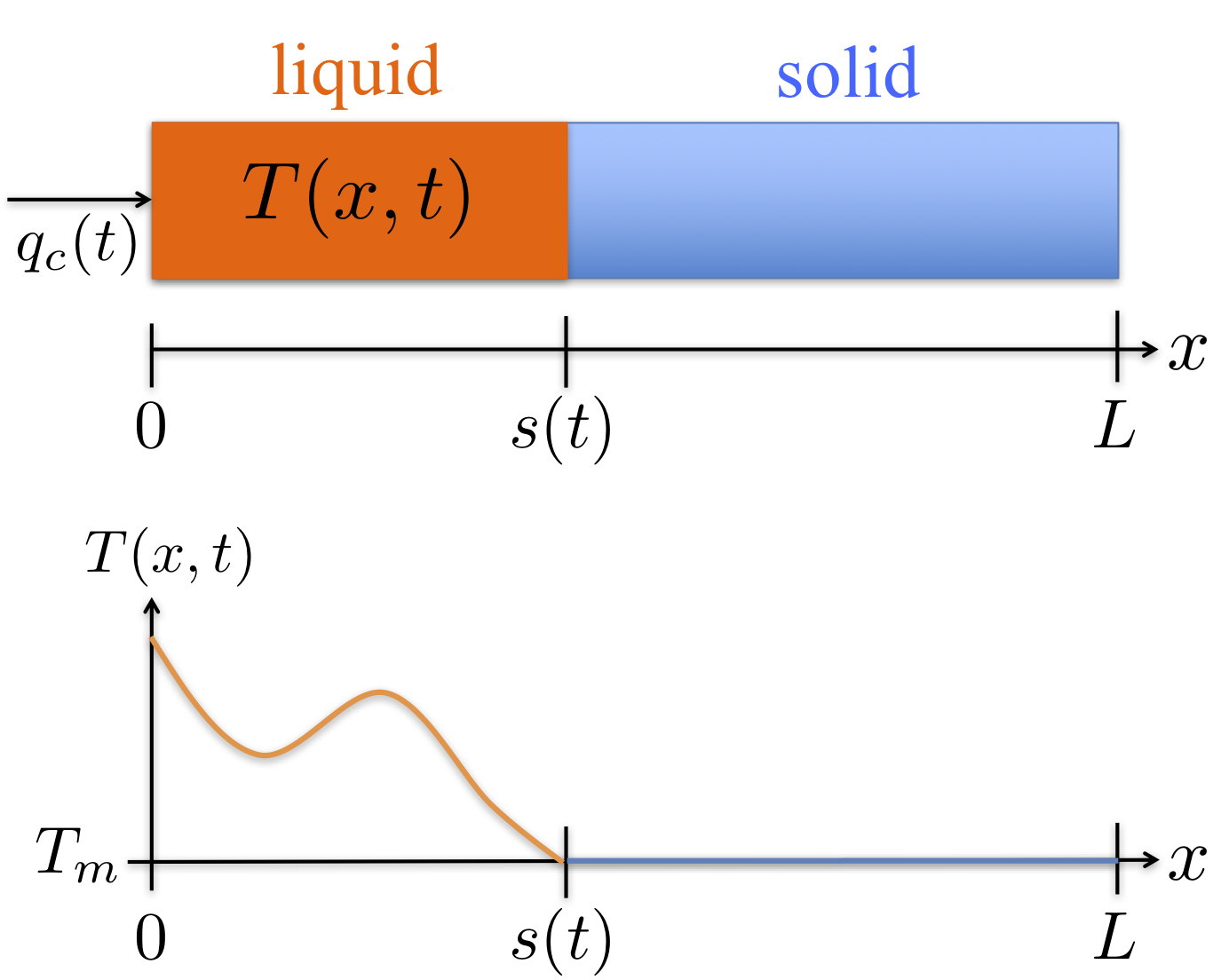}
\vspace{2mm}\\
\caption{Schematic of one-dimensional one-phase Stefan problem. The temperature profile in the solid phase is assumed to be a uniform melting temperature. }
\label{fig:stefan_schematic}
\end{center}
\end{figure}
We begin with a consideration of a pure one-component (liquid only) material in one dimension as depicted  in Fig. \ref{fig:stefan_schematic}, namely, with the assumption that the solid phase happens to be at the freezing temperature (and not cooler than freezing). This assumption is made for simplicity, so that our presentation can focus on controlling just the position of the liquid-solid boundary, and not also the temperature in the ``distal'' solid compartment. 
%Later in the article, in Section \ref{sec-two-phase}, we extend the consideration from a single-phase Stefan system to a two-phase Stefan system. 

The essence of the dynamics of the Stefan system is the evolution in time of the moving liquid-solid interface, which we denote by $s(t)$. The $[0,s(t)]$ portion of the domain is liquid, at a varying temperature $T(x,t)$, whereas the remainder of the domain $[0,L]$, namely, the subdomain $[s(t),L]$ is solid, at the melting/freezing temperature $T_m$, i.e., on the verge of becoming liquid. 

Considering a material which, in its liquid phase, has density $\rho$ and  heat capacity $C_{{\rm p}}$, 
the local energy conservation law  is given by
\begin{align}\label{local-energy}
\rho C_{{\rm p}} T_t(x,t) = -q_x(x,t), \quad x \in (0, s(t)) 
\end{align}
where $q(x,t)$ is a heat flux profile and $T(x,t)$ is a temperature profile. Moreover, the local energy balance at the position of the liquid-solid interface $x=s(t)$ involved with the latent heat leads to the dynamics of the moving boundary\index{moving boundary} 
\begin{align} \label{latent-balance}
\rho \Delta H^* \dot{s}(t) = q(s(t),t) . 
\end{align} 
The thermal conduction for a melting component obeys the well known Fourier's Law 
\begin{align}\label{fourier}
q(x,t) = -k T_x (x,t), \quad x \in[0, s(t)] 
\end{align}
where $k$ is the thermal conductivity.  

At the boundary $x=0$, heat flux enters as an external source which can be manipulated as a controlled variable, denoted as $q_{c}(t)$. The boundary condition  at $x=s(t)$ must maintain the melting temperature $T_{\rm m}$, which is the constant threshold level to cause the phase change from the solid to liquid under a static pressure. 

By combining the energy conservation \eqref{local-energy} and the thermal condition \eqref{fourier} with these two boundary conditions, 
the time evolution of the temperature profile in the material's domain can be obtained by 
\begin{align} \label{ch1:stefanPDE}
T_t(x,t)=&\alpha T_{xx}(x,t), \hspace{3mm} x\in (0, s(t)), \\
\label{ch1:stefancontrol} -k T_x(0,t)=&q_{{\mathrm c}}(t),  \\
 \label{ch1:stefanBC}T(s(t),t)=&T_{{\mathrm m}} , 
 \end{align} 
 where $ \alpha :=\fr{k}{\rho C_{{\rm p}}}$. The initial condition is defined as an arbitrary spatial function for the temperature profile 
 $T(x,0)=T_0(x) > T_{\rm m}$, along with a positive valued interface position $s(0) = s_0$.
 
  If we disregard the dynamics of the moving boundary $s(t)$, the linear PDE model \eqref{ch1:stefanPDE}--\eqref{ch1:stefanBC} may deceive us into thinking that the Stefan model is linear.  However, by combining the latent heat energy balance \eqref{latent-balance} and the thermal conduction \eqref{fourier}, we arrive at the so called "Stefan condition\index{Stefan condition}," given by the following nonlinear ODE
\begin{align}\label{ch1:stefanODE}
 \dot{s}(t)=-\beta T_x(s(t),t), 
\end{align}
where $\beta :=\frac{k}{\rho \Delta H^*}$. 

There are two requirements for the validity of the model \eqref{ch1:stefanPDE}-\eqref{ch1:stefanODE}:
\begin{align}\label{temp-valid}
T(x,t) \geq& T_{{\rm m}}, \quad  \forall x\in(0,s(t)), \quad \forall t>0, \\
\label{int-valid}0 < s(t)<  &L, \quad \forall t>0. 
\end{align}
First, the trivial: the liquid phase is not frozen, i.e., the liquid temperature $T(x,t)$ is greater than the melting temperature $T_{\rm m}$. Second, equally trivially, the material is not entirely in one phase, i.e., the interface remains inside the material's domain. These physical conditions are also required for the existence and uniqueness of  solutions \cite{alexiades2018mathematical}.
Hence, we assume the following for the initial data. 

\begin{assum}\label{ass:initial} 
$0 < s_0 < L$, $T_0(x) \in C^1([0, s_0];[T_{\rm m}, +\infty))$ with $T_0(s_0) = T_{\rm m}$.
 \end{assum}

\begin{lem}\label{lem1}
With Assumption \ref{ass:initial}, if $q_{\rm c}(t)$ is a bounded piecewise continuous non-negative heat function, i.e.,
\begin{align}
q_{{\rm c}}(t) \geq 0,  \quad \forall t\geq 0,  
\end{align}
then there exists a unique classical solution for the Stefan problem \eqref{ch1:stefanPDE}--\eqref{ch1:stefanODE}, which satisfies \eqref{temp-valid}, and 
%for all $t \geq 0$, and
\begin{align} \label{eq:sdot-pos} 
    \dot s(t) \geq 0, \quad \forall t \geq 0. 
\end{align}

\end{lem}

The definition of the classical solution of the Stefan problem is given in 
%can be seen in literature, for instance see 
Appendix A of \cite{Shumon19journal}. The proof of Lemma \ref{lem1} is by maximum principle for parabolic PDEs and Hopf's lemma, as shown in \cite{Gupta03}. 

%%%%%%%%%%%%%%%%%%%%%%%%%%%%%%%%%%%%%%%%%%%%%%%%%%%%%%%%%%%%%%%%%%%%%%%%%%%%%%%%%%%%%%%%%%%%%%%%%%%%%%%%

\section{State Estimator Design by PDE Backstepping} \label{sec:estimator}

In this section, we focus on the state estimation of the Stefan system by a Luenberger-type observer approach with designing the observer gain via the backstepping method.  
%As a comparison, we also introduce the Extended Kalman Filter (EKF) for a reduced-order model of the Stefan system. 

\subsection{Estimation of temperature profile by measuring boundary temperature and interface position} 

\begin{figure}
\begin{center}
\includegraphics[width=0.8\linewidth]{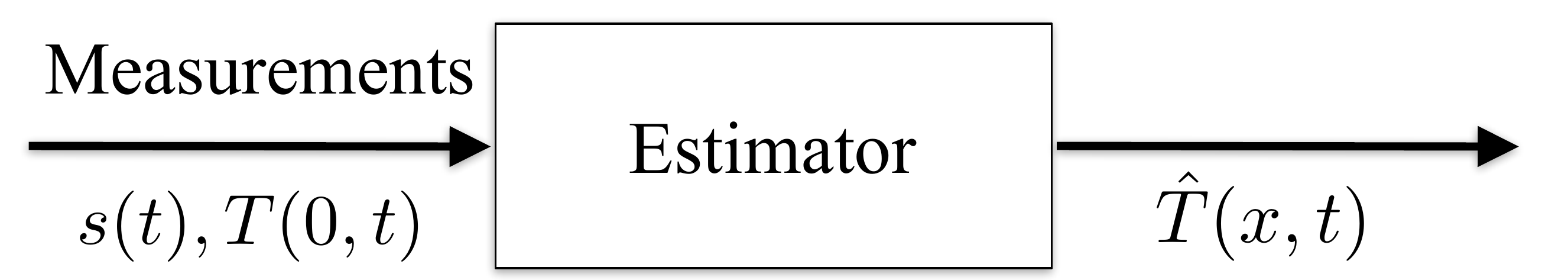}
\caption{The estimation problem measuring a boundary temperature and the interface position.} 
\end{center}
\end{figure} 

\subsubsection{Estimator design and main theorem} 

First, we develop the temperature profile estimation design under the available measurements of a boundary temperature and the interface position. This setup is practical, as the boundary temperature can be measured by a thermocouple and the interface position can be measured by a camera or some optical sensor. For the Stefan system in \eqref{ch1:stefanPDE}--\eqref{ch1:stefanODE}, the following observer is designed with the statement on the theorem. 

\begin{thm}\label {ch34-observer}
Consider the plant \eqref{ch1:stefanPDE}--\eqref{ch1:stefanODE} with the measurements 
\begin{align} 
Y_1(t)=s(t), \quad Y_2(t)=T(0,t), 
\end{align} 
and the following observer 
\begin{align}
 \label{practical-observerPDE}\hat{T}_t(x,t)=&\alpha \hat{T}_{xx}(x,t) + p_1(x,Y_1(t))\left(Y_{2}(t) - \hat{T}(0,t)\right),  \\
 \label{practical-observerBC1}\hat{T}_x(0,t)=& - \frac{q_{{\mathrm c}}(t)}{k} + p_2(Y_1(t)) \left(Y_{2}(t) - \hat{T}(0,t)\right), \\
\label{practical-observerBC2}\hat{T}(Y_1(t),t)=&T_{{\mathrm m}}, 
\end{align}
where $x\in[0,Y_1(t)]$, and the observer gains are 
\begin{align}\label{practical-p1x}
p_1(x,Y_1(t)) =& \lambda Y_1(t)(Y_1(t) - x)  \frac{I_2 \left( \sqrt{ \frac{\lambda}{\alpha} \left\{Y_1(t)^2 - (x-Y_1(t))^2 \right\}}\right)}{ Y_1(t)^2 - (x-Y_1(t))^2 }, \\
p_2(Y_1(t)) = &  -\frac{\lambda}{2\alpha} Y_1(t) \label{practical-p2x}
\end{align}
with a gain parameter $\lambda>0 $. Assume that the model validity condition $T(x,t) \geq T_m$ is satisfied. Then, for all  $\lambda >0$, the observer error system is exponentially stable in the sense of the norm $
||T-\hat{T}||_{{\cal H}_1}^2 $. 
\end{thm}

%\subsection*{Design procedure and proof of Theorem \ref{ch34-observer}} 
Let $\tilde{u}(x,t)=T(x,t)-\hat{T}(x,t)$ be an estimation error variable. Then, we have a system for error variable as
\begin{align}\label{ch3-errorPDE2} 
\tilde{u}_t(x,t)=& \alpha \tilde{u}_{xx}(x,t)-p_1(x,s(t))\tilde{u}(0,t), \quad 0<x<s(t)\\
\label{ch3-errorBC12} \tilde{u}_x(0,t)=&-p_2(x,s(t))\tilde{u}(0,t), \\
\label{ch3-errorBC22} \tilde{u}(s(t),t)=&0
\end{align}
Due to the moving boundary nature, we cannot establish a good target system using the same form of the transformation. Instead, we first consider the observer design for the analogous system on the fixed domain and develop the observer gain with the associated backstepping transformation. After that, we apply the analogous gain and transformation on the moving boundary to the error system \eqref{ch3-errorPDE2}--\eqref{ch3-errorBC22}, and prove the stability of the associated target system on the moving boundary.  

\subsubsection{Fixed domain design} 
Consider the analogous estimation error system on the fixed domain $x \in (0, D)$ given by 
\begin{align}
\tilde{u}_t(x,t)= & \alpha \tilde{u}_{xx}(x,t)-p_1(x,D)\tilde{u}(0,t), \quad 0<x<D\\
\tilde{u}_x(0,t)= & -p_2(D)\tilde{u}(0,t), \\
\tilde{u}(D,t)=& 0 , 
\end{align}
Introduce the transformation 
\begin{align} \label{ch34-DBST}
\tilde{u} (x,t) = w(x,t)+\int_{0}^{x} P(x,y)w(y,t) dy, 
\end{align}
which transforms into 
\begin{align}
w_t(x,t)=&\alpha w_{xx}(x,t)-\lambda w(x,t), \\
w_x(0,t)=&0, \\
w(D,t)= &0. 
\end{align}
Taking time and spatial derivatives of \eqref{ch34-DBST}, the conditions for the gain kernel and the observer gain are obtained by 
\begin{align}
P_{xx}(x,y) -P_{yy}(x,y)= &- \frac{\lambda}{\alpha} P(x,y), \\
P(x,x) =& - \frac{\lambda}{2\alpha} (x-D), \\
 P(D,y) =& 0, \\
 p_1(x,D) =& -\alpha P_{y}(x,0) , \label{ch34-gain1}\\
p_{2}(D) =& -P(0,0) \label{ch34-gain2}. 
\end{align}
The solution to the gain kernel PDE is derived as  
\begin{align}
P(x,y) = \lambda' (D-x) \frac{I_1 \left(\sqrt{ \lambda' \left\{(D-y)^2 - (D-x)^2 \right\}}\right)}{ \sqrt{ \lambda' \left\{(D-y)^2 - (D-x)^2 \right\}}}. 
\end{align}
By the differentiation formula for the Bessel functions\index{Bessel functions}, we have $\frac{d}{dz} \left(\frac{I_1(z)}{z} \right) = \frac{I_2(z)}{z}$. Using this formula and some calculus, the observer gains in \eqref{ch34-gain1} and \eqref{ch34-gain2} are described by 
\begin{align}
p_1(x,D) &=  \alpha \lambda'^2 D(D-x)  \frac{I_2 \left(z \right)}{ z^2 }, \quad z = \sqrt{ \lambda' \left\{D^2 - (D-x)^2 \right\}} , \\
p_2(D) & = - \frac{\lambda}{ 2 \alpha} D. 
\end{align}
Then, in the similar manner, we have 
\begin{align}
w (x,t) = \tilde{u}(x,t)+\int_{0}^{x} Q(x,y)\tilde{u}(y,t) dy, 
\end{align}
which leads to the conditions of  
\begin{align}
Q_{xx}(x,y) - Q_{yy}(x,y) =& \frac{\lambda}{\alpha} Q(x,y) , \\
Q(x,x) =& \frac{\lambda}{2\alpha} (x-D), \\
Q(D,y) =& 0. 
\end{align}
The solution is 
\begin{align}
Q(x,y) =& P(x,y, -\lambda ) = -\lambda' (D-x) \frac{J_1 \left(\sqrt{ \lambda' \left\{(D-y)^2 - (D-x)^2 \right\}}\right)}{ \sqrt{ \lambda' \left\{(D-y)^2 - (D-x)^2 \right\}}} . 
\end{align}

\subsubsection{Analogous observer design on moving boundary domain }
Referring to the result of fixed domain, we apply the backstepping observer design 
\begin{align}
\hat{u}_t(x,t)=&\alpha \hat{u}_{xx}(x,t)+p_1(x,s(t))(u(0,t)-\hat{u}(0,t)), \quad 0<x<s(t)\\
\hat{u}(s(t),t)=&0, \\
\hat{u}_x(0,t)=&-q_c(t)/k + p_2(s(t))(u(0,t)-\hat{u}(0,t)), 
\end{align}
with gains 
\begin{align}
p_1(x,s(t)) = &\frac{ \lambda^2}{\alpha} s(t)(x-s(t))  \frac{I_2 \left( \sqrt{ \frac{\lambda}{\alpha} \left\{s(t)^2 - (x-s(t))^2 \right\}}\right)}{  \sqrt{ \frac{\lambda}{\alpha} \left\{s(t)^2 - (x-s(t))^2 \right\}}}, \\
p_2(s(t)) = &  -\frac{\lambda}{2\alpha} s(t). 
\end{align}
Now, we look at the original model in moving boundary coordinate. Consider the invertible transformation 
\begin{align}
\tilde{w} (x,t) = \tilde{u}(x,t)+\int_{0}^{x} Q(x-s(t),y-s(t)) \tilde{u}(y,t) dy, \\
\tilde{u} (x,t) = \tilde{w}(x,t)+\int_{0}^{x} P(x-s(t),y-s(t) ) \tilde{w}(y,t) dy . 
\end{align}
Then, the target system has the form of 
\begin{align}\label{ch34-tar1} 
\tilde{w}_t(x,t)= &\alpha \tilde{w}_{xx}(x,t)-\lambda\tilde{w}(x,t)  \notag\\
& - \dot{s}(t) \int_{0}^{x} q(\bar{x},\bar{y}) \left(\tilde{w}(y,t)+\int_{0}^{y} P(\bar{y},\bar{z}) \tilde{w}(z,t) dz\right) dy, \\
\tilde{w}(s(t),t)= &0, \\
\tilde{w}_x(0,t)= &0 , \label{ch34-tar3}
\end{align}
where $\bar{x} = x-s(t)$, $\bar{y} = y-s(t)$, and $q(\bar x, \bar y) = Q_{x}(x,y) + Q_{y}(x,y)$. 

We prove that the target $\tilde w$-system in \eqref{ch34-tar1}--\eqref{ch34-tar3} is stable under $\dot s(t) > 0$ and $s(t) < s_r$. Consider the Lyapunov function 
\begin{align}
\tilde{V}_1=\frac{1}{2} || \tilde{w}||^2
\end{align}
The time derivative is given by 
\begin{align}
\dot{\tilde{V}}_1= & -\alpha \int_0^{s(t)} \tilde{w}_{x}(x,t)^2 dx - \lambda  \int_0^{s(t)}\tilde{w}(x,t)^2 dx\notag\\
&-\dot{s}(t)\int_0^{s(t)} \tilde{w}(x,t) \left(\int_{0}^{x} q(\bar{x},\bar{y})\left(\tilde{w}(y,t)+\int_{0}^{y} P(\bar{y},\bar{z}) \tilde{w}(z,t) dz\right) dy\right)dx . \label{ch34-Vtildedot} 
\end{align}
Define 
\begin{align} 
\bar{q} =& \max_{(x,y)\in[0, s_r]} q(\bar{x},\bar{y})^2, \\
 \bar{p} = &\max_{(y,z)\in[0, s_r]} P(\bar{y},\bar{z}). 
 \end{align} 
Applying Young's Cauchy Schwarz inequalities to the second line of \eqref{ch34-Vtildedot} with the help of $\dot s(t)>0 $ and $s(t) \leq s_r$, the following inequality is derived 
 \begin{align}
 \dot{\tilde{V}}_1&\leq -\alpha || \tilde{w}_{x}||^2 - \lambda ||\tilde{w}||^2 +\frac{\dot{s}(t)}{2}\left( 1 +2\bar{q} s_r^2 ( 1+\bar{p}^2 s_r^2 )\right)||\tilde{w}||^2 . 
 \end{align}
Consider 
\begin{align}
\tilde{V}_2=\frac{1}{2}\int_0^{s(t)} \tilde{w}_{x}(x,t)^2 dx . 
\end{align}
The time derivative is obtained by 
\begin{align}\label{ch34-V2dot}
\dot{\tilde{V}}_2&=-\alpha || \tilde{w}_{xx}||^2  -\lambda ||\tilde{w}_{x}||^2 -\frac{\dot{s}(t)}{2}\tilde{w}_{x}(s(t),t)^2+ \dot{s}(t) \int_0^{s(t)} \Phi(w(x,t),s(t),x) , 
\end{align}
where 
\begin{align}  
\Phi : = &\int_0^{s(t)} \tilde{w}_{xx}(x,t) \left( \int_{0}^{x} q(\bar{x},\bar{y})\left(\tilde{w}(y,t)+\int_{0}^{y} P(\bar{y},\bar{z}) \tilde{w}(z,t) dz\right) dy \right) dx. 
\end{align} 
Doing integration by parts twice leads to 
\begin{align}
\Phi =&\tilde{w}_{x}(s(t),t)\left( \int_{0}^{s(t)} q(0,\bar{y}) \left(\tilde{w}(y,t)+\int_{0}^{y} P(\bar{y},\bar{z}) \tilde{w}(z,t) dz\right) dy \right)\notag\\
&+ \tilde{w}(0,t)\left( \frac{d}{dx} \left(\int_{0}^{x} q(\bar{x},\bar{y})\left(\tilde{w}(y,t)+\int_{0}^{y} P(\bar{y},\bar{z}) \tilde{w}(z,t) dz\right) dy\right) \right)|_{x=0} \notag\\
& + \int_0^{s(t)}\tilde{w}(x,t)\left( \frac{d^2}{dx^2} \left(\int_{0}^{x} q(\bar{x},\bar{y})\left(\tilde{w}(y,t)+\int_{0}^{y} P(\bar{y},\bar{z}) \tilde{w}(z,t) dz\right) dy\right)dx \right)  . 
\end{align}
We calculate 
\begin{align}
\frac{d}{dx}& \int_{0}^{x} q(\bar{x},\bar{y})\left(\tilde{w}(y,t)+\int_{0}^{y} P(\bar{y},\bar{z}) \tilde{w}(z,t) dz\right) dy |_{x=0} = q(s(t),s(t))\tilde{w}(0,t)
\end{align}
Here, we see that $q(s(t), s(t)) = Q_x(0,0) + Q_y(0,0) = \left(\frac{d}{dx} Q(x,x) \right)|_{x=0} = \frac{\lambda}{2\alpha} $. Moreover, 
\begin{align} 
& \frac{d^2}{dx^2} \left(\int_{0}^{x} q(\bar{x},\bar{y})\left(\tilde{w}(y,t)+\int_{0}^{y} P(\bar{y},\bar{z}) \tilde{w}(z,t) dz\right) dy\right) \notag\\
= & q(\bar{x},\bar{x})\tilde{w}_{x}(x,t) + ( q(\bar{x},\bar{x})P(\bar x, \bar x)  -  2 q_{\bar x}(\bar{x},\bar{x}) - q_{\bar y}(\bar{x},\bar{x}) ) \tilde{w}(x,t) \notag\\
& -\int_{0}^{x} \left( ( 2 q_{\bar x}(\bar{x},\bar{x})+ q_{\bar y}(\bar{x},\bar{x}))P(\bar{x},\bar{z}) +q(\bar{x},\bar{x}) P_{\bar x}(\bar{x},\bar{z}) + q_{\bar x}(\bar{x},\bar{z}) \right) \tilde{w}(z,t) dz \notag\\
& - \int_{0}^{x} q_{\bar x}(\bar{x},\bar{y}) \left( \int_{0}^{y} P(\bar{y},\bar{z}) \tilde{w}(z,t) dz\right) dy. 
\end{align} 
In addition, we have 
\begin{align} 
 \int_0^{s(t)}  q(\bar{x},\bar{x})\tilde{w}(x,t)\tilde{w}_{x}(x,t) dx  =& - \frac{\lambda }{4 \alpha} \tilde{w}(0,t)^2 + \frac{1}{2} (q_{\bar x}(\bar x, \bar x) + q_{\bar y}(\bar x, \bar x)) \tilde{w}(x,t)^2 dx . 
 \end{align} 
 Thus, 
\begin{align} 
\Phi =&\tilde{w}_{x}(s(t),t)\left( \int_{0}^{s(t)} q(0,\bar{y}) \left(\tilde{w}(y,t)+\int_{0}^{y} P(\bar{y},\bar{z}) \tilde{w}(z,t) dz\right) dy \right)+ \frac{\lambda}{4\alpha} \tilde{w}(0,t)^2 \notag\\
&  + \int_0^{s(t)} ( q(\bar{x},\bar{x})P(\bar x, \bar x)  -  \frac{3}{2} q_{\bar x}(\bar{x},\bar{x}) - \frac{1}{2}q_{\bar y}(\bar{x},\bar{x}) ) \tilde{w}(x,t)^2 dx \notag\\
&+ \int_0^{s(t)}\tilde{w}(x,t) I(\tilde w(x,t),x,s(t)) dx,  \label{ch3-phieq}
\end{align}
where 
\begin{align} 
&I(w(x,t),x,s(t)) \notag\\
= & -\int_{0}^{x} \left( ( 2 q_{\bar x}(\bar{x},\bar{x})+ q_{\bar y}(\bar{x},\bar{x}))P(\bar{x},\bar{z}) +q(\bar{x},\bar{x}) P_{\bar x}(\bar{x},\bar{z}) + q_{\bar x}(\bar{x},\bar{z}) \right) \tilde{w}(z,t) dz \notag\\
& - \int_{0}^{x} q_{\bar x}(\bar{x},\bar{y}) \left( \int_{0}^{y} P(\bar{y},\bar{z}) \tilde{w}(z,t) dz\right) dy. 
\end{align} 
Applying Young's and Cauchy-Schwarz inequality, we can show that there exist positive constants $M_1$, $M_2$, $M_3$ such that  
 \begin{align} 
 \tilde{w}_{x}(s(t),t)\left( \int_{0}^{s(t)} q(0,\bar{y}) \left(\tilde{w}(y,t)+\int_{0}^{y} P(\bar{y},\bar{z}) \tilde{w}(z,t) dz\right) dy \right) & \notag\\
 &\hspace{-30mm} \leq\frac{1}{2}  \tilde{w}_{x}(s(t),t)^2 + M_1 || \tilde w||^2, \\
 \int_0^{s(t)} ( q(\bar{x},\bar{x})P(\bar x, \bar x)  -  \frac{3}{2} q_{\bar x}(\bar{x},\bar{x}) - \frac{1}{2}q_{\bar y}(\bar{x},\bar{x}) ) \tilde{w}(x,t)^2 dx \leq & M_2 || \tilde w||^2, \\
\int_0^{s(t)}\tilde{w}(x,t) I(\tilde w(x,t),x,s(t)) dx  \leq& M_3 || \tilde w ||^2 
\end{align} 
Furthermore, by Agmon's inequality, it holds $ w(0,t)^2 \leq 4 s_r || \tilde w||^2$. Therefore, applying all these inequalities to \eqref{ch3-phieq} leads to 
\begin{align} 
\Phi \leq \frac{1}{2}  \tilde{w}_{x}(s(t),t)^2 + \frac{\lambda s_r }{\alpha}|| \tilde{w}_{x}||^2 + (M_1 + M_2 + M_3) || \tilde w||^2 . \label{Phiineq}
\end{align} 
Applying \eqref{Phiineq} to \eqref{ch34-V2dot}, we arrive at 
\begin{align}
\dot{\tilde{V}}_2\leq & -\alpha || \tilde{w}_{xx}||^2 -\lambda || \tilde{w}_{x}||^2  +  \dot{s}(t) \left ( \frac{\lambda s_r }{\alpha}|| \tilde{w}_{x}||^2 + (M_1 + M_2 + M_3) || \tilde w||^2  \right) 
\end{align}
Thus, defining $\tilde{V} = \tilde{V}_1 + \tilde{V}_2$ and $b =  \frac{\alpha}{4 s_r^2} + \lambda$, we can see that there exists a positive constant $a>0$ such that the following inequality holds 
\begin{align}
\dot{\tilde{V}} \leq  - b \tilde V + a \dot{s}(t) \tilde V, 
\end{align}
from which we conclude Theorem \ref{ch34-observer}.

\subsection{Estimation of Both Temperature Profile and Moving Interface by Measuring Only a Boundary Temperature} \label{sec:3-5}  
\begin{figure}
\begin{center}
\includegraphics[width=0.8\linewidth]{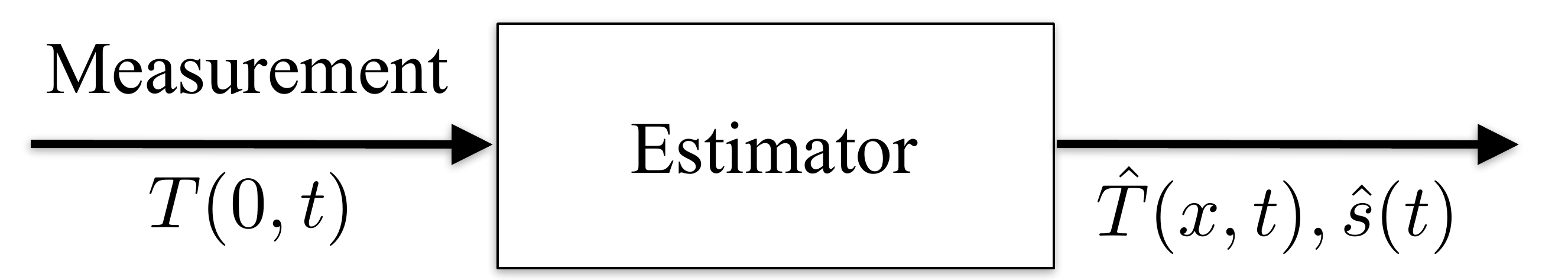}
\caption{The estimation problem measuring only the boundary temperature at heat inlet. The problem is challenging due to the requirement to estimate the interface position. } 
\end{center}
\end{figure} 

The most challenging setup for the state estimation of 1D one-phase Stefan problem is to estimate both temperature profile and moving interface position given a measured value of single boundary temperature. This is not only mathematically challenging, but also practically important in some industrial processes such as steel casting as developed in \cite{petrus12}. We have not established a mathematically rigorous result for this problem, however, we suggest an observer design and investigate the performance in numerical simulation. We consider the plant \eqref{ch1:stefanPDE}--\eqref{ch1:stefanODE} with the measurement
\begin{align} 
Y(t)=T(0,t), 
\end{align} 
and the following PDE observer 
\begin{align}
 \label{challenge-observerPDE}\hat{T}_t(x,t)=&\alpha \hat{T}_{xx}(x,t) + p_1(x,\hat{s}(t))\left(Y(t) - \hat{T}(0,t)\right),  \\
 \label{challenge-observerBC1}\hat{T}_x(0,t)=& - \frac{q_{{\mathrm c}}(t)}{k} + p_2(\hat{s}(t)) \left(Y(t) - \hat{T}(0,t)\right), \\
\label{challenge-observerBC2}\hat{T}(\hat{s}(t),t)=&T_{{\mathrm m}}, 
\end{align}
with the ODE observer 
\begin{align} \label{challenge-observerODE}
\dot{\hat{s}}(t) = - \beta \hat{T}_{x}(\hat{s}(t),t) + l (Y(t) - \hat{T}(0,t)) 
\end{align} 
where $x\in[0,Y_1(t)]$, the observer gains $p_1$, $p_2$ are given in \eqref{practical-p1x}, \eqref{practical-p2x} with a gain parameter $\lambda>0 $, and $l>0$ is a gain parameter. 

\begin{figure}[htb!]
\begin{center} 
\includegraphics[width=0.45\linewidth]{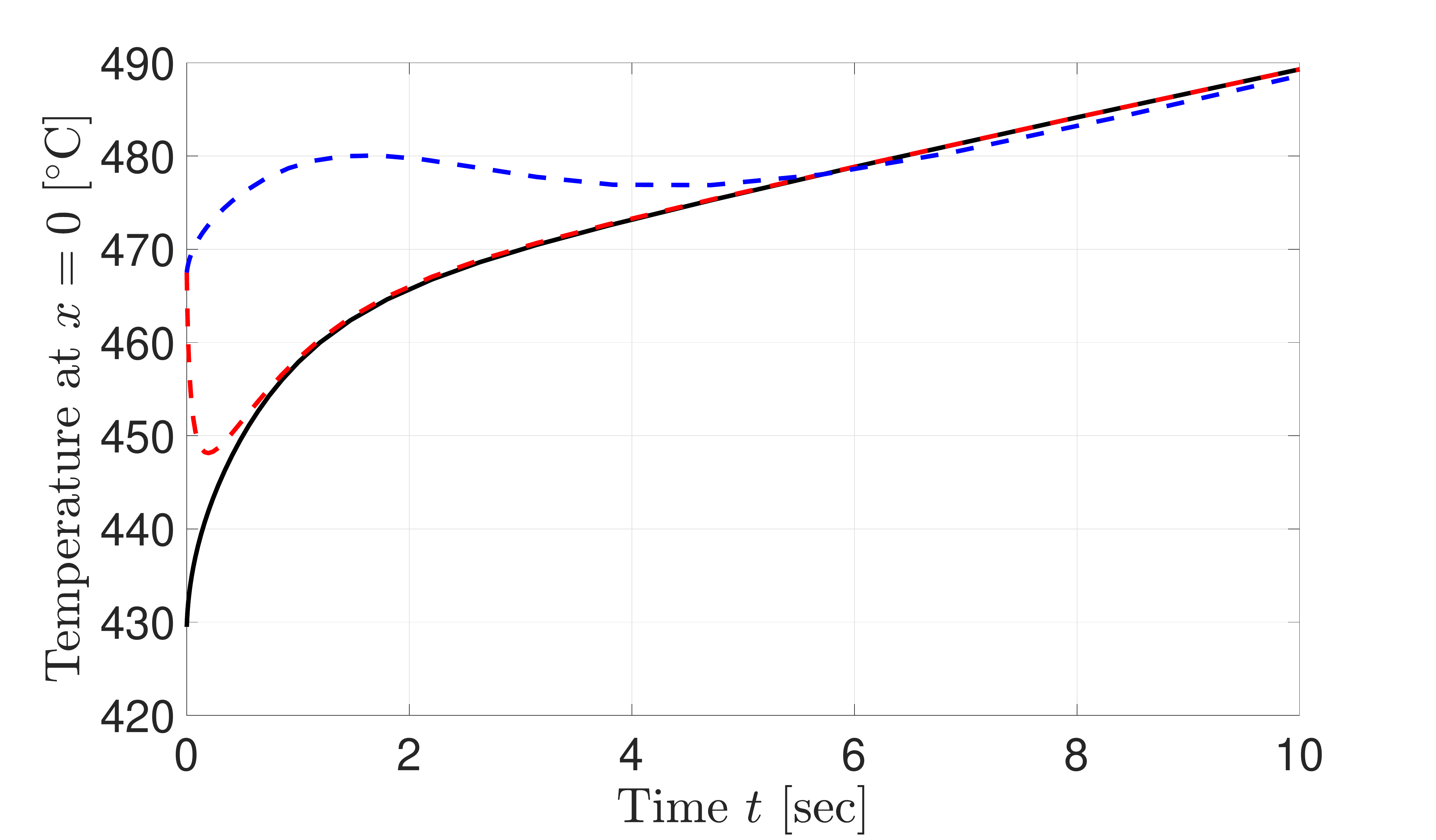}
\includegraphics[width=0.45\linewidth]{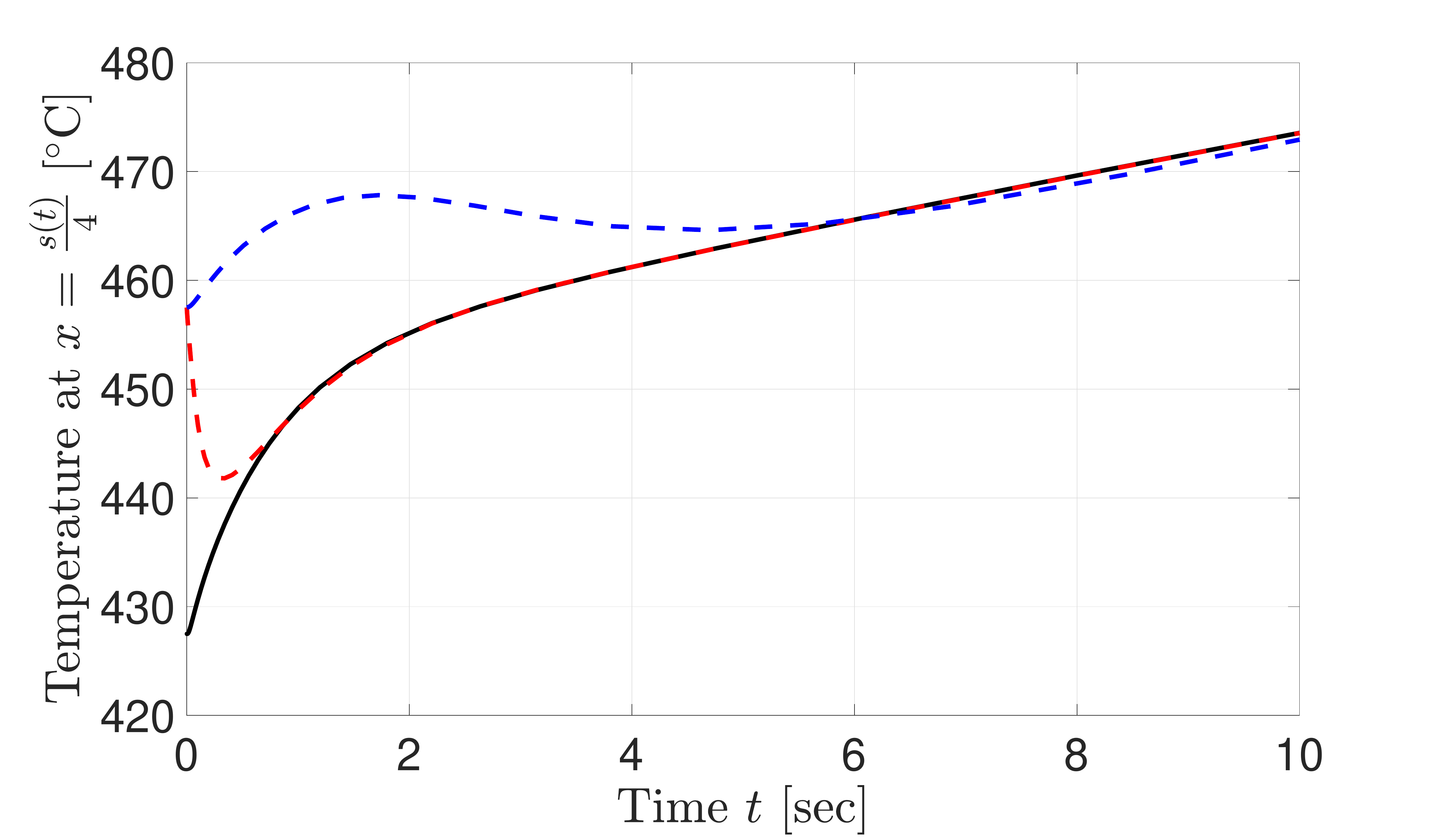}
\includegraphics[width=0.45\linewidth]{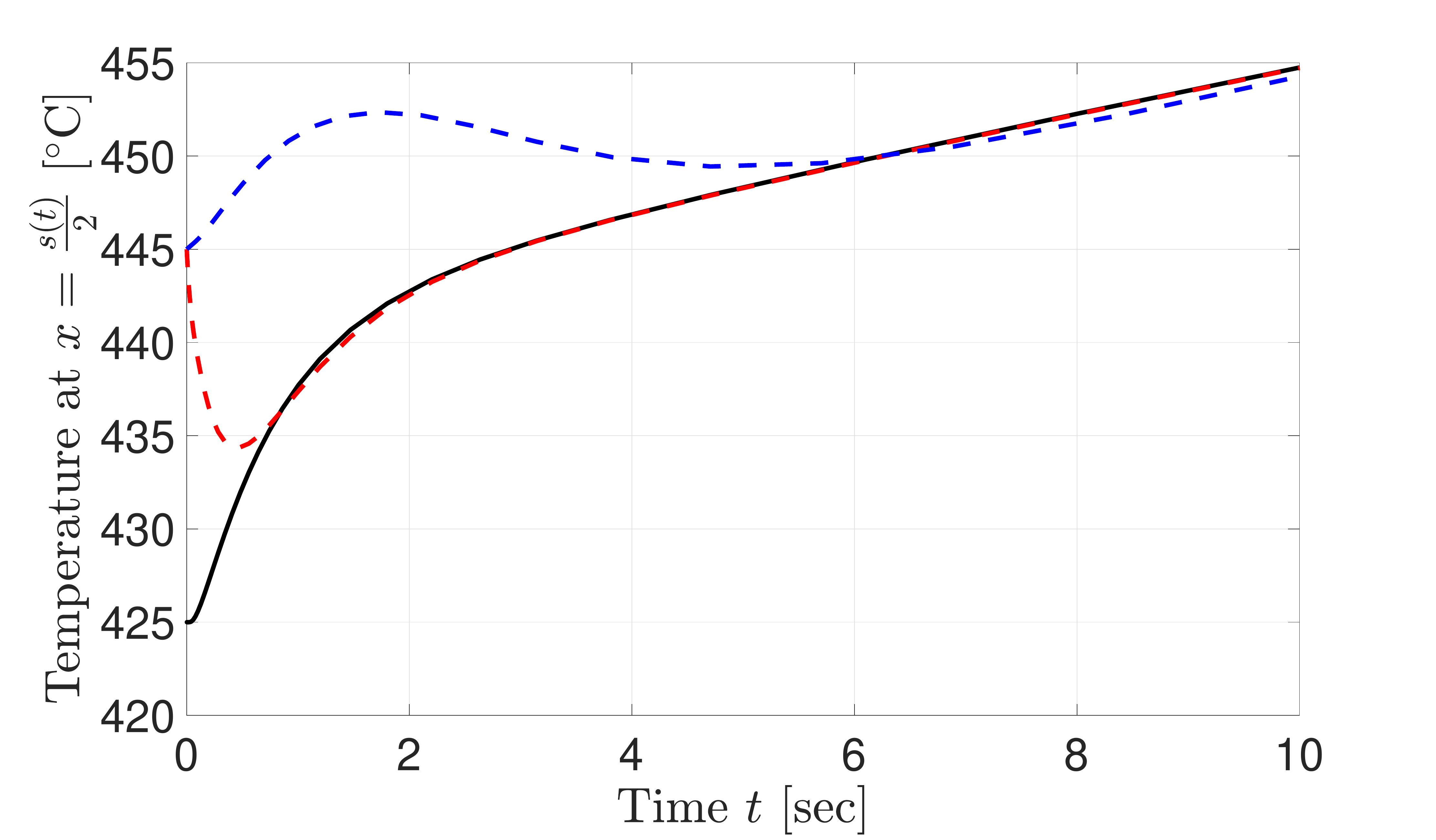}
\includegraphics[width=0.45\linewidth]{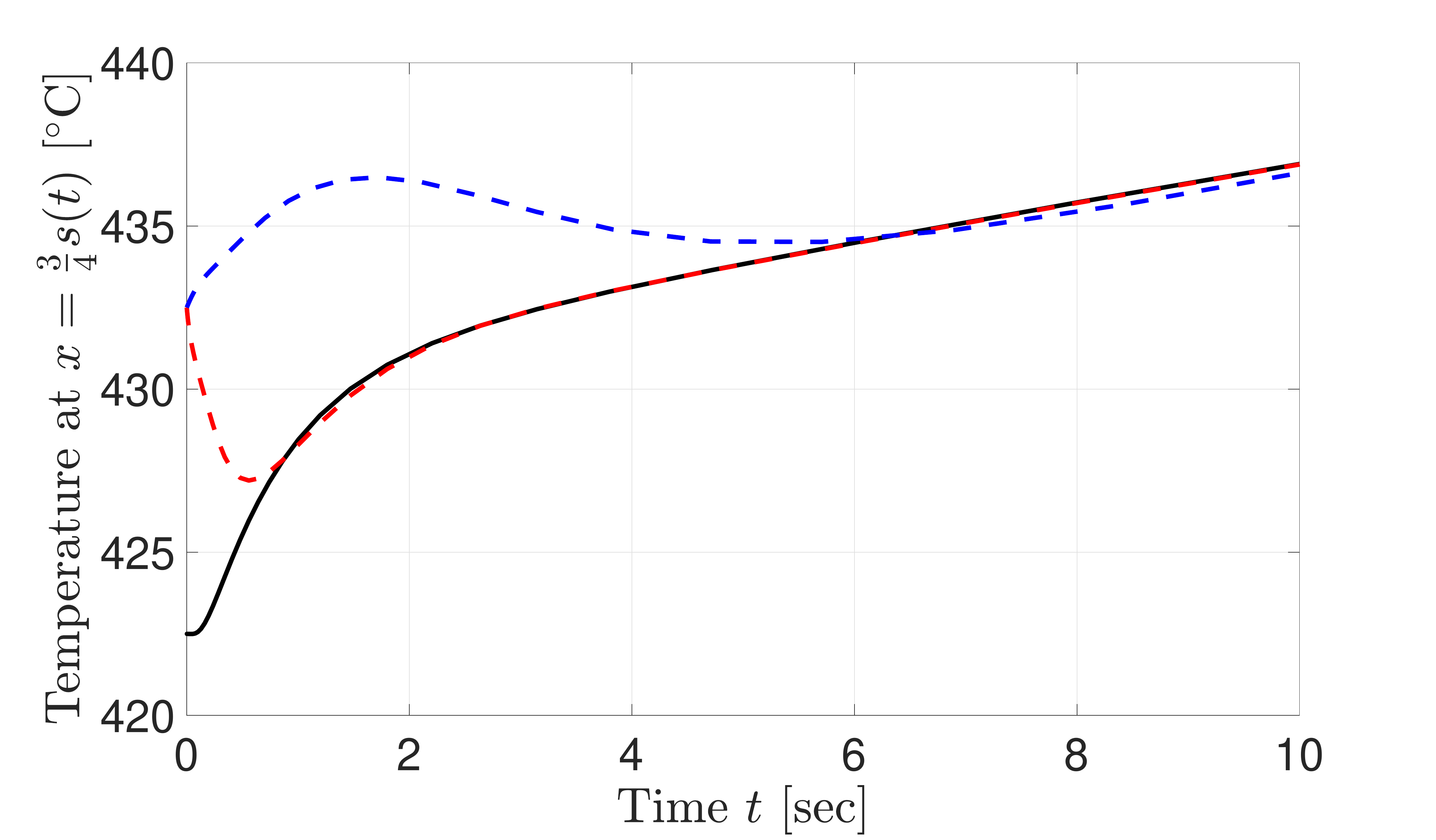}
\includegraphics[width=0.7\linewidth]{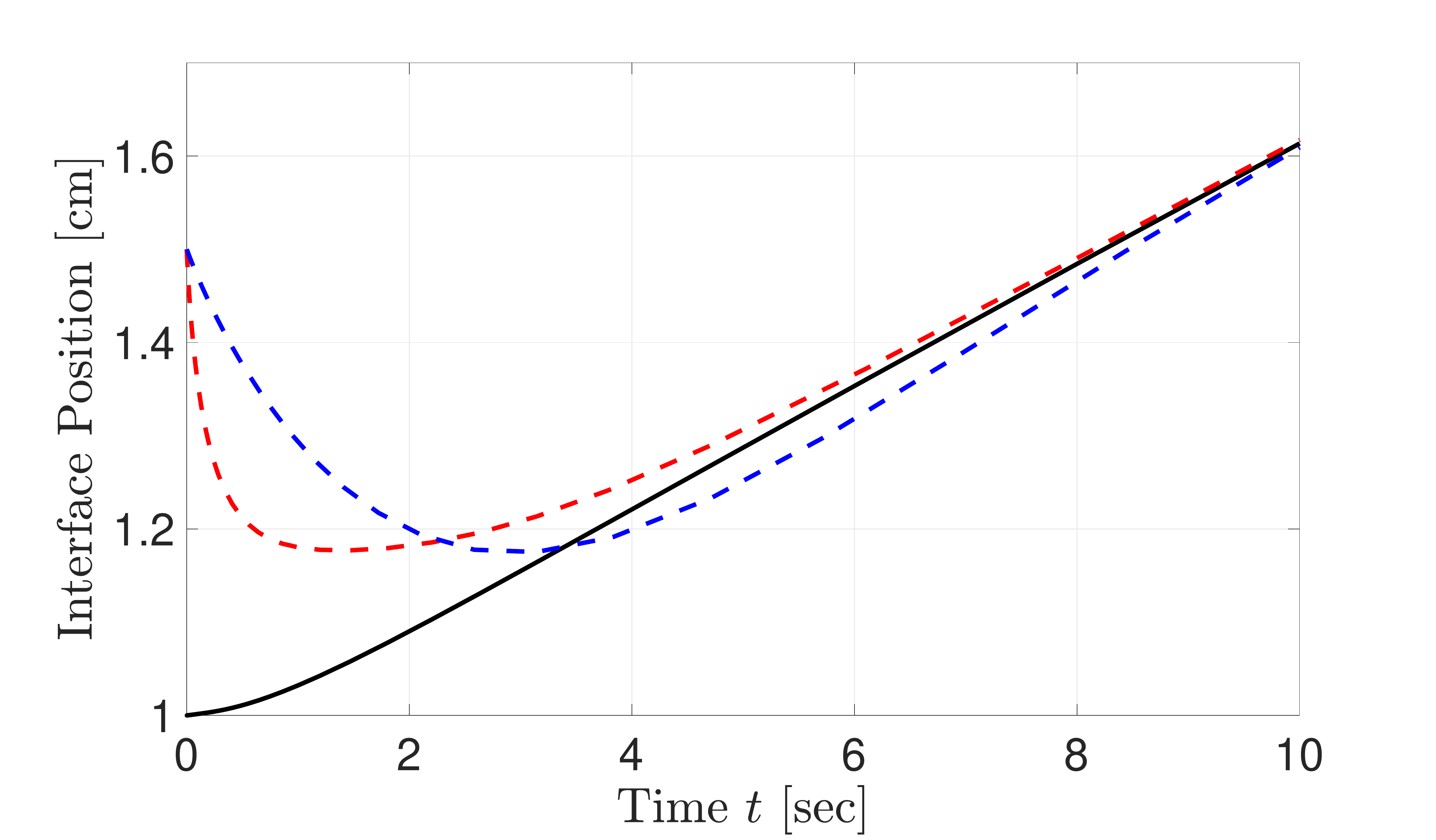}
\caption{Dynamics of the true states (black solid), proposed estimation method (red dash), and estimation using the method in \cite{petrus12} (blue dash), respectively. The upper four figures show the temperature at $x=0, s(t)/4, s(t)/2, 3 s(t)/4$, and the lower figure shows the interface position, respectively. We can observe that for all the states, proposed estimation method achieves faster convergence to the true states.}
\label{fig:ch3_practical} 
\end{center} 
\end{figure}

We study the performance of the observer in numerical simulation using parameters of zinc and tuning the gain parameters $\lam$ and $l$. Also, we compare with the observer design suggested in \cite{petrus12}, which is given by a copy of the plant for PDE observer and the same structure in \eqref{challenge-observerODE} for ODE observer. Fig. \ref{fig:ch3_practical} depicts the simulation results and its comparison, as stated in its caption, which illustrates the better performance of the proposed estimation compared to the method in \cite{petrus12}. 

%%%%%%%%%%%%%%%%%%%%%%%%%%%%%%%%%%%%%%%%%%%%%%%%%%%%%%%%%%%%%%%%%%%%%%%%%%%%%%%%%%%%%%%%%%%%%%%%%%%%

\section{Sea Ice} 

\subsection{Importance of the Arctic Sea Ice for Global Climate Modeling}

The Arctic sea ice\index{Arctic sea ice} has been studied intensively  in the field of climate and geoscience. One of the main reasons is due to ice-albedo feedback which influences climate dynamics through the high reflectivity of sea ice. The other reason is the rapid decline of the Arctic sea ice extent in the recent decade shown in several observations. These observations motivate the investigation of future sea ice amount. Several studies have developed a computational model of the Arctic sea ice and performed numerical simulations of the model with initial sea ice temperature profile. However, the spatially distributed temperature in sea ice is difficult to recover in realtime using a limited number of thermal sensors. Hence, the online estimation of the sea ice temperature profile based on some available measurements is crucial for the prediction of the sea ice thickness. 

A thermodynamic model for the Arctic sea ice was firstly developed in \cite{maykut71} (hereafter MU71), in which the authors investigated the correspondence between the annual cycle pattern acquired from the simulation and empirical data of \cite{unter1961}. The model involves a temperature diffusion equation evolving on a spatial domain defined as the sea ice thickness. Due to melting or freezing phenomena, the aforementioned spatial domain is time-varying.  Such a model is of the Stefan type \cite{Gupta03} and involves a PDE with a state-dependent moving boundary driven by a Neumann boundary value.

Refined models of MU71 have been suggested in literature. For instance, \cite{Semtner1976} proposed a numerical model to achieve faster and accurate computation of MU71 by discretizing the temperature profile into some layers and neglecting the salinity effect. An energy-conserving model of MU71 was introduced in \cite{Bitz1999} by taking into account an internal brine pocket melting on surface ablation and the vertically varying salinity profile. Their thermodynamic model was demonstrated by \cite{Bitz2001} using  a global climate model\index{global climate model} with a  Lagrangian ice thickness distribution. Combining these two models, \cite{Winton2000} developed an energy-conserving three-layer model of sea ice by treating the upper half of the ice as a variable heat capacity layer. 

Remote sensing techniques have been employed to obtain the Arctic sea ice data in several studies. In \cite{Hall2004}, the authors suggested an algorithm to calculate sea ice surface temperature using the satellite measured brightness temperatures, which provided an excellent measurement of the actual surface temperature of the sea ice during the Arctic cold period. The Arctic sea ice thickness data were acquired in \cite{Kwok2009} through a satellite called "ICESat" during 2003-2008 and compared with the data in \cite{Rothrock2008} observed by a submarine during 1958-2000. More recent data describing the evolution of the sea ice thickness have been collected between 2010 and 2014 from the satellite called "CryoSat-2" \cite{Kwok2015}. 

On the other hand, state estimation has been studied as a specific type of data assimilation which utilizes the numerical model along with the measured value. For finite dimensional systems associated with noisy measurements, a well-known approach is the Kalman Filter. Another well-known method is the Luenberger type state observer, which reconstructs the state variable from partially measured variables. For the application to sea ice, \cite{Fenty2013} developed an adjoint-based method as an iterative state and parameter estimation for the coupled sea ice-ocean in the Labrador Sea and Baffin Bay to minimize an uncertainty-weighted model-data misfit in a least-square sense as suggested in \cite{Wunsch2007}, using Massachusetts Institute of Technology general circulation model (MITgcm) developed in \cite{Marshall1997}. In \cite{Fenty2015}, the same methodology was applied to reconstruct the global ocean and ice concentration. Their sea ice model is based on the zero-layer approximation of the numerical model in \cite{Semtner1976}, which is a crude model lacking internal heat storage and promoting fast melting.

\subsection{Thermodynamic Model of Arctic Sea Ice}

\begin{figure}[t]
\begin{center}
\includegraphics[width=0.6\linewidth]{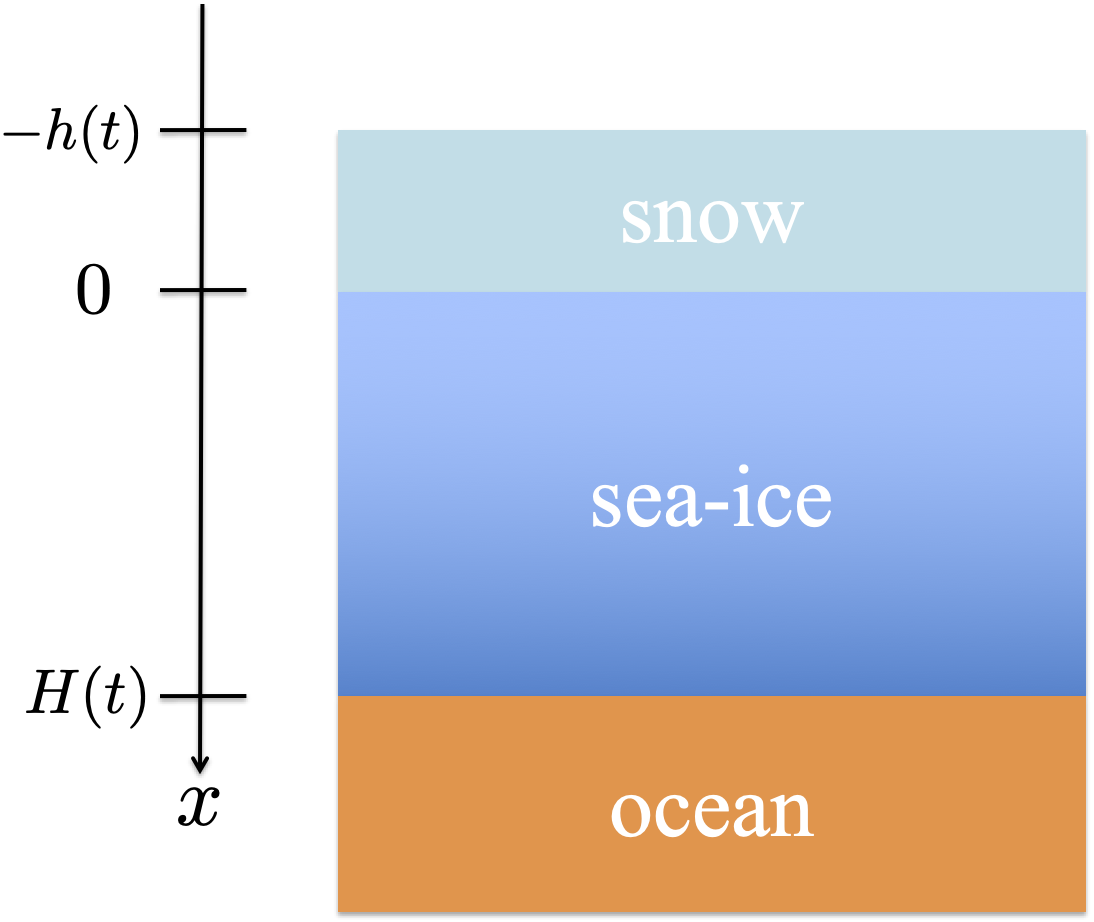}\\
\caption{Schematic of the vertical one-dimensional model of the Arctic sea ice. }
\label{Ch7:fig:seaice}
\end{center}
\end{figure}

The thermodynamic model of MU71 describes the time evolution of the sea ice temperature profile in the vertical axis along with its thickness, which also evolves in time due to accumulation or ablation caused by energy balance.

Fig. \ref{Ch7:fig:seaice} provides a schematic of the Arctic sea ice model. During the seasons other than summer (July and August), the sea ice is covered by snow, and the surface position of the snow also evolves in time. Let $T_{{\rm s}}(x,t)$, $T_{{\rm i}}(x,t)$ denote the temperature profile of snow and sea ice, and $h(t)$ and $H(t)$ denote the thickness of snow and sea ice. The total incoming heat flux from the atmosphere is denoted by $F_{{\rm a}}$, and the heat flux from the ocean is denoted by $F_{{\rm w}}$. The Arctic sea ice model suggested by MU71 gives governing equations of a Stefan-type free boundary problem formulated as 
\begin{align}\label{Ch7:scsys1}
 F_{{\rm a}} - I_0- \sigma (T_{{\rm s}}(-h(t)&,t)+273)^4+ k_{{\rm s}} \frac{\partial T_{{\rm s}}}{\partial x}(-h(t),t)\notag\\
=&
\left\{
\begin{array}{ll}
0, \quad &{\rm if}\quad  T_{{\rm s}}(-h(t),t) < T_{{\rm m}1}, \\
- q \dot{h}(t) , \quad &{\rm if} \quad T_{{\rm s}}(-h(t),t) = T_{{\rm m}1},
\end{array}
\right.\\
\label{Ch7:scsys2} \rho_{{\rm s}} c_0 \frac{\partial T_{{\rm s}}}{\partial t}  (x,t)=& k_{{\rm s}}  \frac{\partial^2 T_{{\rm s}}}{\partial x^2}(x,t), \quad \forall x \in (-h(t),0),\\
\label{Ch7:scsys3} T_{{\rm s}}(0,t) =& T_{{\rm i}}(0,t), \\
\label{Ch7:scsys3-2} k_{{\rm s}} \frac{\partial T_{{\rm s}}}{\partial x}(0,t) =& k_0 \frac{\partial T_{{\rm i}}}{\partial x}(0,t),\\
\label{Ch7:scsys4} \rho c_{{\rm i}} (T_{{\rm i}}, S)\frac{\partial T_{{\rm i}} }{\partial t}(x,t)=& k_{{\rm i}}(T_{{\rm i}}, S) \frac{\partial^2 T_{{\rm i}} }{\partial x^2}(x,t)+ I_0 \kappa_{{\rm i}} e^{-\kappa_{{\rm i}} x}, \quad \forall x \in (0,H(t)),\\
\label{Ch7:scsys5}T_{{\rm i}}(H(t),t) =& T_{{\rm m}2}, \\
\label{Ch7:scsys6}q\dot{H}(t) =& k_{{\rm i}} \frac{\partial T_{{\rm i}}}{\partial x}(H(t),t)-F_{{\rm w}},
\end{align}
where $I_0$, $\sigma$, $k_{{\rm s}}$, $\rho_{{\rm s}}$, $c_0$, $k_0$, $\rho$, $T_{{\rm m}1}$, $T_{{\rm m}2}$, and $q$ are solar radiation penetrating the ice, Stefan-Boltzmann constant, thermal conductivity of snow, density of snow, heat capacity of pure ice, thermal conductivity of pure ice, density of pure ice, melting point of surface snow, melting point of bottom sea ice, and latent heat of fusion, respectively.  The total heat flux from the air is given by 
\begin{align} 
F_{{\rm a}} = (1-\alpha ) F_{{\rm r}}  + F_{{\rm L}}  + F_{{\rm s}} +F_{{\rm l}}, 
\end{align} 
where $F_{{\rm r}}$, $F_{{\rm L}}$, $F_{{\rm s}}$, $F_{{\rm l}}$, and $\alpha $ denote the incoming solar short-wave radiation, the long-wave radiation from the atmosphere and clouds, the flux of sensible heat, the latent heat in the adjacent air, and the surface albedo, respectively.  The heat capacity and thermal conductivity of the sea ice are affected by the salinity as 
\begin{align}\label{Ch7:saltcapa}
c_{{\rm i}}(T_{{\rm i}}, S(x)) = c_{0} +  \frac{ \gamma_1S(x)}{T_{{\rm i}}(x,t)^2}, \quad k_{{\rm i}} (T_{{\rm i}}, S(x)) = k_{0} + \frac{  \gamma_2 S(x)}{T_{{\rm i}}(x,t) },
\end{align}
where $S(x)$ denotes the salinity in the sea ice. $\gamma_1$ and $\gamma_2$ represent the weight parameters. 
The thermodynamic model \eqref{Ch7:scsys1}-\eqref{Ch7:scsys6} allows us to predict the future thickness $(h(t),H(t))$ and the temperature profile $(T_{{\rm s}}, T_{{\rm i}})$ given the accurate initial data. However, from a practical point of view, it is not feasible to obtain a complete temperature profile due to a limited number of thermal sensors. To deal with the problem, the estimation algorithm is designed so that the state estimation converges to the actual state starting from an  initial estimate.

\subsection{Annual Cycle Simulation of Sea Ice Thickness}

For the computation, we use boundary immobilization method and finite difference semi-discretization \cite{kutluay97} with 100-point mesh in space, and the resulting approximated ODEs are calculated by using MATLAB ode15 solver.

%\subsection{Input Parameters}
\subsubsection*{Input Parameters} 
The input parameters are taken from  \cite{maykut71} in SI units and  Table \ref{Ch7:table:1} shows the monthly averaged values of heat fluxes coming from the atmosphere for each month. Table \ref{Ch7:table:2} shows the physical parameters of snow and sea ice. Following \cite{Bitz1999}, the salinity profile is described by
\begin{align} 
S(x) = A\left[1 - {\rm cos} \left\{ \pi \left(\frac{x}{H(t)}\right)^{\frac{n}{m+\frac{x}{H(t)}}} \right\}\right], 
\end{align} 
where $A = 1.6$, $n = 0.407$, and $m = 0.573$.

\begin{table}[htb!]
\begin{center}
\caption{Average monthly values for the energy fluxes. }

    \begin{tabular}{|c |c c c c c|}
    \hline
     $\textbf{Symbol}$&$F_{{\rm r}}$  & $F_{{\rm L}}$ & $F_{{\rm s}}$ & $ F_{{\rm l}}$ & $\alpha$ \\ \hline
    $\textbf{Unit}$ & W/m$^2$ & W/m$^2$ & W/m$^2$ & W/m$^2$ & \\ \hline
    $\textbf{Jan.}$ & 0 & 168 & 19.0 & 0 & $\cdots$ \\
    $\textbf{Feb.}$ & 0 &  166 & 12.3 & -0.323 & $\cdots$\\
    $\textbf{Mar.}$ & 30.7 & 166 & 11.6 & -0.484 & 0.83 \\
    $\textbf{Apr.}$ &  160& 187 & 4.68 & -1.45& 0.81 \\ 
    $\textbf{May.}$ & 286 &  244 &-7.26 & -7.43 & 0.82 \\
    $\textbf{Jun.}$ & 310 & 291 & -6.30 & -11.3 &  0.78 \\
    $\textbf{Jul.}$ &  220 & 308 &-4.84 & -10.3 & 0.64 \\
    $\textbf{Aug.}$ & 145 & 302 & -6.46 & -10.7 &  0.69 \\
    $\textbf{Sep.}$ & 59.7 & 266 & -2.74 &  -6.30 & 0.84 \\
    $\textbf{Oct.}$ & 6.46 & 224 & 1.61 & -3.07 & 0.85 \\
    $\textbf{Nov.}$ & 0 & 181 & 9.04 &  -0.161 & $\cdots$ \\
    $\textbf{Dec.}$ & 0 & 176& 12.8 & -0.161 & $\cdots$ \\ \hline
    \end{tabular}
     \label{Ch7:table:1}
%    \end{center}
 %   \end{table}
  %  \begin{table}[t]
  %  \begin{center}
   \\
\caption{Physical parameters of snow and sea ice. }
    \begin{tabular}{|clc|c|c|}
    \hline
     $\textbf{Symbol}$&  $\textbf{Meaning}$ & $\textbf{Unit}$ & $\textbf{Value}$ \\ \hline
     $\rho_{{\rm s}}$ & density (snow) & kg/m$^3$ & 330 \\ 
    $k_{s}$ & conductivity (snow) & W/m/$^{\circ}$C & 0.31\\
    $\rho$ & density (ice) & kg/m$^3$ & 917 \\ 
    $ c_{0}$ & heat capacity (ice) & J/kg/$^{\circ}$C & 2110  \\  
    $k_{0}$ & conductivity (ice) & W/m/$^{\circ}$C & 2.034\\ 
    $\gamma_1$ & weight of heat capacity & kJ $^{\circ}$C/kg & 18.0 \\ 
    $\gamma_2$ & weight of conductivity & W/m &  0.117\\ 
     $I_0$ & solar radiation & W/m$^2$ & 1.59 \\ 
    $\kappa_{{\rm i}}$ & penetration rate  & /m & 1.5 \\ 
    $T_{{\rm m}1}$ & melting temperature of sea ice at surface& $^{\circ}$C & -0.1 \\
    $T_{{\rm m}2}$ & melting temperature of sea ice at bottom& $^{\circ}$C & -1.8 \\  \hline
    \end{tabular}
\label{Ch7:table:2}
\end{center}
\end{table}

\begin{figure*}[htb!]
\begin{center}
\subfloat[Thickness evolution of the snow and sea ice. ]{\includegraphics[width=0.8\linewidth]{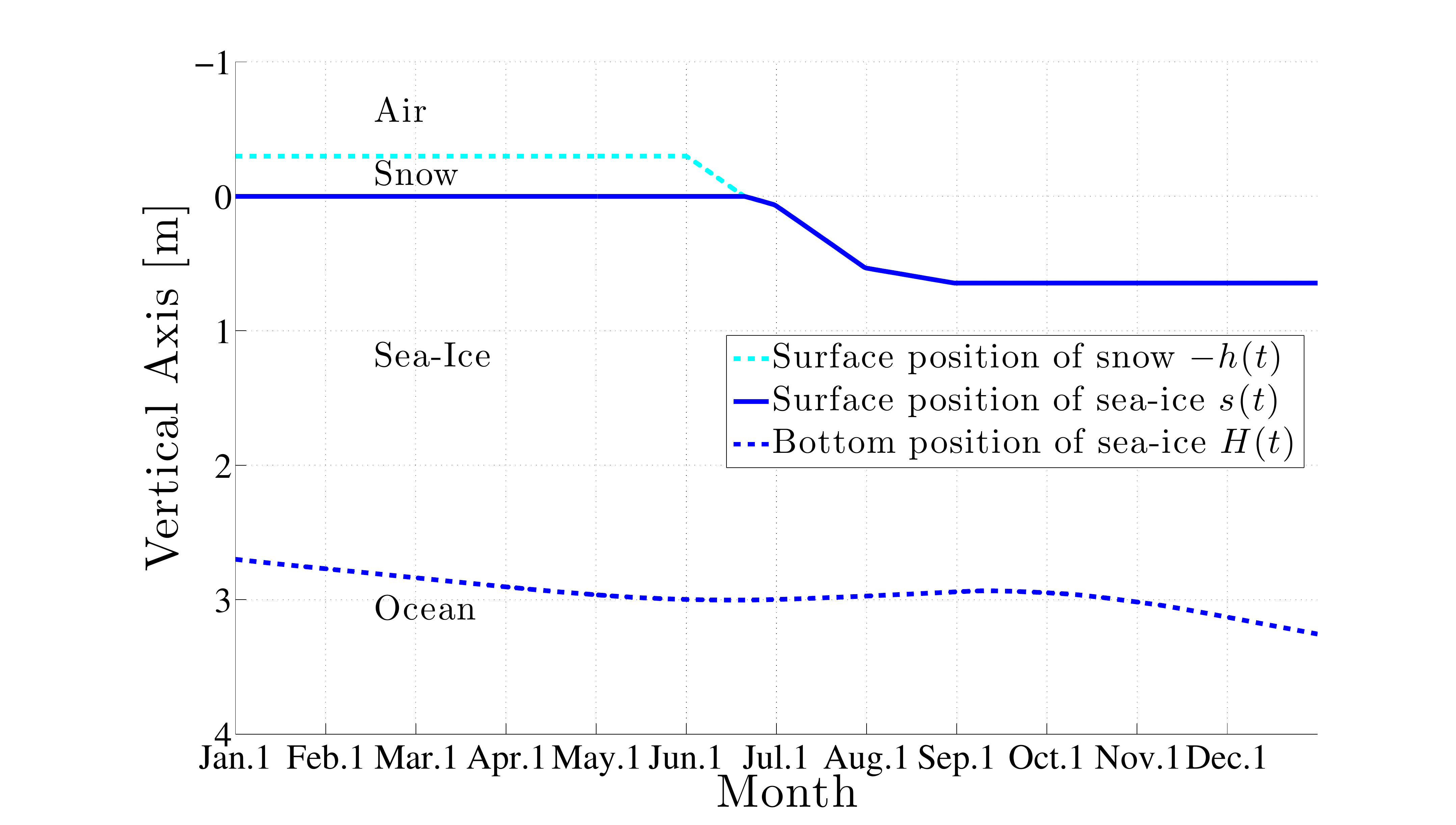}}\\
\subfloat[Time evolution of temperature profile in sea ice.  ]{\includegraphics[width=0.8\linewidth]{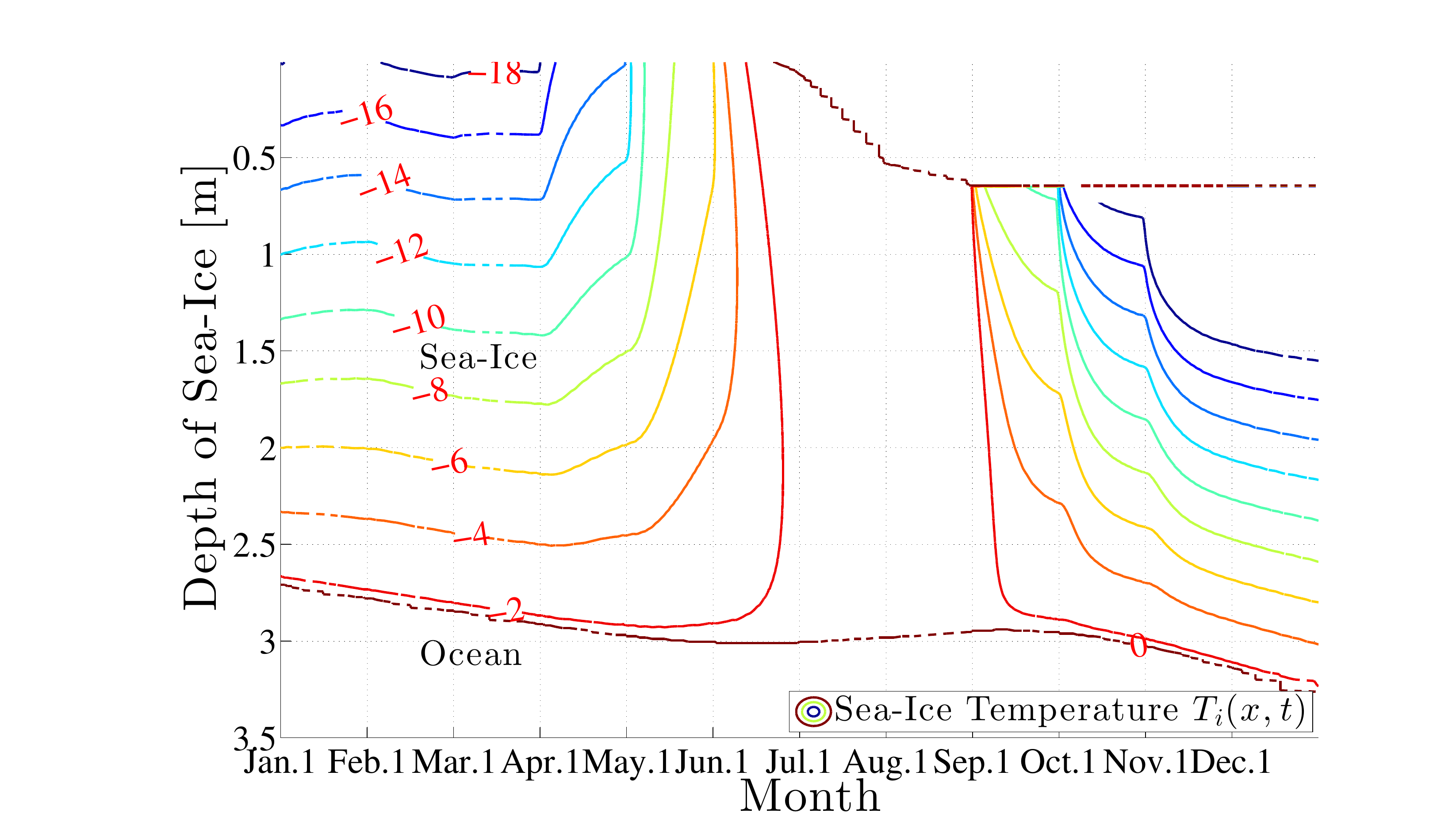}}
\caption{Simulation tests of the plant \eqref{Ch7:scsys1}--\eqref{Ch7:scsys6} on annual cycle. Both (a) and (b) are in good agreement with the simulation results in \cite{maykut71}.}
\label{Ch7:fig:1}
\end{center}
\end{figure*}

%\subsection{Simulation Test of MU71}
\subsubsection*{Simulation Test of MU71} 
Using the given data, firstly the simulation of \eqref{Ch7:scsys1}-\eqref{Ch7:scsys6} is performed and showed in Fig. \ref{Ch7:fig:1} to recover the evolution of $h(t)$ and $H(t)$ in the annual season as in \cite{maykut71}. The dynamic behavior of the snow surface and the bottom of sea ice are shown in Fig. \ref{Ch7:fig:1} (a), and the time evolution of the temperature profile in sea ice is illustrated in Fig. \ref{Ch7:fig:1} (b). We can see that both of Fig. \ref{Ch7:fig:1} (a) and (b) have a good agreement with the simulation results shown in \cite{maykut71}. 

\subsection{Temperature Profile Estimation}

In this section, we derive the estimation algorithm utilizing some available measurements and show the exponential convergence of the designed estimation to a simplified sea ice model. The ice thickness and surface temperature are measured in several studies \cite{Hall2004,Kwok2009,Rothrock2008}. It is indeed typical to check observability\index{observability} before observer design, at least for systems on a constant domain (see \cite{Moura2014} for instance). 
%To study the local observability, one may (i) approximate the PDEs by ODEs (nonlinear for the sea-ice model) using the finite difference method, (ii) linearize the resulting ODEs and the output around the reference states, producing a constant system's matrix $A$ and measurement matrix $C$, (iii) define the observability matrix ${\mathcal O}:=(C; CA; CA^2; \cdots; CA^{n-1} )$, and check if $rank({\mathcal O}) = n$ holds. 
Here, we start with the observer design that is accompanied by a proof of exponential stability, which ensures the states' detectability\index{detectability}. 

\subsubsection*{Simplification of the Model}\label{Simplification}
For the sake of the design and stability proof, we give a simplification on the system \eqref{Ch7:scsys1}-\eqref{Ch7:scsys6}. 
The effect of the salinity profile on the physical parameters is assumed to be sufficiently small so that it can be negligible, i.e.
%\begin{align}\label{nosalt}
$S(x) = 0$.
%\end{align}
Therefore, the heat equation of the sea ice temperature \eqref{Ch7:scsys4} is rewritten as
\begin{align} \label{Ch7:simpsys4} 
\frac{\partial T_{{\rm i}} }{\partial t}(x,t) =& D_{{\rm i}} \frac{\partial^2 T_{{\rm i}} }{\partial x^2}(x,t) + \bar{I}_0 \kappa_{{\rm i}} e^{-\kappa_{{\rm i}} x}, \hspace{1mm} \forall x  \in (0,H(t)),
\end{align}
where the diffusion coefficient is defined as $D_{{\rm i}} = k_0/\rho c_0$. Next, we impose the following assumptions. 
\begin{assum} \label{Ch7:ass:Ht} 
The thickness $H(t)$ is positive and upper bounded, i.e. there exists $\bar{H}>0$ such that $0< H(t) < \bar{H}$, for all $t\geq 0$. 
\end{assum} 
\begin{assum}  \label{Ch7:ass:dotH} 
$\dot H(t)$ is bounded, i.e., there exists $M>0$ such that $| \dot{H}(t) | < M$ , for all $t\geq 0$.  
\end{assum} 

According to \cite{Kwok2009}, the observation data of the sea ice's thickness from the 1950s to 2008 show that the maximum value including the uncertainty is less than 5[m]. Moreover, the largest variation of the thickness in a snow-covered season of a year essentially happens from December to March as an accumulation, and most of the literature shows at most 20 [cm] accumulation per month. Hence, conservatively it is plausible to set $\bar H = 10$ [m], and $M = 50[\textrm{cm/Month}] = 1.9 \times  10^{-7} \textrm{[m/s]}$. 
%A mathematically rigorous proof to ensure Assumptions \ref{ass:Ht} and \ref{ass:dotH} has not been established yet. }
%We formulate the observer structure for sea ice temperature estimation based on the simplified sea ice model composed of \eqref{simpsys4}, \eqref{scsys5}, and \eqref{scsys6}. 

Mathematically, the existence of the classical solution of the simple Stefan problem given by \eqref{Ch7:simpsys4} and \eqref{Ch7:scsys5}--\eqref{Ch7:scsys6} has been established in literature.  
%\cite{Gupta03}, under the conditions of $T_{\rm i}(0,t) \leq T_{\rm m 2}$ for all $t \in (0, \bar t)$ where $\bar t \in (0, \infty]$, $T_{\rm i}(x,0) \leq T_{\rm m2}$ and continuously differentiable for all $x \in [0, H(0)]$. 
We refer the readers to follow \cite{Gupta03} for the detailed explanation. The solution of the original sea ice model \eqref{Ch7:scsys1}--\eqref{Ch7:scsys6} has not been studied due to its high complexity.
% The classical solution is defined by  (i) $\fr{\pa^2 T_{\rm i}}{\pa x^2}$ and $\fr{\pa T_{\rm i}}{\pa t}$ are continuous for $0<x<H(t)$, $0<t<\bar t $; (ii) $T_{\rm i}$ and $\fr{\pa T_{\rm i}}{\pa x}$ are continuous for $0\leq x \leq H(t)$, $0<t<\bar t$; (iii) $T_{\rm i}$ is also continuous for $t=0$, $0<x\leq H(0)$; (iv) $H(t)$ is continuously differentiable for $0\leq t<\bar t$.  }

% a rigorous proof for the existence of the solution of the original sea-ice model \eqref{scsys1}--\eqref{scsys6} has not been established: however, 

\subsection*{Observer Structure}
Suppose that the sea ice thickness and the ice surface temperature are obtained as measurements ${\mathcal Y}_1(t)$ and ${\mathcal Y}_2(t)$, i.e.
\begin{align}\label{Ch7:meas1}
 {\mathcal Y}_{1}(t) =& H(t), \quad {\mathcal Y}_{2}(t) = T_{{\rm i}}(0,t) . 
 \end{align}
The state estimate $\hat{T}_{{\rm i}}$ of the sea ice temperature is governed by a copy of the plant \eqref{Ch7:simpsys4} and \eqref{Ch7:scsys5}-\eqref{Ch7:scsys6} plus the error injection of $H(t)$, namely, as follows: 
\begin{align}
 \label{Ch7:obsys1}\frac{\partial \hat{T}_{{\rm i}}}{\partial t}(x,t) =& D_{{\rm i}} \frac{\partial^2 \hat{T}_{{\rm i}}}{\partial x^2}(x,t) + \bar{I}_0 \kappa_{{\rm i}} e^{-\kappa_{{\rm i}} x} - p_1(x,t)\left({\mathcal Y}_1(t) - \hat H(t) \right), \hspace{2mm} \forall  x \in (0, H(t)) \\
\label{Ch7:obsys2}
 \hat{T}_{{\rm i}}(0,t) =& {\mathcal Y}_{2}(t) - p_2(t) \left({\mathcal Y}_1(t) - \hat H(t) \right), \\
 \label{Ch7:obsys3}\hat{T}_{{\rm i}}(H(t),t) =& T_{{\rm m}2} - p_3(t) \left({\mathcal Y}_1(t) - \hat H(t) \right),  \\
\label{Ch7:obsys4}  \dot{\hat H}(t) = p_4(t) & \left({\mathcal Y}_1(t) - \hat H(t) \right)+  \beta \frac{\partial \hat T_{{\rm i}}}{\partial x}({\mathcal Y}_{1}(t),t)- \frac{F_{{\rm w}}}{q},  
\end{align}
where $\beta := \frac{k_{{\rm i}}}{q}$. 
Next, we define the estimation error states as
\begin{align}
 \tilde T(x,t) := -(T_{{\rm i}}(x,t) - \hat{T}_{{\rm i}}(x,t)), \quad  \tilde H(t) := H(t) - \hat H (t),
\end{align}
where the negative sign is added to be consistent with the description developed in Section \ref{sec:estimator} for the liquid phase. 
Subtraction of the observer system \eqref{Ch7:obsys1}-\eqref{Ch7:obsys4} from the system \eqref{Ch7:simpsys4} and \eqref{Ch7:scsys5}-\eqref{Ch7:scsys6} yields the estimation error system as
\begin{align}
\label{Ch7:errorsys3}
\frac{\partial \tilde{T} }{\partial t}(x,t) = D_{{\rm i}}& \frac{\partial^2 \tilde T}{\partial x^2}(x,t) - p_1(x,t) \tilde H(t),  \hspace{1mm} \forall x \in (0, H(t))\\
\label{Ch7:errorsys4} \tilde{T}(0,t) =& - p_2(t) \tilde H(t),\\
\label{Ch7:errorsys5} \tilde{T}(H(t),t) =& - p_3(t) \tilde H(t),\\
\label{Ch7:errorsys6} \dot{\tilde{H}}(t) =& - p_4(t)\tilde{H}(t) -\beta \frac{\partial \tilde T}{\partial x}(H(t),t) . 
\end{align}
Our goal is to design the observer gains $p_1(x,t)$, $p_2(t)$, $p_3(t)$, $p_4(t)$ so that the temperature error $\tilde{T}$ converges to zero. The main theorem is stated as follows. 
\begin{thm}\label{Ch7:theorem}
Let Assumptions \ref{Ch7:ass:Ht} and \ref{Ch7:ass:dotH} hold. Consider the estimation error system \eqref{Ch7:errorsys3}-\eqref{Ch7:errorsys6} with the design of the observer gains
\begin{align}\label{Ch7:p1gain} 
p_1(x,t) =& \frac{c\lambda x}{\beta} \frac{ I_1 \left( z\right)}{z} + \left( \frac{\ep H(t)}{D_{{\rm i}}} - \frac{3}{\beta} \right) \lambda^2 x   \frac{ I_2 \left( z\right)}{z^2}   + \frac{\lambda^3 x^3}{D_{{\rm i}}\beta}   \frac{ I_3 \left( z\right)}{z^3} , \\
\label{Ch7:p2gain} p_2(t) =& 0, \\
\label{Ch7:p3gain} p_3 (t)=&  -\frac{\lambda}{2 \beta} H(t) - \ep , \\
\label{Ch7:p4gain} p_4(t) =&  c - \frac{\lambda}{2} \left( 1 - \frac{\lambda H(t)^2}{8 D_{{\rm i}}} \right) + \frac{\beta \lambda }{2D_{{\rm i}}} \ep H(t), 
\end{align}
where $\lambda>0$, $c>0$, and $\ep>0$ are positive free parameters, $z$ is defined by  
\begin{align} \label{Ch7:zdef} 
z: = \sqrt{ \bar \lambda (H(t)^2-x^2)}, 
 \end{align} 
where $\bar{\lambda}:= \frac{\lambda}{D_{{\rm i}}}$, and $I_j(\cdot)$ denotes the modified Bessel function of the $j$-th kind. Then, there exist positive constants $c^*>0$ and $\tilde M>0$ such that, for all $c>c^*$, the norm 
\begin{align} \label{Ch7:Phi-def} 
 \Phi(t) : = \int_{0}^{H(t)} \tilde{T}(x,t)^2 dx + \tilde{H}(t)^2 
\end{align} 
satisfies the following exponential decay 
\begin{align} \label{Ch7:Phi-decay} 
{\Phi}(t) \leq \tilde M {\Phi}(0) e^{- \min\{\lambda, c\} t} , 
\end{align} 
namely, the origin of the estimation error system is exponentially stable in the spatial $L_2$ norm.  
\end{thm}

%The technical novelty of the proposed observer is that the design only requires the measured surface temperature and the thickness, which is practically implementable. On the other hand, the estimation design developed in \cite{Shumon17journal,Shumon17seaice} additionally requires the temperature gradient at the liquid-solid (i.e., ice-ocean) interface, which is nearly impossible to measure.

\begin{remark} \label{Ch7:rem:1} 
The observer gains \eqref{Ch7:p1gain}-\eqref{Ch7:p4gain} include the thickness $H(t)$, so the gains are not precomputed offline, but are easily calculated \emph{online}, along with the state estimation. Owing to the slow dynamics of the sea ice model, the computation time is much less than the time step size, which enables the real-time computation of the proposed observer.  
\end{remark}

\begin{remark} \label{Ch7:rem:2} 
The measurements \eqref{Ch7:meas1} are assumed to be noiseless; however, in practice, the measured data accompany with some noise. 
%Though it is possible to implement the proposed observer with the noisy measurements as itself, 
Preferably the observer needs pre-filtering to deal with the noisy measurements.   
\end{remark} 

%The measurements are also assumed to be obtained in continuous time. 
To handle the discrete-time measurements in practice as in \cite{petrus17}, the designed observer should be discretized in time such as Euler or Runge-Kutta methods\index{Runge-Kutta method} so that the estimation can be computed at every sampling of the discrete-time measurements.  
The free parameters $\lambda$, $c$, and $\ep$ have their physical units [1/s], [1/s], and [$^\circ$C/m], respectively. Hence we can see the consistency of the physical units in the estimation error system \eqref{Ch7:errorsys3}--\eqref{Ch7:errorsys6} together with \eqref{Ch7:p1gain}--\eqref{Ch7:p4gain}. 
%For instance, the gain $p_{3}(t)$ in \eqref{Ch7:p3gain} leads to a unit [$^\circ$C/m] which matches with the unit in \eqref{Ch7:errorsys5}.

%The estimation of the snow temperature profile is also achievable using the same observer structure as in Theorem \ref{Ch7:theorem} with measurements of the snow thickness and snow surface temperature. To avoid the lengthy statement, we do not provide it in this paper.

\subsection*{Gain Derivation via State Transformation} 
For the estimation error system \eqref{Ch7:errorsys3}--\eqref{Ch7:errorsys6}, we apply the following invertible transformations: 
\begin{align}\label{Ch7:trs1}
\tilde T(x,t) =&  w(x,t) - \int_{x}^{H(t)} q(x,y) w(y,t) dy - \psi (x,H(t)) \tilde H(t), \\
\label{Ch7:trs2}w(x,t) = & \tilde{T}(x,t) - \int_{x}^{H(t)} r(x,y) \tilde{T}(y,t) dy - \phi (x, H(t)) \tilde H(t), 
\end{align}
which map the estimation error system \eqref{Ch7:errorsys3}-\eqref{Ch7:errorsys6} into the following target system\index{target system}: 
\begin{align}\label{Ch7:tarsys4}
 w_{t}(x,t) =& D_{{\rm i}} w_{xx}(x,t) - \lambda w(x,t) - \dot{H}(t)f(x,H(t))\tilde H(t), \hspace{2mm} \forall x \in (0,H(t))\\
\label{Ch7:tarsys5}w(0,t) =& 0,\\
\label{Ch7:tarsys6}w(H(t),t) =&\ep  \tilde {H}(t),\\
\label{Ch7:tarsys7}\dot{\tilde{H}}(t) =& - c\tilde{H}(t) -\beta w_{x}(H(t),t) ,
\end{align}
where $f(x,H(t))$ is to be determined. Taking the first and second spatial derivatives of the transformation \eqref{Ch7:trs1}, we get  
\begin{align}\label{Ch7:trs1x}
\tilde T_x(x,t) =&  w_x(x,t) + q(x,x) w(x,t)  \notag\\
&- \int_{x}^{H(t)} q_x(x,y) w(y,t) dy - \psi_x (x,H(t)) \tilde H(t), \\
\label{Ch7:trs1xx}
\tilde T_{xx}(x,t) =&  w_{xx}(x,t) + q(x,x) w_x(x,t)+\left(q_x(x,x) + \frac{d}{dx} q(x,x) \right) w(x,t)   \notag\\
&- \int_{x}^{H(t)} q_{xx}(x,y) w(y,t) dy- \psi_{xx}(x,H(t)) \tilde H(t). 
\end{align} 
Next, taking the time derivative of \eqref{Ch7:trs1} along the solution of the target system \eqref{Ch7:tarsys4}--\eqref{Ch7:tarsys7}, using integration by parts, and substituting the boundary condition \eqref{Ch7:tarsys6}, we get 
%we get 
%\begin{align} 
%\tilde T_t(x,t) 
%=& D_{{\rm i}} w_{xx}(x,t) - \lambda w(x,t) - \dot{H}(t)f(x,H(t))\tilde H(t) \notag\\
%&- D_{{\rm i}} \int_{x}^{H(t)} q(x,y) w_{yy}(y,t) dy \notag\\
%& + \lambda \int_{x}^{H(t)} q(x,y) w(y,t) dy \notag\\
%&+ \dot{H}(t) \int_{x}^{H(t)} q(x,y) f(y,H(t)) dy \tilde H(t) \notag\\
%&- \dot{H}(t)  q(x,H(t)) w(H(t),t) \notag\\
%& + \psi (x,H(t)) ( c\tilde{H}(t) +\beta  w_{x}(H(t),t)) \notag\\
%&- \dot{H}(t) \psi_{H}(x,H(t)) \tilde H(t). 
%\end{align} 
%Taking integration by parts and substituting the boundary condition \eqref{Ch7:tarsys6}, we get 
\begin{align} 
\tilde T_t(x,t) 
 =& D_{{\rm i}} w_{xx}(x,t) + D_{{\rm i}} q(x,x)w_{x}(x,t)  - (\lambda + D_{{\rm i}} q_{y}(x,x) )  w(x,t)  \notag\\
& + ( \beta \psi (x,H(t))- D_{{\rm i}}  q(x,H(t))) w_{x}(H(t),t) \notag\\
& + ( D_{{\rm i}} \ep q_{y}(x,H(t)) + c \psi (x,H(t)) ) \tilde{H}(t)  \notag\\
&+  \int_{x}^{H(t)} (\lambda q(x,y) - D_{{\rm i}} q_{yy}(x,y) ) w(y,t) dy \notag\\
&-  \dot{H}(t) \tilde H(t) \left( \ep q(x,H(t)) + \psi_{H}(x,H(t)) \right. \notag\\
&\hspace{-0mm} \left. + f(x,H(t)) - \int_{x}^{H(t)} q(x,y) f(y,H(t)) dy \right)  . \label{Ch7:trs1t}
\end{align} 
Thus, by \eqref{Ch7:trs1xx} and \eqref{Ch7:trs1t}, we have 
\begin{align} 
&\tilde T_t(x,t) - D_{{\rm i}} \tilde T_{xx}(x,t) + p_1(x,t) \tilde H(t) \notag\\
=&  - \left(\lambda + 2D_{{\rm i}} \frac{d}{dx}q(x,x) \right) w(x,t)  \notag\\
& + ( \beta \psi (x,H(t))- D_{{\rm i}}  q(x,H(t))) w_{x}(H(t),t) \notag\\
& + \left( D_{{\rm i}} \ep q_{y}(x,H(t)) + D_{{\rm i}} \psi_{xx}(x,H(t))  + c \psi (x,H(t)) + p_1(x,t)\right) \tilde{H}(t)  \notag\\
&+  \int_{x}^{H(t)} (\lambda q(x,y) + D_{{\rm i}} q_{xx}(x,y) - D_{{\rm i}} q_{yy}(x,y) ) w(y,t) dy \notag\\
&-  \dot{H}(t) \tilde H(t) \left( \ep q(x,H(t)) + \psi_{H}(x,H(t)) \right. \notag\\
& \left. + f(x,H(t)) - \int_{x}^{H(t)} q(x,y) f(y,H(t)) dy \right)  . 
\end{align} 
Substituting $x=0$ and $x = H(t)$ into \eqref{Ch7:trs1}, we get 
\begin{align} 
\tilde T(0,t) + p_2(t) \tilde H(t) =& - \int_{0}^{H(t)} q(0,y) w(y,t) dy\notag\\
 &+ ( p_2(t) - \psi (0,H(t)) ) \tilde H(t), \\
\tilde T(H(t),t)  + p_3(t) \tilde H(t) =&(\ep  - \psi (H,H) + p_3(t)) \tilde H(t). 
\end{align} 
Moreover, substituting $x = H(t)$ into \eqref{Ch7:trs1x} yields 
\begin{align}
&\dot{\tilde{H}}(t) + p_4(t) \tilde H(t) + \beta 	\tilde T_x(H(t),t) \notag\\
=& ( p_4(t) - c+ \beta (\ep q(H(t),H(t)) - \psi_x (H(t),H(t)))  ) \tilde H(t). 
\end{align}
Therefore, for the equations \eqref{Ch7:errorsys3}--\eqref{Ch7:errorsys6} to hold, the gain kernel functions must satisfy the following conditions: 
\begin{align} \label{Ch7:qeq1} 
	q_{xx}(x,y) - q_{yy}(x,y) = & - \bar \lambda q(x,y) , \\
\label{Ch7:qeq2} \frac{d}{dx} q(x,x) =& - \frac{\bar \lambda}{2} , \quad q(0,y) = 0, \\
\label{Ch7:psieq1} \beta \psi (x,H(t)) =& D_{{\rm i}} q(x,H(t)), 
\end{align} 
and the observer gains must satisfy 
\begin{align} 
\label{Ch7:p1cond}  p_1(x,t) = & - D_{{\rm i}}( \ep q_{y}(x,H(t))+ \psi_{xx}(x,H)) - c \psi (x,H), \\
 \label{Ch7:p2cond}  p_2(t) = & \psi (0,H(t)), \\
 \label{Ch7:p3cond}   p_3(t)  =& \psi (H(t),H(t)) - \ep, \\
 \label{Ch7:p4cond}   p_4(t) = & c - \beta (\ep q(H(t),H(t)) - \psi_x (H(t),H(t))), 
\end{align}
and the function $f(x,H(t))$ must satisfy 
\begin{align} \label{Ch7:fcond} 
 &f(x,H) +\ep q(x,H) + \psi_{H}(x,H)= \int_{x}^{H} q(x,y) f(y,H) dy . 
\end{align}
The solutions to \eqref{Ch7:qeq1}--\eqref{Ch7:psieq1} are uniquely given by 
\begin{align} \label{Ch7:qxy} 
q(x,y) =& - \bar{\lambda} x \frac{ I_1 \left( \sqrt{\bar{\lambda} (y^2-x^2)} \right)}{\sqrt{\bar{\lambda} (y^2-x^2)}}, \\
\psi(x,H(t)) =&  - \frac{\lambda}{\beta} x \frac{ I_1 \left( z \right)}{z}, \label{Ch7:psixy}
\end{align}
where $z$ is defined by \eqref{Ch7:zdef}. Then, 
%by taking the derivatives of the obtained gain kernels, we get 
%\begin{align} \label{Ch7:psixx} 
%\psi_{xx}(x,H(t)) =&  \frac{\lambda^2}{D_{{\rm i}}\beta} x \left( 3 \frac{ I_2 \left( z\right)}{z^2} -\bar{\lambda} x^2 \frac{ I_3 \left( z\right)}{z^3} \right), \\
%q_y(x,H(t)) = &   - \frac{\lambda^2}{D_{{\rm i}}^2}H(t) x \frac{ I_2 \left( z \right)}{ z^2} , \label{Ch7:qyH} 
%\end{align}
%where $z$ is defined by \eqref{Ch7:zdef}. 
using \eqref{Ch7:qxy}--\eqref{Ch7:psixy}, the conditions \eqref{Ch7:p1cond}--\eqref{Ch7:p4cond} are led to the explicit formulations of the observer gains given as  \eqref{Ch7:p1gain}--\eqref{Ch7:p4gain}. In the similar manner, the conditions for the gain kernel functions of the inverse transformation \eqref{Ch7:trs2} are given by 
\begin{align}
\label{Ch7:req1}  r_{xx}(x,y)- r_{yy}(x,y) =& \bar{\lambda} r(x,y), \\
\label{Ch7:req2}  \frac{d}{dx} r(x,x) =& \fr{\bar \lambda}{ 2} , \quad r(0,y)=0, \\
 \label{Ch7:phieq1} \beta \phi (x,H(t))  =& D_{{\rm i}}  r(x,H(t)), 
\end{align}
and, the function $f(x,H(t))$ is obtained by 
\begin{align} \label{Ch7:fsol} 
	f(x,H(t)) =  r(x,H(t)) p_3(H(t)) + \phi_{H}(x,H(t)). 
\end{align}
The solutions to \eqref{Ch7:req1}--\eqref{Ch7:phieq1} are given by 
\begin{align}
\label{Ch7:rsol}	r(x,y) =& \bar{\lambda} x \frac{ J_1 \left( \sqrt{\bar{\lambda} (y^2-x^2)} \right)}{\sqrt{\bar{\lambda} (y^2-x^2)}}, \hspace{1mm}\phi(x,H) = \frac{\lambda}{\beta} x \frac{ J_1 \left( z \right)}{z},  
\end{align}
where $J_1$ is Bessel function of the first kind. Using the solutions \eqref{Ch7:rsol}, the function $f(x,H(t))$ is obtained explicitly by \eqref{Ch7:fsol}, which also satisfies the condition \eqref{Ch7:fcond}. Hence, the transformation from $(\tilde T, \tilde H)$ to $(w, \tilde H)$ is invertible.

\subsection*{Stability Analysis}\label{Ch7:sec:Stability}
We prove the exponential stability of the origin of the estimation error system \eqref{Ch7:errorsys3}-\eqref{Ch7:errorsys6} in the spatial $L_2$ norm. First, we show the exponential stability of the origin of the target system \eqref{Ch7:tarsys4}-\eqref{Ch7:tarsys7}. We consider the following Lyapunov functional 
\begin{align}\label{Ch7:V1}
V = \frac{1}{2} ||w||^2 + \frac{\ep}{2 \beta} \tilde{H}(t)^2.  
\end{align}
Taking the time derivative of \eqref{Ch7:V1} together with the solution of \eqref{Ch7:tarsys4}-\eqref{Ch7:tarsys7} yields 
\begin{align}
\dot{V} =&- D_{{\rm i}} || w_x ||^2 - \lambda || w ||^2 - \frac{\ep c}{\beta} \tilde H(t)^2 +\frac{\dot{H}(t)}{2}\ep^2 \tilde{H}(t)^2 \notag\\
 &\hspace{-0mm} - \dot{H}(t)\tilde H(t)  \int_{0}^{H(t)} w(x,t) f(x,H(t))dx. \label{Ch7:V1dot} 
\end{align} 
Applying Young's and Cauchy-Schwarz inequalities to the last term in \eqref{Ch7:V1dot} with the help of Assumption \ref{Ch7:ass:dotH}, 
%we have 
%\begin{align} 
%&- \dot{H}(t)\tilde H(t)  \int_{0}^{H(t)} w(x,t) f(x,H(t))dx \notag\\
%\leq & M |\tilde H(t)| \cdot  \left( \int_{0}^{H(t)} f(x,H(t))^2 dx \right)^{1/2} \cdot ||w || \notag\\
%\label{Ch7:CSineq}  \leq &  \frac{ M^2  \bar f }{2 \lambda} \tilde H(t)^2  + \frac{\lambda}{2} ||w || ^2 , 
%\end{align} 
%where $\bar f :=\max_{H(t) \in (0, \bar H)} \int_{0}^{ H(t)}  f(x,H(t))^2 dx $. Applying \eqref{Ch7:CSineq} to \eqref{Ch7:V1dot}, 
%we obtain the following inequality: 
%\begin{align} 
%\dot{V}_1 \leq & - D_{{\rm i}} || w_x ||^2 - \frac{\lambda}{2} || w ||^2  + \ep  \tilde H(t) w_{x}(H(t),t) \notag\\
%& + \left(\frac{ M^2 \bar f }{2 \lambda} + \frac{M \ep^2}{2} \right) \tilde H(t)^2,  \label{Ch7:V1dotineq} 
%\end{align}
%where $\bar f :=\max_{H \in (0, \bar H)} \int_{0}^{ H}  f(x,H)^2 dx $. Next, we consider 
%\begin{align} \label{Ch7:V2}
%V_{2} = \frac{1}{2} \tilde{H}(t)^2. 
%\end{align} 
%Taking the time derivative of \eqref{Ch7:V2} along with \eqref{Ch7:tarsys7} yields 
%\begin{align} 
%\dot{V}_{2} = &- c \tilde H(t)^2 - \beta \tilde H(t)  w_{x}(H(t),t)  . \label{Ch7:V2dot} 
%\end{align} 
%Let $V$ be the Lyapunov functional defined by 
%\begin{align} \label{Ch7:Vdef} 
%V = V_1 + \frac{\ep}{\beta} V_2 . 
%\end{align} 
%Combining \eqref{Ch7:V1dotineq} and \eqref{Ch7:V2dot}, 
and choosing the gain parameter $c$ to satisfy 
\begin{align}  \label{Ch7:c:cond} 
 c >  \frac{\beta  M^2 \bar f}{\ep \lambda}  + \beta M  \ep, 
 \end{align}
%  the time derivative of \eqref{Ch7:Vdef} is shown to satisfy the following inequality: 
one can obtain the following inequality: 
%\begin{align}
%\dot{V} \leq &- D_{{\rm i}} || w_x ||^2 - \frac{\lambda}{2} || w ||^2  \notag\\
%&- \left( \frac{\ep}{\beta} c  -  \frac{ M^2 \bar f }{2 \lambda} - \frac{M \ep^2}{2}\right) \tilde H(t)^2 .  \label{Ch7:Vdot1} 
%\end{align} 
%Hence, by 
 %the inequality \eqref{Ch7:Vdot1} is led to 
 \begin{align} 
 \dot{V}  \leq& - \min\{\lambda, c\}  V . \label{Ch7:Vdotineq} 
\end{align} 
Applying comparison principle to the differential inequality \eqref{Ch7:Vdotineq}, we get 
\begin{align}\label{Ch7:Vfin}
&V(t) \leq V(0) e^{- \min\{\lambda, c\} t }.
\end{align}
Hence, the target system \eqref{Ch7:tarsys4}-\eqref{Ch7:tarsys7} is exponentially stable at the origin. Due to the invertibility of the transformations \eqref{Ch7:trs1} and \eqref{Ch7:trs2}, there exist positive constants $\underline{M}>0$ and $\bar M>0$ such that for the norm $\Phi(t)$ defined in \eqref{Ch7:Phi-def} the inequalities hold $\underline{M} \Phi(t)  \leq V(t) \leq \bar M  \Phi(t)$. Hence, we obtain \eqref{Ch7:Phi-decay} by defining $\tilde M = \bar M/\underline{M}$, which completes the proof of Theorem \ref{Ch7:theorem}. Note that the designed backstepping observer achieves faster convergence with a possibility of causing overshoot since the overshoot coefficient $\bar M/ \underline{M}$ is a monotonically increasing function in the observer gains' parameters $(\lambda, c)$.

While we have focused on the simplified PDE \eqref{Ch7:simpsys4} to derive a rigorous proof of the proposed state estimation design \eqref{Ch7:obsys1}-\eqref{Ch7:obsys4} with observer gains given by \eqref{Ch7:p1gain}--\eqref{Ch7:p4gain}, simulation studies are performed by applying  the estimation design to the original thermodynamic model \eqref{Ch7:scsys1}-\eqref{Ch7:scsys6} including salinity. 
%\subsection{Simulation Results of State Estimation}

  \subsection{Numerical Tests of the Sea Ice Estimation} 
\subsubsection*{Initial conditions} 
The simulation results of temperature estimation $\hat{T}_{{\rm i}}$ computed by \eqref{Ch7:obsys1}-\eqref{Ch7:obsys4} along with the available measurements obtained by the online calculation of \eqref{Ch7:scsys1}-\eqref{Ch7:scsys6} are shown in Fig. \ref{Ch7:fig:2}. 
%As stated in Remark \ref{Ch7:rem:2}, we suppose that the noiseless measurements are available for the observer through pre-filtering in practice. 
Here the initial temperature profiles are formulated as 
\begin{align} 
T_{{\rm s}}(x,0) =& \frac{k_{0}(T_{{\rm m}1}-T_{0})}{k_{s} H_{0}} x + T _{0}, \\
T_{{\rm i}}(x,0) =& \frac{T_{{\rm m}1} - T_{0}}{H_{0}} x + T_{0} + a \sin \left( \frac{4 \pi x}{H_0} \right),
\end{align}  
where $T_0 = T_{{\rm i}}(0,0)$ which is obtained by solving fourth order algebraic equation from \eqref{Ch7:scsys1} and the input data, and $a$ is set as $a = 1$ [C$^\circ$]. The estimated initial temperature is chosen as 
\begin{align} 
\hat T_{{\rm i}}(x,0) = \frac{T_{{\rm m}1} - T_{0}}{H_{0}^2 (1-2d)} (x^2 - 2dH_0 x) + T_{0}
\end{align} 
 with setting $d=1/4$. Hence, the initial temperature estimate is lower than the actual temperature. This initial condition satisfies the boundary conditions \eqref{Ch7:obsys2} and \eqref{Ch7:obsys3}. The initial state of the estimated ice thickness $\hat H(0)$ is set as that of the true thickness, i.e., $\hat H(0) = H(0)$, which is feasible because the thickness is actually measured. 
 
 \subsubsection*{Tuning method for gain parameters} 
 The design parameters $(\lambda,c,\ep)$ are selected as follows:\\ 
(i) Choose $\ep \approx \beta $ for the norm \eqref{Ch7:V1} to be similarly weighted. \\
(ii) Select $\lambda$ to be the inverse of a desired time constant (i.e., the time at 63\% decay of the estimation error is achieved): here we set as one day, leading to $\lambda \approx \frac{1}{24 \times 3600} = 1.2 \cdot 10^{-5}$.\\
(iii) Select $c$ sufficiently larger than $\lambda$ so that the decay rate $\min\{\lambda, c\}$ is not reduced and \eqref{Ch7:c:cond} is satisfied. \\
Finally, these parameters are varied around these reference values until we observe a smooth and sufficiently fast convergence. Throughout the simulation, we see that the minimum value of the time step size in ode solver is more than 1 minute, while the computation time of each time update is less than 0.1 seconds, which shows its real-time implementability as addressed in Remark \ref{Ch7:rem:1}. 

%\subsection*{Simulation Results} 

\subsubsection*{Numerical simulation of state estimation} 
\begin{figure*}[htb!]
\begin{center}
\subfloat[Open-loop estimation, i.e., $p_1(x,t) = 0$ and $p_i(t) = 0$ for $i=2,3,4$.  ]{\includegraphics[width=0.8\linewidth]{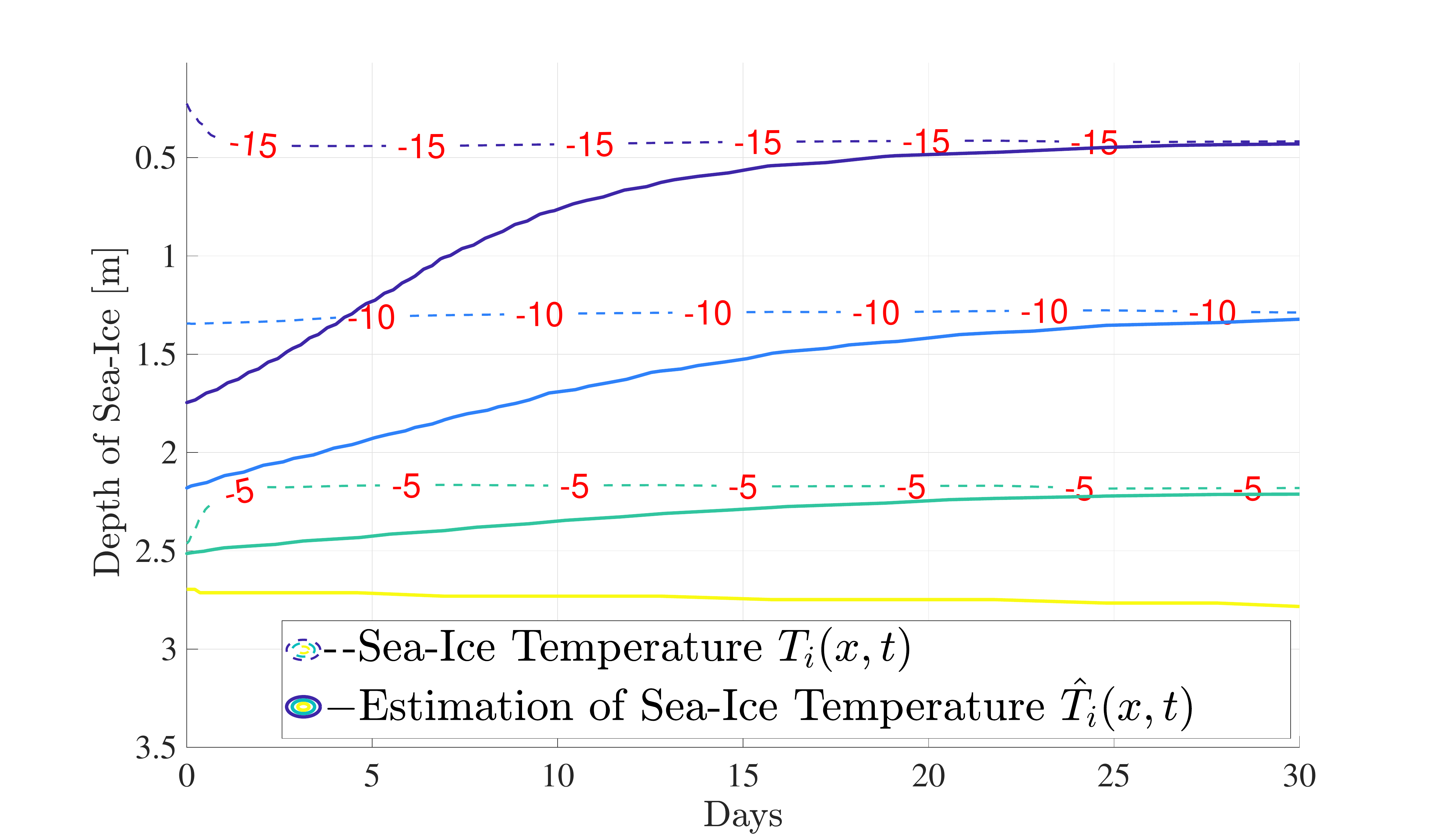}}\\
\subfloat[The proposed estimation with the observer gains given in \eqref{Ch7:p1gain}--\eqref{Ch7:p4gain}. ]{\includegraphics[width=0.8\linewidth]{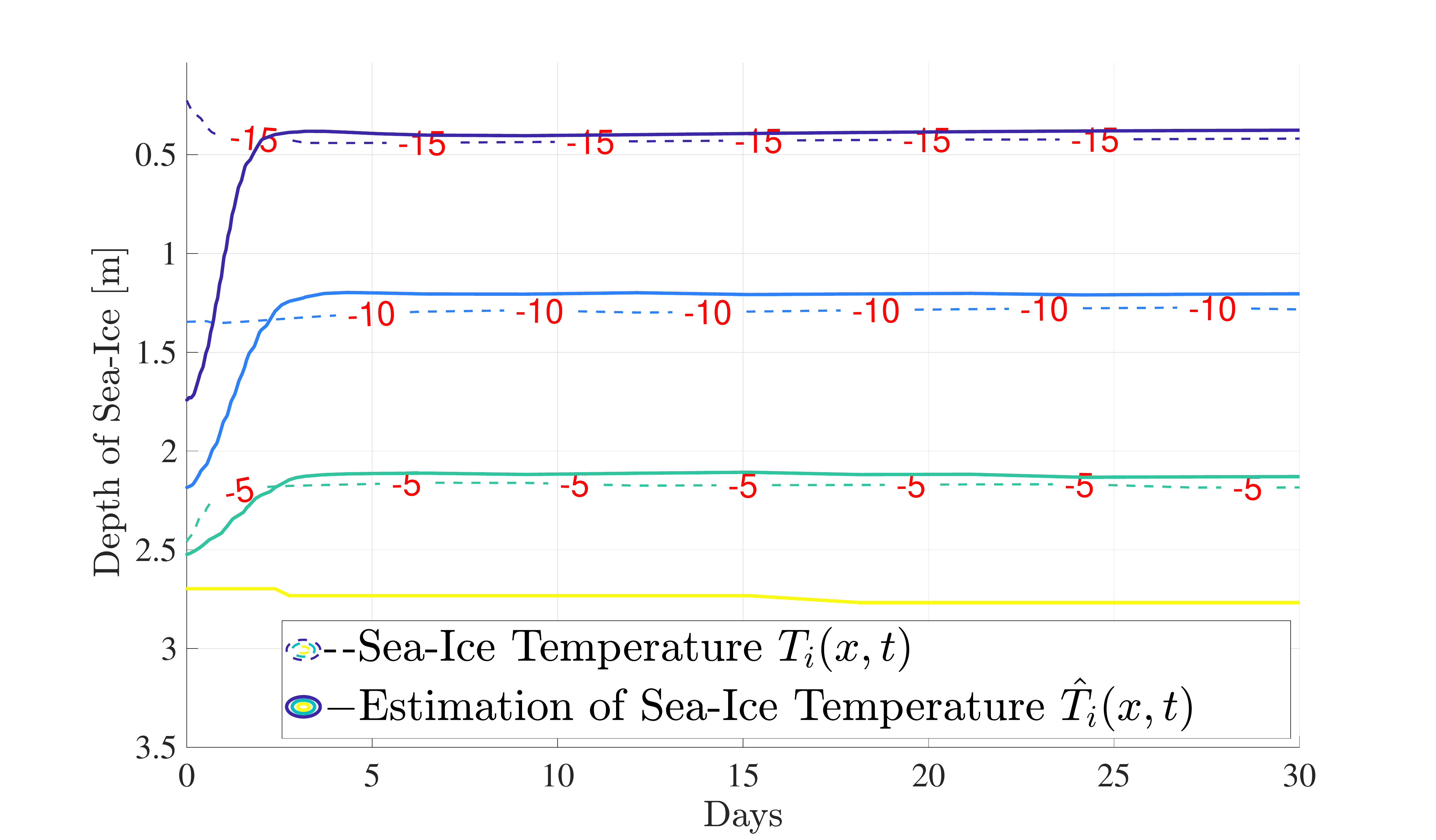}}
\caption{Simulation results of the plant \eqref{Ch7:scsys1}--\eqref{Ch7:scsys6} and the estimator \eqref{Ch7:obsys1}-\eqref{Ch7:obsys4} using parameters in Table \ref{Ch7:table:1}. The designed backstepping observer achieves faster convergence to the actual state than the straightforward open-loop estimation\index{open-loop estimation}.}
\label{Ch7:fig:2}
\end{center}
\end{figure*}
\begin{figure*}[htb!]
\begin{center}
\subfloat[The proposed estimation with larger value of $\lambda$ than Fig. \ref{Ch7:fig:2} (b). The overshoot beyond the true temperature is observed during the first two days. ]{\includegraphics[width=0.8\linewidth]{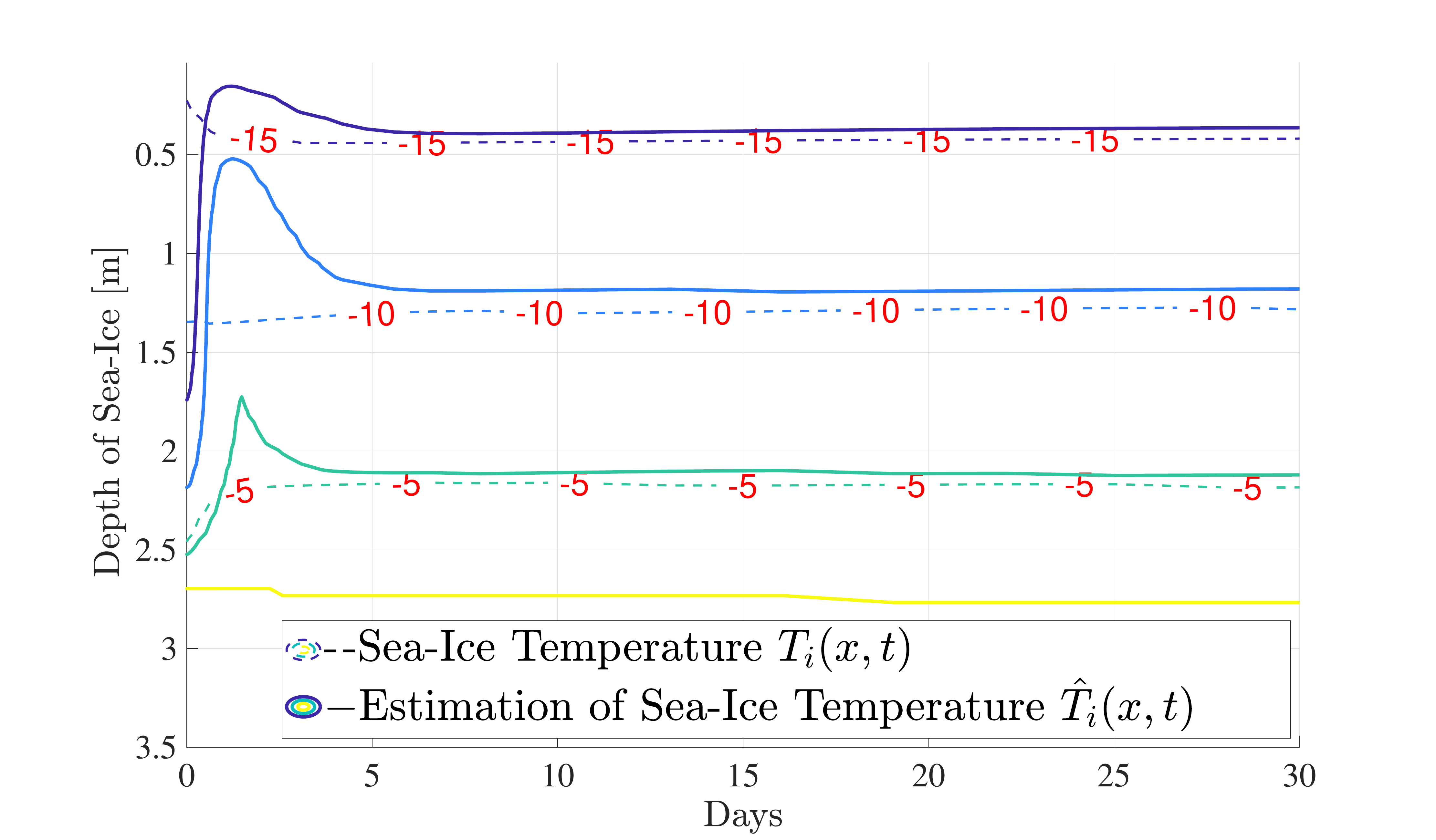}} \\
\subfloat[The proposed estimation with smaller value of $\lambda$ than Fig. \ref{Ch7:fig:2} (b). The convergence speed gets slower than the result of Fig. \ref{Ch7:fig:2} (b). ]{\includegraphics[width=0.8\linewidth]{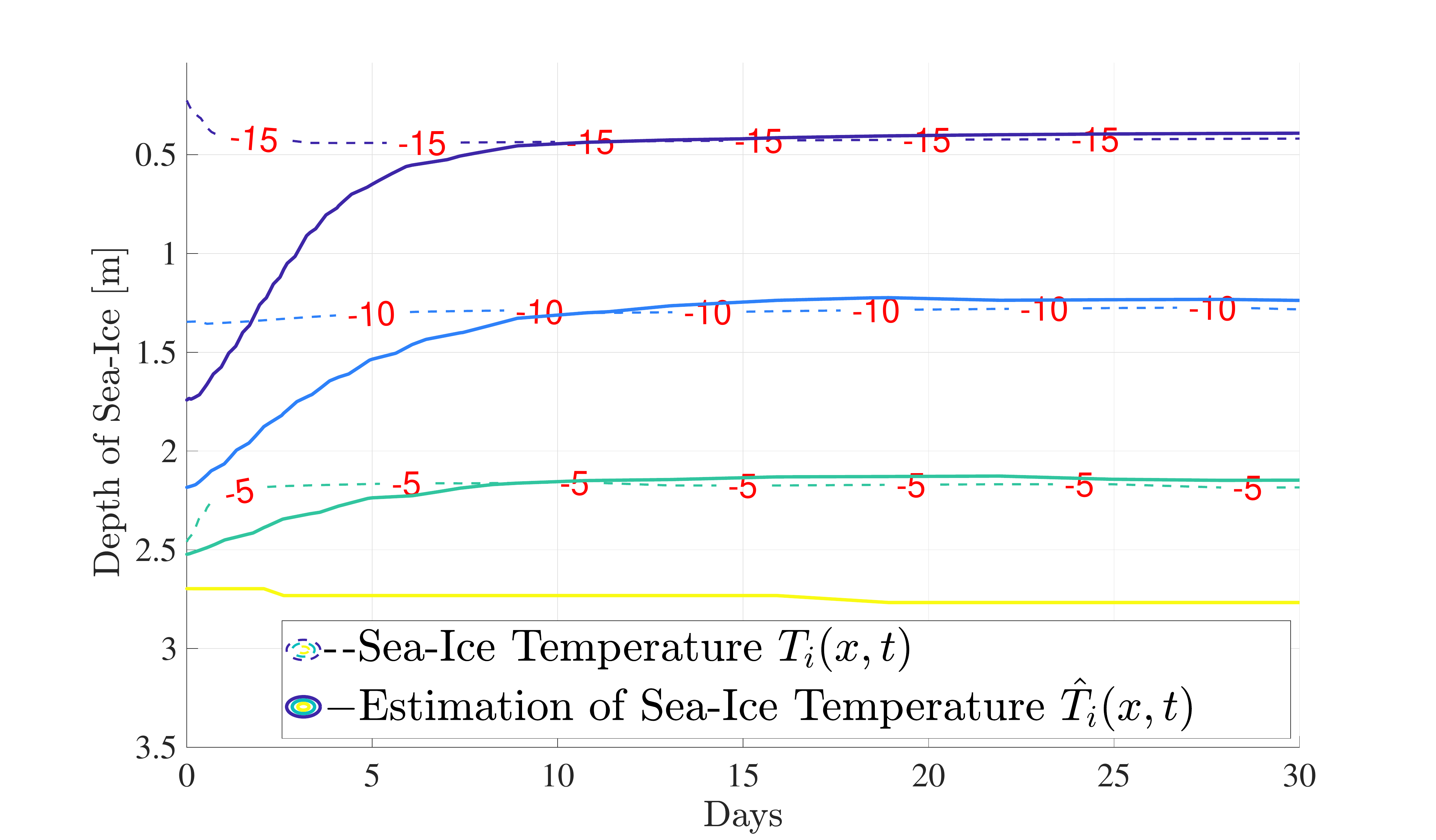}}
\caption{Simulation results of the plant \eqref{Ch7:scsys1}--\eqref{Ch7:scsys6} and the bacsktepping estimator \eqref{Ch7:obsys1}-\eqref{Ch7:obsys4} with some chosen free parameters.}
\label{Ch7:fig:2-2}
\end{center}
\end{figure*}

\begin{figure*}[htb!]
\begin{center}
\subfloat[Temperature profile of both true and estimate on January 1st. ]{\includegraphics[width=0.5\linewidth]{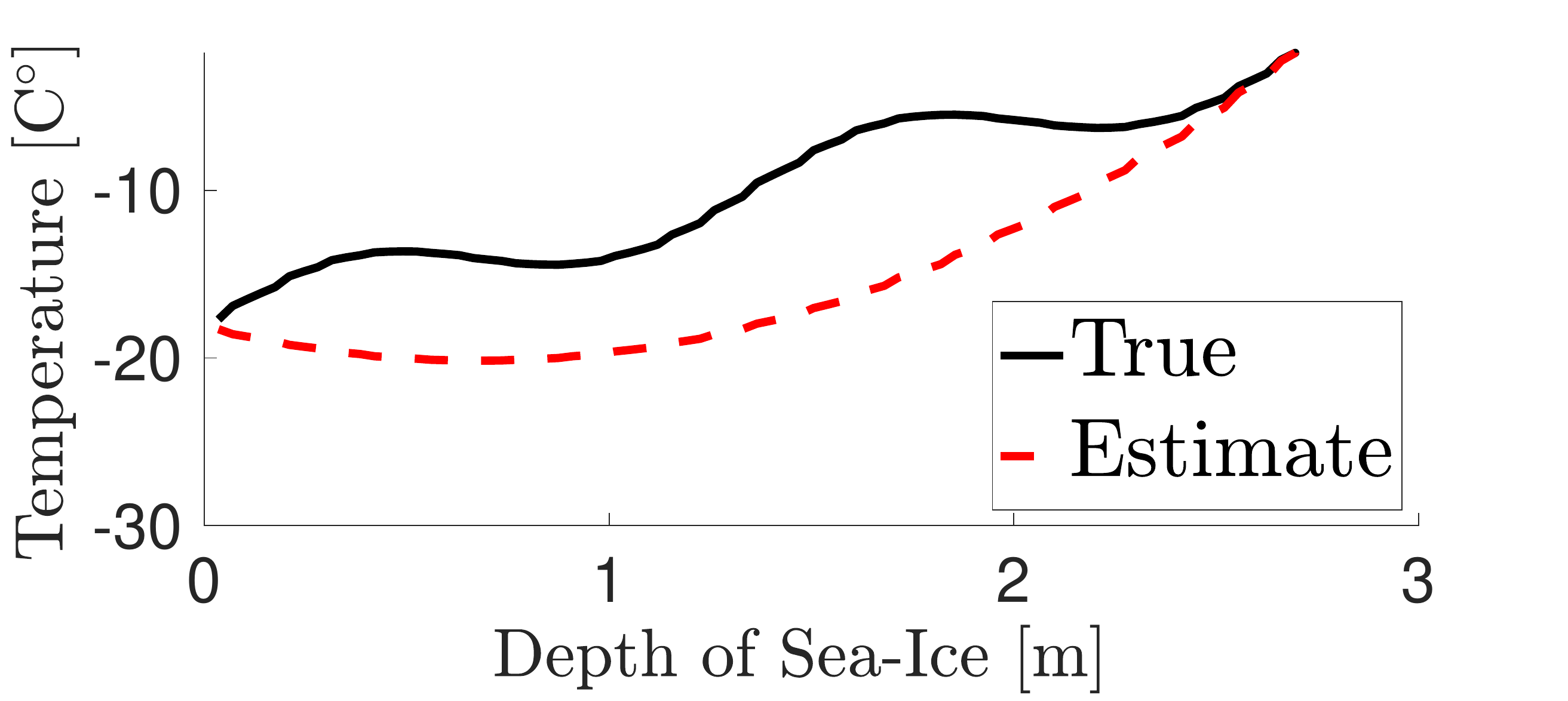}}
\subfloat[Temperature profile of both true and estimate on January 2nd. ]{\includegraphics[width=0.5\linewidth]{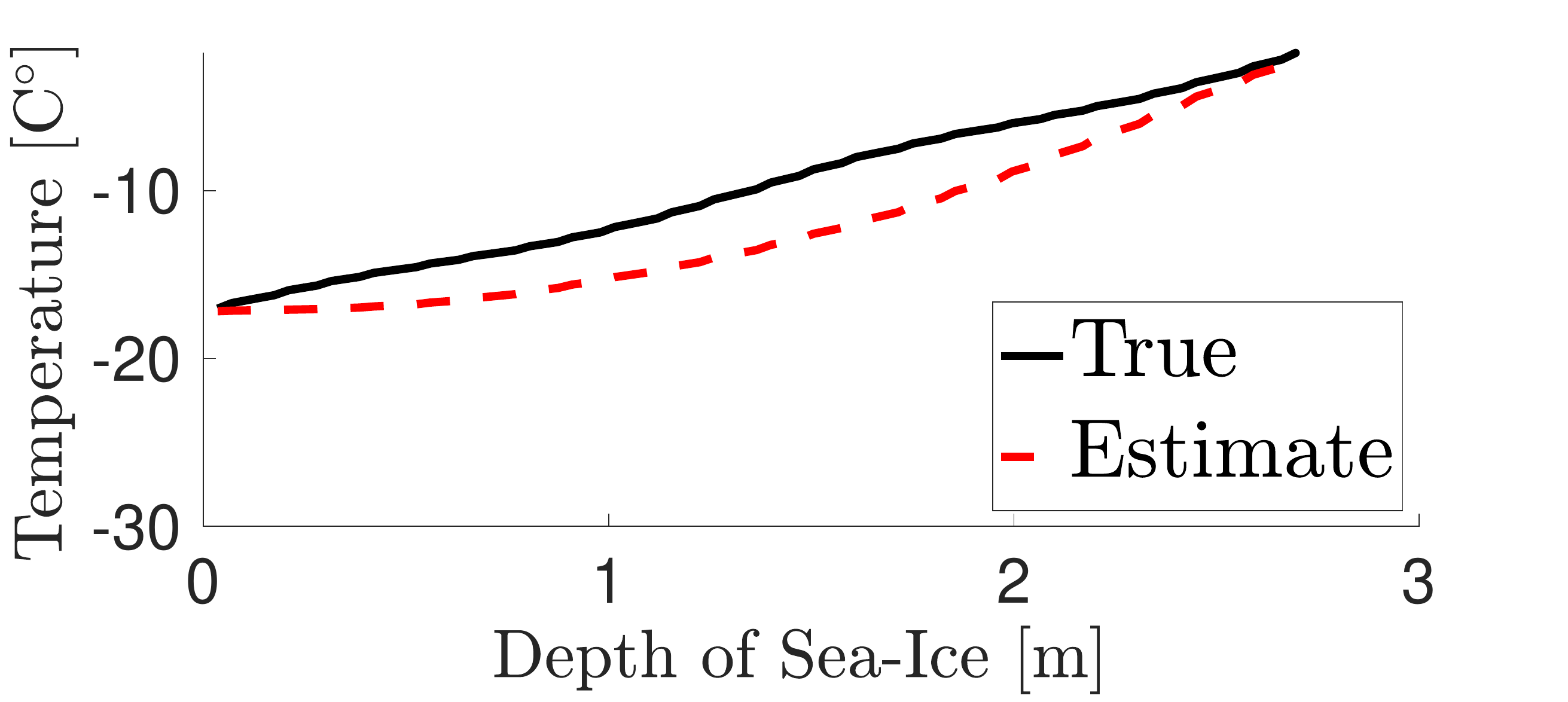}}\\
\subfloat[Temperature profile of both true and estimate on January 3rd. ]{\includegraphics[width=0.5\linewidth]{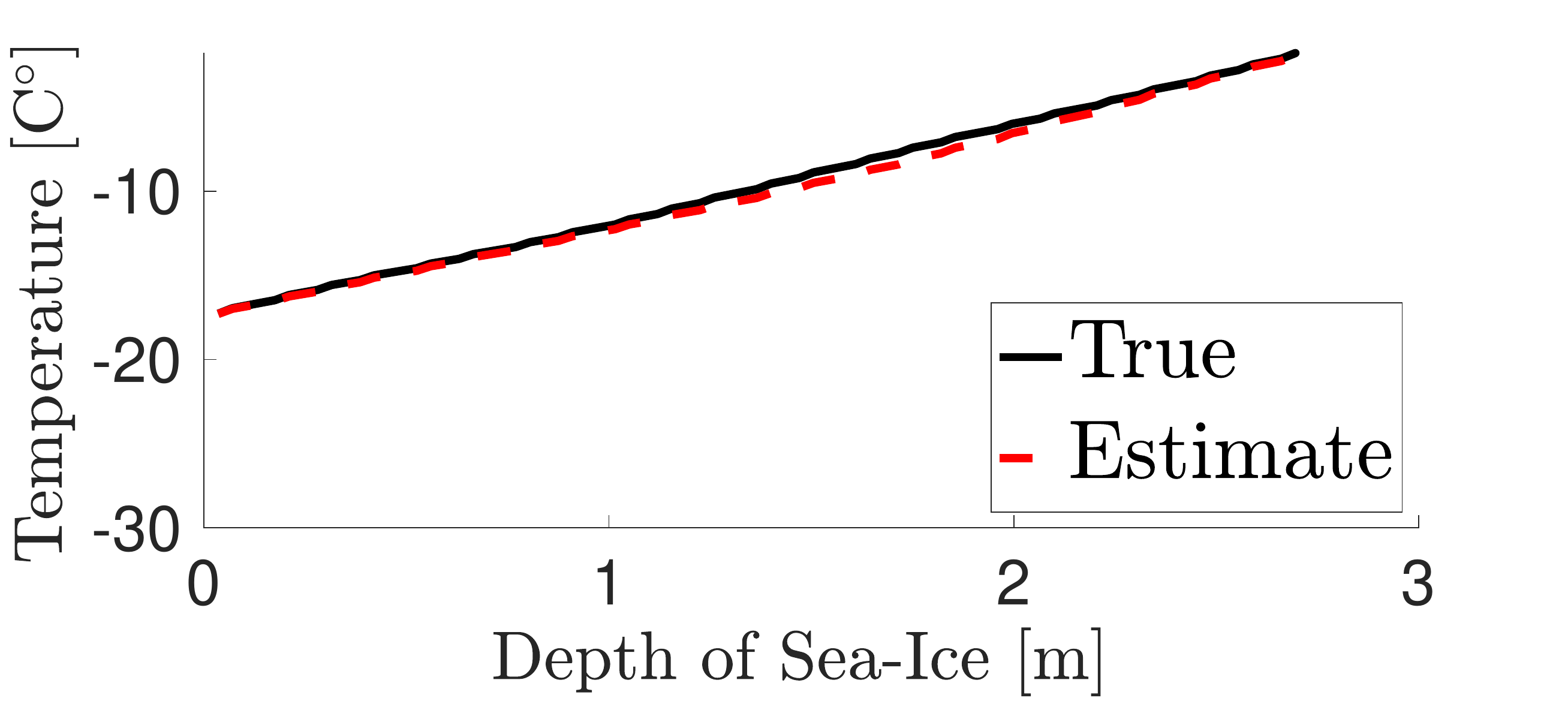}}
\subfloat[The time evolution of thickness estimation error $\tilde H(t)$. ]{\includegraphics[width=0.5\linewidth]{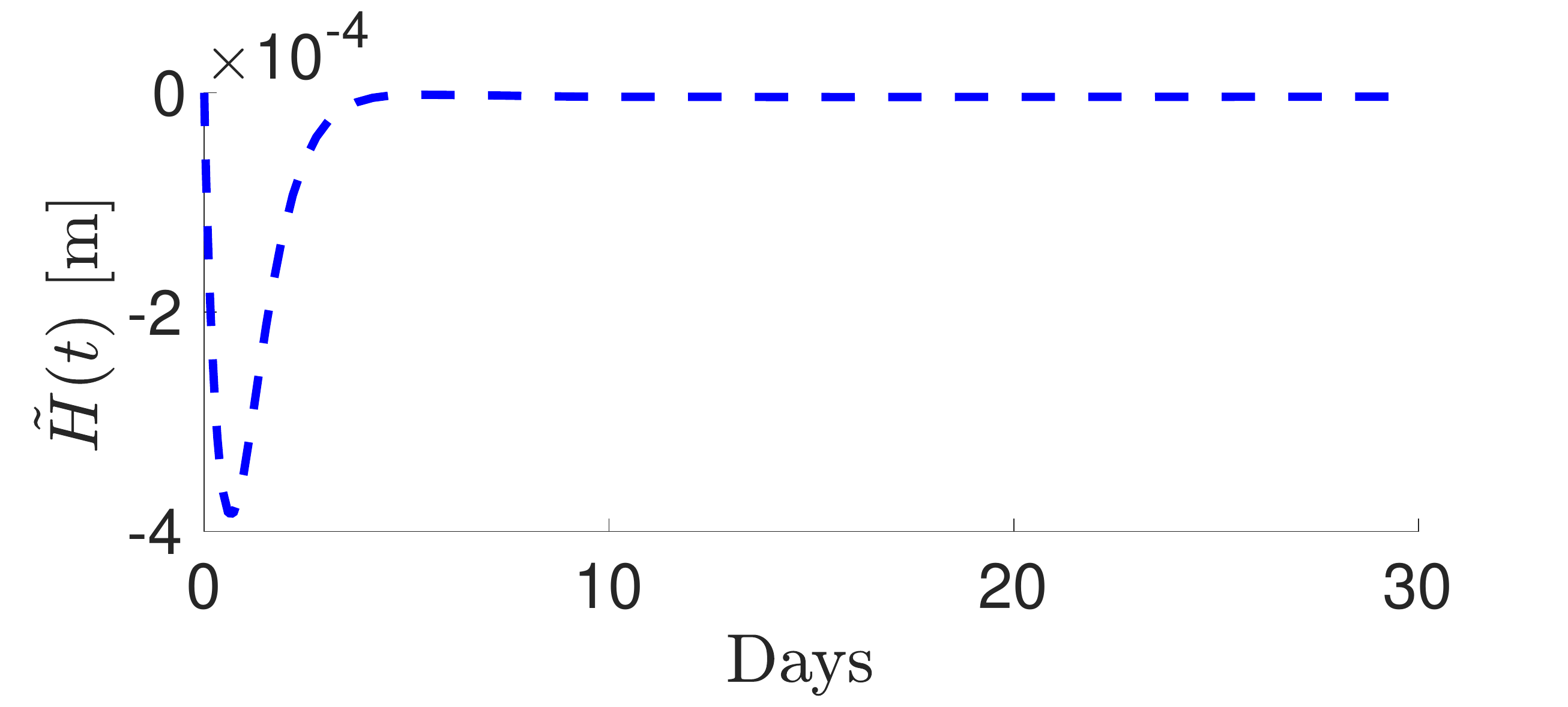}}
\caption{Simulation result of the plant \eqref{Ch7:scsys1}--\eqref{Ch7:scsys6} and the estimator \eqref{Ch7:obsys1}-\eqref{Ch7:obsys4} with parameters used in Fig. \ref{Ch7:fig:2} (b). }
\label{Ch7:fig:3}
\end{center}
\end{figure*}

         The contour plot of the simulation results of $T_{{\rm i}}(x,t)$ and $\hat{T}_{{\rm i}}(x,t)$ for open-loop estimation by setting all the observer gain to be zero is depicted in Fig. \ref{Ch7:fig:2} (a),  and those for the proposed estimation are depicted in Fig. \ref{Ch7:fig:2} (b) and Fig. \ref{Ch7:fig:2-2} (a)-(b) with observer gains \eqref{Ch7:p1gain}--\eqref{Ch7:p4gain}, respectively, by using input data on January. For the proposed estimation, we fix the parameters of $c = $3.0 $\times$ 10$^{-5}$ and $\ep = $ 1.0 $\times$ 10$^{-8}$, and use the parameter of $\lambda = $5.0 $\times$ 10$^{-6}$ in Fig. \ref{Ch7:fig:2} (b), $\lambda = $1.0 $\times$ 10$^{-5}$ in Fig. \ref{Ch7:fig:2-2} (a), and $\lambda = $5.0 $\times$ 10$^{-7}$ in Fig. \ref{Ch7:fig:2-2} (b). The figures show that the backstepping observer gain makes the convergence speed of the estimation to the actual value approximately 5 to 10 times faster at every point in sea ice. As seen in Fig. \ref{Ch7:fig:2-2}, while the larger value of $\lambda$ makes the convergence speed faster, it causes more overshoot beyond the actual temperature. Hence, the tradeoff between the convergence speed and overshoot can be handled by tuning the gain parameter $\lambda$ appropriately, thereby the parameters used in (b) achieve the desired performance. The overshoot behavior is noted at the end of Section \ref{Ch7:sec:Stability} from a theoretical perspective. Consequently, the stability properties stated in Theorem \ref{Ch7:theorem} for the simplified model can be observed in numerical results of the proposed estimation applied to the original model \eqref{Ch7:scsys1}-\eqref{Ch7:scsys6}. To visualize the convergence of the estimated temperature profile used in \ref{Ch7:fig:2} (b) more clearly, Fig. \ref{Ch7:fig:3} illustrates the profiles of both true temperature (black solid) and estimated temperature (red dash) on January 1st to 3rd in (a)--(c), respectively. We observe that the estimated temperature profile becomes almost the same as the true temperature profile on January 3rd, which is two days after the estimation algorithm runs. Moreover, Fig. \ref{Ch7:fig:3} (d) depicts the time evolution of $\tilde H(t)$, which is an estimation error of the ice's thickness. We observe that the error is ``enlarged" from $\tilde H(0) = 0$ due to the error of temperature profile, and returns to zero after the temperature profiles become almost indistinguishable on January 3rd, from which the necessity of the estimator of the ice's thickness is ensured while the thickness is actually measured.

\begin{figure*}[htb!]
\begin{center}
\subfloat[Estimated temperature converges to the true temperature with a modest error. ]{\includegraphics[width=4.0in]{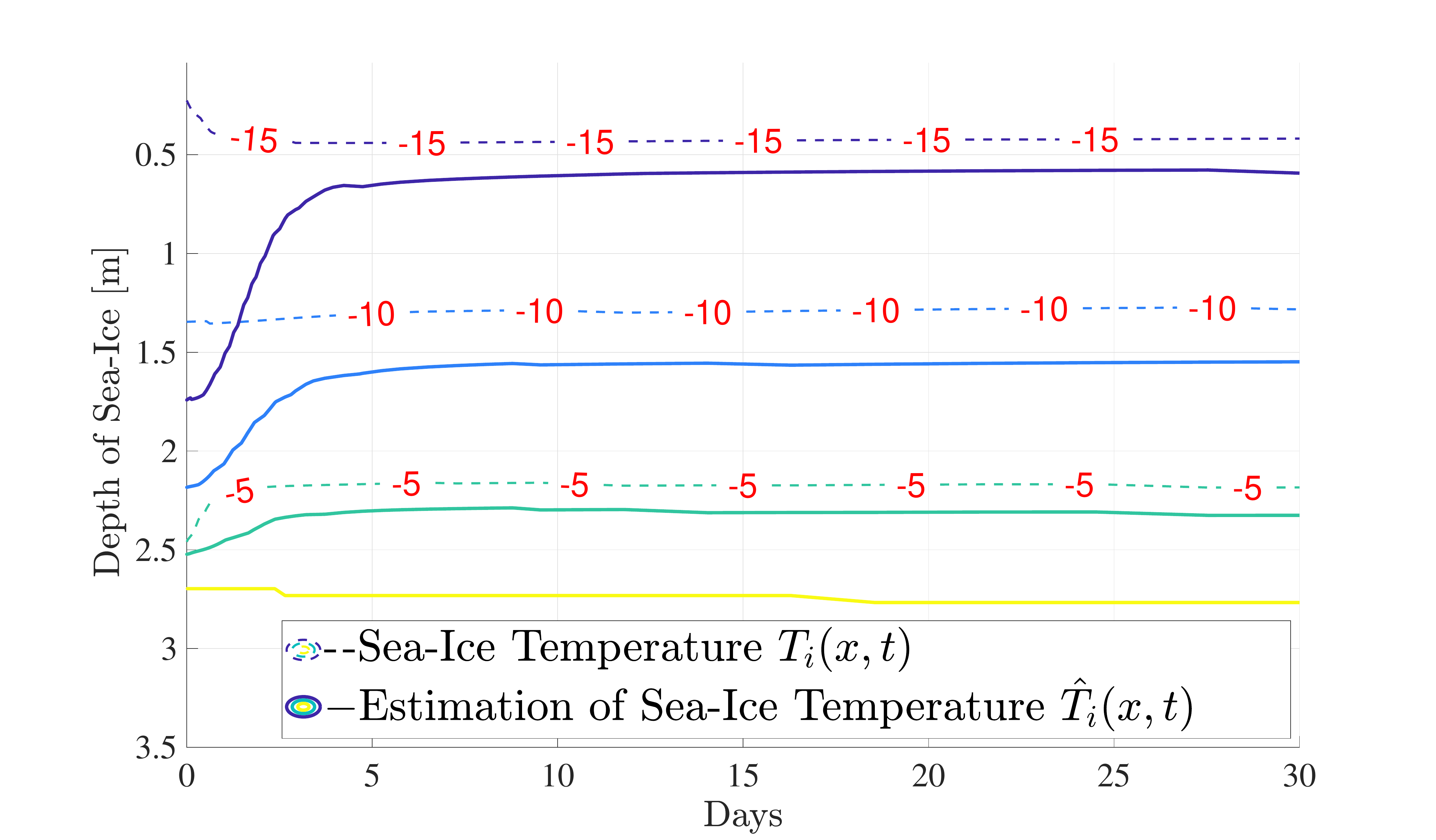}}\\
\subfloat[$\tilde H(t)$ dynamically varies first and stays at a value near zero after 5 days. ]{\includegraphics[width=4.0in]{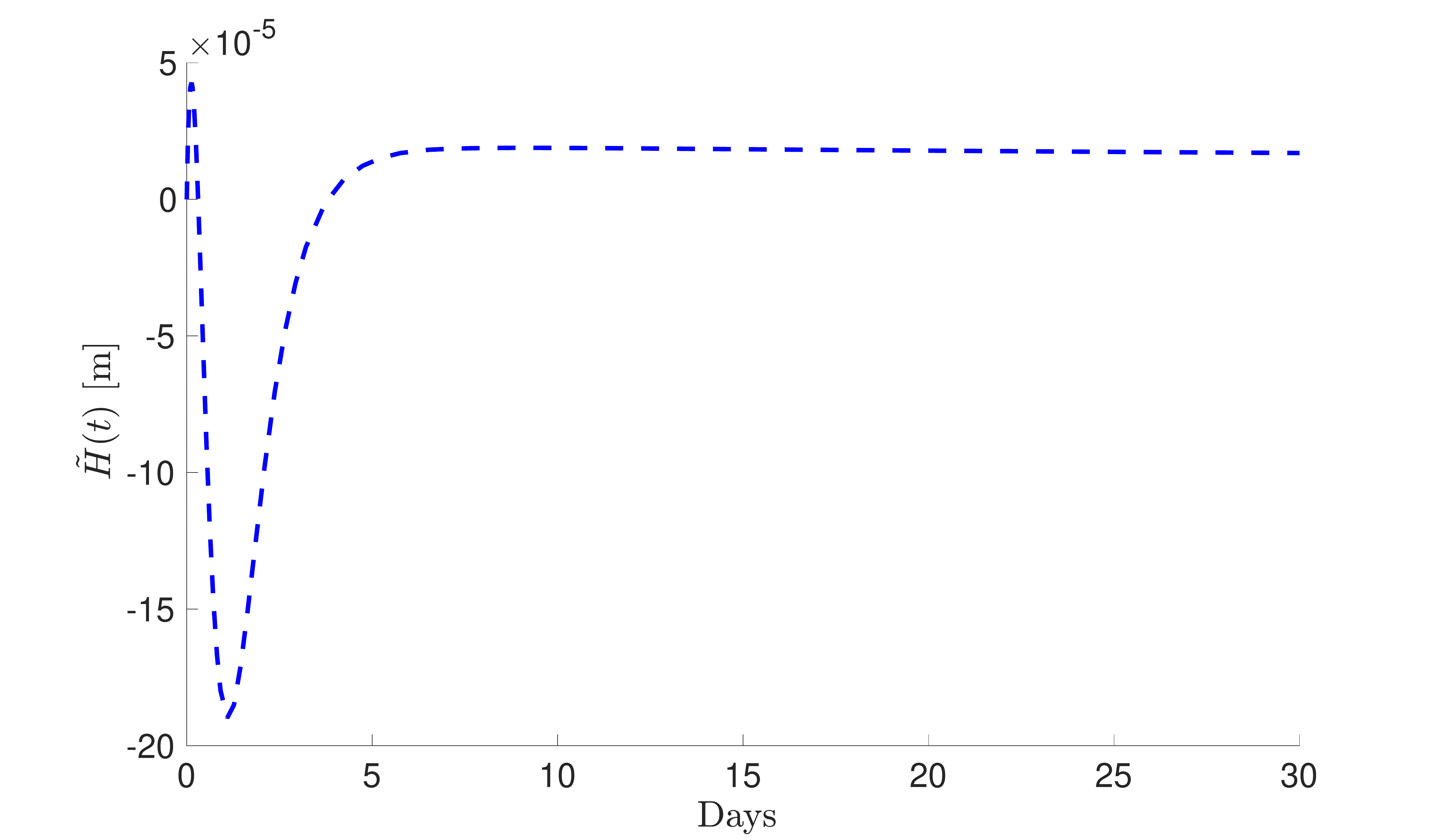}}
\caption{Robustness of the proposed estimation with significant parametric errors: 30[\%] in diffusion coefficient $D_{\rm i}$, 30[\%] in latent heat parameter $\beta$, and 30[\%] in heat flux $F_{\rm w}$ from the ocean.  }
\label{Ch7:fig:4}
\end{center}
\end{figure*}

Finally, we have studied the robustness of the proposed observer by varying the parameters $D_{\rm i}$, $\beta$, and $F_{\rm w}$ in the observer \eqref{Ch7:obsys1}-\eqref{Ch7:obsys4} and the gains \eqref{Ch7:p1gain}--\eqref{Ch7:p4gain} to $D_{\rm i} (1 + \delta_1)$, $\beta(1+\delta_2)$, and $F_{\rm w}(1+\delta_3)$ with setting $\delta_1 = 0.3$, $\delta_2 = -0.3$, and $\delta_3 = 0.4$. 
%namely, the parameters have error 30[\%] in diffusion coefficient $D_{\rm i}$, 30[\%] in latent heat parameter $\beta$, and 40[\%] in heat flux $F_{\rm w}$ from the ocean, respectively. 
Fig. \ref{Ch7:fig:4} (a) shows the contour plots of estimated and true temperature profiles and Fig. \ref{Ch7:fig:4} (b) shows the evolution of $\tilde H(t)$. From both figures, we can see that the observer states converge and stay around the true states with a modest error after 5 days, which illustrates robust performance of the proposed observer under the parameters' uncertainties. 

%%%%%%%%%%%%%%%%%%%%%%%%%%%%%%%%%%%%%%%%%%%%%%%%%%%%%%%%%%%%%%%%%%%%%%%%%%%%%%%%%%%%%%%%%%%%%%%%%%%%%%%%%

section{Lithium-Ion Batteries}
\label{Ch8:intro} % Always give a unique label
% use \chaptermark{}
% to alter or adjust the chapter heading in the running head

%{\bf Motivation:}
%However, none of the existing estimation algorithms consider the phase transition phenomena in LFP or any other electrode material with this characteristic until our recent contribution in \cite{koga2017dscc}.
 
% The estimation problem for batteries with LFP electrodes have been relatively less studied.
% The first approach for state estimation of LFP was proposed in \cite{schwunk2013} using particle filter. A simultaneous state and parameter estimation of LFP was designed in \cite{li2014} via Sequential Monte Carlo filter. However, the methodologies used in literature are commonly based on 
%% model-free estimation,
% and the electrochemical model-based estimation for phase transition electrodes have not been designed. 

%The \emph{core-shell} model proposed for phase transition electrodes is described by a parabolic PDE with a state-dependent moving boundary. 
%This is the so called \emph{Stefan problem}, 
%derived  originally to model liquid-solid phase transition phenomena\cite{gupta2003}. 
%Recently, a control and state estimation technique for the \emph{Stefan problem} was developed in \cite{koga2017,koga2019delay,koga2019arctic}. There, the authors introduced a backstepping design \cite{krstic2008} for this problem and showed the exponential stability of the closed-loop system under some particular (physical) constraints on the moving interface.

\subsection{Battery Management Systems} 

Battery management is crucial for safe and efficient use of numerous kinds of electronics such as smartphones and laptops, and electric vehicles. Among several chemical materials used for electrodes of lithium-ion batteries\index{lithium-ion batteries}, Lithium Iron Phosphate\index{Lithium Iron Phosphate} (LFP) has several attractive features as an active material in lithium-ion batteries such as thermal safety, high energy, and power 
density \cite{padhi1997}.
LFP and other common active materials show unique charge-discharge characteristics due to an underlying crystallographic solid-solid phase transition. Electrochemical models\index{electrochemical models} for lithium-ion batteries with single phase materials do not allow to capture these unique characteristics and thus a mathematical description of phase transitions needs to be added to these models.
Electrochemical models are of interest for the design of accurate estimation algorithms in battery management systems.\index{battery management systems} 
Estimation algorithms based on these models provide visibility into operating regimes that induce degradation enabling a larger domain of operation, therefore, increasing the performance of the battery in terms of energy capacity, power capacity, and fast charge rates \cite{chaturvedi2010,perez2015enhanced}. Electrochemical model-based estimation is challenging for several reasons. First, measurements of lithium concentrations outside specialized laboratory environments is impractical. Second, the concentration dynamics are governed by coupled and nonlinear partial differential algebraic equations (PDAE) \cite{thomas2002mathematical}. Finally, the only measurable quantities (voltage and current) are related to dynamic states through a nonlinear function. 

\begin{figure}[t]
	\begin{center}
	\includegraphics[width=0.9\linewidth]{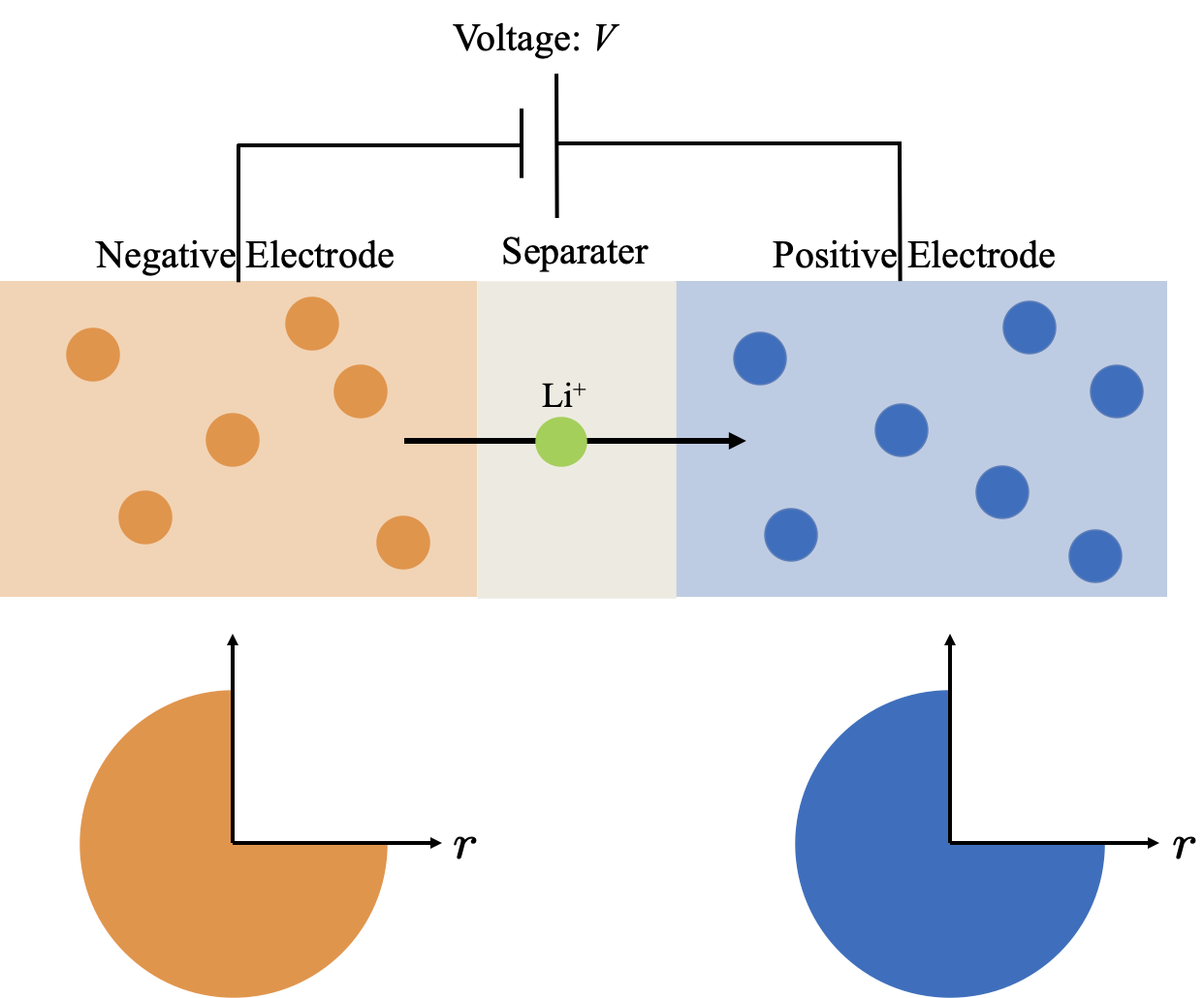}
	\caption{ Schematic of lithium-ion battery and the description of particles in electrochemical models. The concentration dynamics of lithium-ion is governed on the geometry of each particle. }
	\label{Ch8:battery-sch}
	\end{center} 
\end{figure}

%{\bf Relevant Literature:} 
Electrochemical models describe the relevant dynamic phenomena in lithium-ion cells: diffusion, intercalation and electrochemical kinematics (see Figure \ref{Ch8:battery-sch}). These models predict accurately the internal states of the battery, however, their complexity renders a challenging problem for estimation algorithms. For this reason, most approaches develop estimation algorithms  based on simplified models. Among the various simplified models, the single particle model\index{single particle model} (SPM) has been broadly used in the  observer design problem, see \cite{Moura2014,di2010lithium,wang2014adaptive,dey2015nonlinear,perez2015sensitivity,moura2016battery,tang2017state}. The main characteristic of the SPM is the use of a single spherical particle to represent diffusion of lithium ions in the intercalation sites of the porous active materials in the electrodes.

LFP has been extensibility used in lithium ion cells due to its thermal stability, cost effectiveness, non-toxic nature, and long cycle life \cite{padhi1997}. 
An electrochemical model for LFP batteries was proposed in \cite{srinivasan04} based on a \emph{core-shell} model, where the concentration at the core is assumed constant and diffusion is allowed for the phase in the shell.
The LFP model with phase transition electrode was revisited in \cite{zhang2007} with a more complete \emph{ core-shell} model, allowing diffusion in both phases of an LiCoO$_2$ cathode. 

The estimation problem for batteries with LFP electrodes has been relatively less studied;
a particle filter was derived in \cite{schwunk2013}
and a Sequential Monte Carlo filter was derived 
in \cite{li2014}. The \emph{core-shell} model proposed for phase transition electrodes is described by a parabolic PDE with a state-dependent moving boundary.

\subsection{Electrochemical Model with Phase Change Electrode}

\begin{figure}
	\begin{center}
	\includegraphics[width=0.9\linewidth]{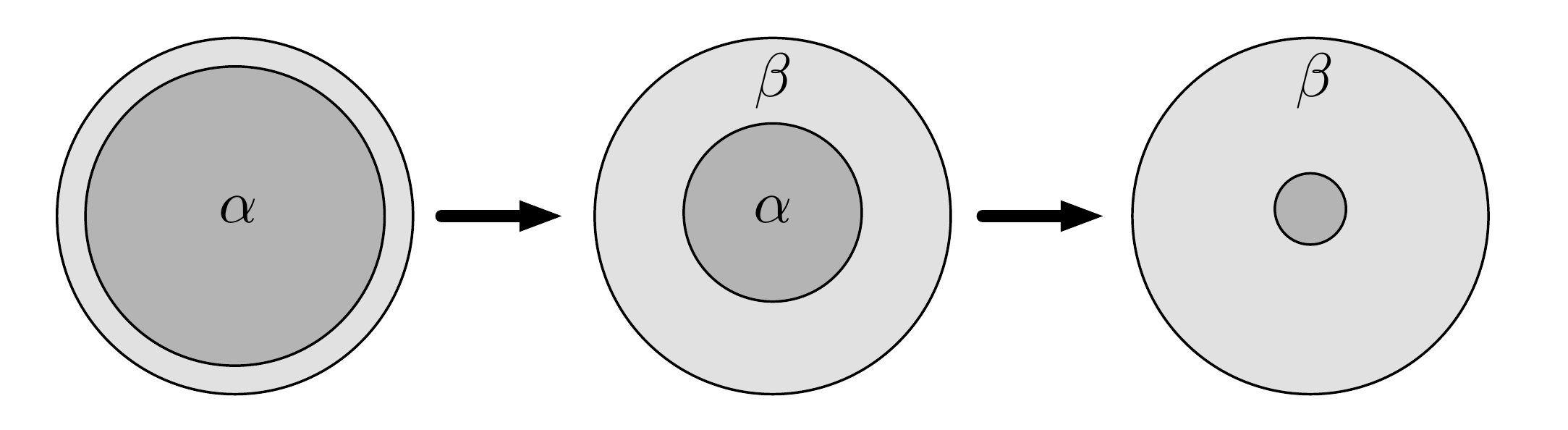}
	\caption{ Phase transition in the positive particle during discharge. 
	The particle starts with a large core of low concentration phase $\alpha$
	and a small shell of high concentration phase $\beta$.
	During discharge there is a positive flux of lithium ion in the surface 
	of the positive particle, increasing the concentration and increasing the 
	size of the $ \beta$-phase shell.}
	\label{Ch8:diagram}
	\end{center} 
\end{figure}

The electrochemical model for lithium-ion cells with a phase transition material in the positive electrode follows \cite{srinivasan04}. 
We restrict the problem to particular initial conditions of the concentration of lithium ions in the particles (i.e. intercalation sites) of the positive electrode and consider only discharge processes.
The initial concentration of lithium ions in the particles of the positive electrode follows a core-shell configuration where the core has a constant distribution of lithium ions in a low concentration phase (the $\alpha$ phase), and the shell has a constant distribution of lithium ions in a high concentration phase (the $\beta$ phase).
During discharge, the fluxes of lithium ions at the surface of the  particles in the positive electrode are positive, thus, increasing the concentration of lithium ions in the shell and the phase boundary is moving to the center, i.e., a \emph{shrinking core} process, as depicted in Figure \ref{Ch8:diagram}.
% \textbf{Add diagram}

\subsubsection*{ Single Particle Model}
The single particle model is a simple electrochemical model that accounts for some phenomena in lithium-ion cells. The main simplification in this 
model comes from the assumption that a single diffusion equation in an spherical particle can be used to model the diffusion of lithium ions in all the intercalation sites of the active material of each electrode.
 In the SPM, the ionic molar fluxes $j_{\mathrm{n},\pm}(t)$ on both electrodes are proportional the current density $I(t)$ applied to the cell
\begin{equation} \label{Ch8:eqn:jn:neg}
	 j_{\mathrm{n},\pm}(t) = \mp \frac{ I(t) }{ a_{\mathrm{s},\pm} FL_{\pm} },  
\end{equation}
where $ a_{\mathrm{s},\pm} = 3 \epsilon_{ \mathrm{s}, \pm }/R_{\mathrm{p},\pm} $ 
is the interfacial area (per unit volume), $\epsilon_{ \mathrm{s}, \pm }$ is the 
volume fraction of active material in each electrode, $ R_{\mathrm{p},\pm}  $ is the averaged radius of the intercalation sites (particles) in the electrodes, $F$ is the Faraday constant, and $L_{\pm}$ is the thickness of each electrode. Throughout this section, the subscripts $+$ and 
$ - $ indicate that the variable corresponds to the  positive or negative particle. 
The concentration dynamics of lithium ions in the negative electrode (single phase) follow the Fick's law\index{Fick's law} for diffusion
\begin{equation} \label{Ch8:eqn:cs:neg:PDE} 
\frac{\partial c_{\mathrm{s},-}}{\partial t}(r,t) 
= \frac{D_{\rm s,-}}{r^{2}} \frac{\partial}{\partial r}
	\left[ r^{2} \frac{\partial c_{\mathrm{s},-} }{\partial r} (r,t) \right] , 
\end{equation}
for $ r \in (0,R_{\rm p,-}) $, $ t>0$ with boundary conditions
\begin{align}
	\label{Ch8:eqn:cs:neg:BC1} \frac{\partial c_{\mathrm{s},-}}{\partial r}(0,t) & = 0, \\
	\label{Ch8:eqn:cs:neg:BC2} D_{\mathrm{s},-} \frac{\partial c_{\mathrm{s},-}}{\partial r}(R_{\mathrm{p},-},t) & = - j_{\mathrm{n},-}(t),
\end{align}
and initial condition $ c_{\rm 0,-} \in \mathcal{C}(0,R_{\rm p,-}) $. 
%The positive electrode undergoes a phase transition between charged and discharged states as seen in LFP. Hence, the dynamics in the positive electrode is separate between $\alpha$ phase in the core and $\beta$ phase in the shell. 
Diffusion in the positive particle follows a core-shell model\index{core-shell model}. In the core of the particle, i.e., for $ r \in (0,r_{\rm p}(t) )$, lithium ions are in the $ \alpha$-phase. The concentration in the core is assumed to be constant and equal to the equilibrium value of the  $\alpha$-phase, i.e., $c_{s,+}(r) = c_{\mathrm{s},\alpha} $ for all $ r \in (0,r_{\rm p}(t) )$ . In the shell of the spherical particle, i.e. for $ r \in (r_{\rm p}(t),R_{p,+})$, the
concentration of lithium ions is in $\beta$-phase.
The concentration dynamics of lithium-ions in the shell of the positive particle follows the Fick's law for diffusion
\begin{equation} \label{Ch8:eqn:cs:pos:PDE} 
\frac{\partial c_{\mathrm{s},+}}{\partial t}(r,t) 
= \frac{ D_{\rm s,+} }{r^{2}} \frac{\partial}{\partial r}
\left[ r^{2} \frac{\partial c_{\mathrm{s},+} }{\partial r} (r,t) \right] , 
\end{equation}
for $r \in (r_{\rm p}(t), R_{\rm p,+}) $ with boundary conditions
\begin{align}\label{Ch8:eqn:cs:pos:BC1} 
	c_{\mathrm{s},+}(r_{\rm p}(t),t) & = c_{\mathrm{s}, \beta}, \\
\label{Ch8:eqn:cs:pos:BC2}	D_{\mathrm{s},+} \frac{\partial c_{\mathrm{s},+}}{\partial r}(R_{\rm p,+},t) & = - j_{\mathrm{n},+}(t),
\end{align}
and initial conditions $ c_{0,+} \in \mathcal{C} (r_{\rm p}(0), R_{\rm p,+}) $. 
The time-evolution of the moving interface $r_{p}(t)$ is not given explicitly. Instead, mass balance\index{mass balance} at the moving interface yields the following state-dependent dynamics:
\begin{align}\label{Ch8:eqn:cs:pos:ODE} 
(c_{\mathrm{s}, \beta} -c_{\mathrm{s}, \alpha} ) 
\frac{d r_{\rm p}(t)}{dt} = - D_{\rm s,+} \frac{\partial c_{\mathrm{s},+}}{\partial r}(r_{\rm p}(t),t). 
\end{align}
Overpotentials $ \eta_{\pm}(t) $ are found by solving the nonlinear algebraic equation
\begin{align}
j_{\mathrm{n},\pm}(t) 
& =\frac{i_{0,\pm}(t)}{F} \left[e^{\frac{\alpha_{\mathrm{a}}F}{RT} \eta_{\pm}(t)} - e^{-\frac{\alpha_{\mathrm{c}}F}{RT} \eta_{\pm}(t)} \right], \label{Ch8:eqn:bv:spm} \\
i_{0,\pm}(t) & = F k_{\pm} \left[ {c}_{\mathrm{ss},\pm}(t) \right]^{\alpha_{\mathrm{c}}} \left[c_{\mathrm{e},0} \left(c_{\mathrm{s},\max,\pm} - c_{\mathrm{ss},\pm}(t)  \right) \right]^{\alpha_{\mathrm{a}}}, \label{Ch8:eqn:i0:spme}
\end{align}
where $ {c}_{\mathrm{ss},\pm}(t) :=  {c}_{\mathrm{s},\pm}(R_{\rm p,\pm},t)$. 
The electric potential in each electrode is given by
\begin{equation}  \label{Ch8:phis}
\phi_{\mathrm{s},\pm}(t) =  \eta_{\pm}(t) + U_{\pm}(c_{\mathrm{ss},\pm}(t)) + R_{\mathrm{f},\pm}F j_{\mathrm{n},\pm}(t).
\end{equation}
Finally, output voltage is computed as the difference
between the electric potential in each electrode
\begin{align}\label{Ch8:eq:voltage}
V(t) = 	\phi_{\mathrm{s},+}(t)  - 	\phi_{\mathrm{s},-}(t).
\end{align}
Equations \eqref{Ch8:eqn:cs:pos:PDE} -\eqref{Ch8:eq:voltage} form a complete description of the single particle model with a phase transition electrode, and it provides the following property on the moving interface during the discharge process. 
\begin{remark}\label{Ch8:rem:shrinking} 
During the single discharge process,  the current density $I(t)$ maintains positive, i.e. $I(t)>0$ for $\forall t>0$. This current positivity ensures the moving interface being shrinking. Furthermore, the initial interface position is less than the cell radius. Hence, 
\begin{align}\label{Ch8:property1}
 \frac{d r_{\rm p}(t)}{dt} < 0, \\
\label{Ch8:property2} 0\leq  r_{\rm p}(t)< R_{\rm p,+}.
 \end{align}
\end{remark}

\subsubsection*{Mass Conservation}\label{Ch8:sec:mass}
In this model, the total amount of lithium ions is conserved.
The mathematical description of this property is given in the following lemma. 
\begin{lem}\label{Ch8:lem:mass}
	The total amount of lithium $n_{\rm Li}$ in solid phase ( moles per unit area ) defined as
	\begin{equation}\label{Ch8:eqn:mass}
	n_{\rm Li}(t) = 
	\epsilon_{ \rm s,-} L_{-}   \overline{ c }_{ \rm {s}, -} (t) + 
	\epsilon_{ \rm s,+} L_{+}  \overline{ c }_{ \rm {s},+} (t) ,
	\end{equation}
	where $ \overline{ c }_{ \rm {s}, - }(t) $ and 
	$ \overline{ c }_{ \rm{s} , +} (t) $ are the volumetric averages of the concentrations
\begin{align}
		\overline{c} _{ \rm {s},-} (t) = & 
		\frac{3}{ R_{ \rm{p},- }^3 } \int_{0}^{ R_{ \rm{p},-} }  c_{ \rm{s},-} (r,t) r^{2} {\rm d}r, \\
		\overline{c}_{ \rm {s}, +} (t) = & 
		\frac{3}{ R_{ \rm{p},+ }^3 } \int_{0 }^{ R_{ \rm{p},+} } c_{ \rm{s},+} (r,t) r^{2} {\rm d}r,
	\end{align}
is conserved, namely $ dn_{ \rm{Li} }(t)/dt = 0$.
	
\end{lem} 

Lemma \ref{Ch8:lem:mass} was derived in \cite{klein2013} for electrodes with a single phase, and we can show that this  result extends to electrodes with phase 
transition materials.

\textbf{Proof:}\\
In our problem formulation there is a single phase in the negative particle 
and there are two phases in the positive particle, i.e., $\alpha$-phase in the core and $\beta$-phase in the shell.
The concentration in $\alpha$-phase at the core is assumed to be constant (at its equilibrium value $ c_{\mathrm{s}, \alpha}$). 
Under these assumptions, the time derivative of \eqref{Ch8:eqn:mass} is given by
	\begin{align}
	\frac{dn_{\rm Li}}{dt}(t) = & - a_{\rm{s},-} L_{-} j_{\rm{n},-}(t) - a_{\rm{s},+} L_{+} j_{\rm{n},+}(t)  - \frac{3 \epsilon_{\rm{s},+} L_{+} }{R_{\rm{p},+}^3} 
													 r_{\rm{p}}^2(t)\notag\\ 
												     & \times \left[  \frac{ dr_{ \rm{p}} }{ dt }(t) 
													 \left[c_{\mathrm{s}, \beta} -c_{\mathrm{s}, \alpha} \right] + D_{\rm s,+} \frac{\partial c_{\mathrm{s},+}}{\partial r}(r_{\rm p}(t),t)\right]. 
\end{align}
Hence, the molar flux equations in \eqref{Ch8:eqn:jn:neg} and the dynamics of the moving interface in \eqref{Ch8:eqn:cs:pos:ODE} lead to $dn_{ \rm{Li} }(t)/dt = 0$.
In a more general formulation introduced in \cite{khandelwal2014,khandelwal2015thermally}, i.e. when both electrodes have multiple phase transitions not necessarily at the equilibrium, mass conservation of lithium ions is guaranteed with the following interface dynamics
\begin{align}
	 \frac{ dr^{\rm [a,b]}_{i} }{ dt }(t)  =& \frac{ 1 }{ c_{\rm{b}} - c_{\rm{a}} }
	 \left[ D_{\rm a}\frac{ \partial c }{ \partial r }( r_{i}^{ \rm [a,b] }(t)^{-} ,t ) 	        - D_{\rm b}\frac{ \partial c }{ \partial r }( r_{i}^{\rm [a,b] }(t)^{+}, t ) \right],
\end{align}
where $r^{\rm [a,b]}_{i}$ is the interface radius between any two phases (phase a and phase b) 
in any electrode.
Each phase has a distinct equilibrium $c_{\rm a}$, $c_{\rm b}$ and diffusion coefficient $D_{a}$, $D_{b}$.
%\end{proof}

\subsection{State-of-Charge\index{state-of-charge} Estimation}

Now, a state estimation algorithm for concentration of lithium ions, in both negative and positive electrodes, is provided in this section from the single particle model. The state observer\index{state observer} for the positive electrode is derived via the backstepping method  for moving boundary PDEs, and the observer for the negative electrode is derived from the mass conservation property. 

\subsubsection*{Observer for Phase Transition Positive Electrode}
The state observer is a copy of the diffusion system \eqref{Ch8:eqn:cs:pos:PDE}-\eqref{Ch8:eqn:cs:pos:BC2} in the positive electrode 
together with output error injection
\begin{align}  \label{Ch8:eqn:cs:obs:PDE} 
\frac{\partial \widehat{c} _{\mathrm{s},+}  }{ \partial t }(r,t) 
= & 
\frac{ D_{\rm s,+} }{ r^{2} } \frac{\partial}{\partial r} 
\left[ r^{2} \frac{ \partial \widehat{ c }_{\mathrm{s},+} }{ \partial r } (r,t) \right]\notag\\
& + 
P(\widehat{ r }_{\rm p}(t),r) \left[c_{ \mathrm{ss,+} }(t) - 
\widehat{ c}_{ \mathrm{s,+} } (R_{\rm p,+},t) \right] , 
\end{align}
for $r \in (\widehat{r}_{\rm p}(t), R_{\rm p,+}) $ 
with boundary conditions
\begin{align} \label{Ch8:eqn:cs:obs:BC1} 
\widehat{c}_{\mathrm{s},+} (\widehat{r_{\rm p}}(t),t)
 = & c_{\beta}, \\
\label{Ch8:eqn:cs:obs:BC2}
 D_{\mathrm{s},+} 
 \frac{\partial \widehat{c}_{\mathrm{s},+}}{\partial r}(R_{\rm p,+},t)
 = & - j_{\mathrm{n},+}(t) \notag\\
& + Q( \widehat{r_{\rm p}}(t)) 
\left[ c_{ \mathrm{ss,+} }(t) - \widehat{ c } _{ \mathrm{s,+} } (R_{\rm p,+},t) \right]  ,
\end{align}
and initial conditions $\widehat{c}_{0,+} \in \mathcal{L}^2(\widehat{r_{\rm p}}(0), R_{\rm p,+}) $ and $\widehat{r_{\rm p}}(0) \in (0, R_{\rm p,+}) $.
Observer gains are given by
\begin{align}\label{Ch8:Pgain}
P( \widehat{r_{\rm p}}(t),r )
& = D_{s,+} \overline{\lambda}^2\frac{ R_{ \rm p,+ } }{r } l(t) s(t)
\frac{ I_{2} \left( z(t) \right) }{ z(t) },  \\
\label{Ch8:Qgain}
Q( \widehat{r_{\rm p}}(t) )
& = \frac{ D_{s,+}}{ R_{ \rm p,+ } }\left( \frac{\overline{\lambda}}{2}s(t) + 1 \right) ,
\end{align}
where $ I_2(\cdot) $ is a modified Bessel function\index{modified Bessel function} of the second kind and 
\begin{align}
\overline{\lambda}  & = \frac{\lambda}{D_{s,+}} , \\
s(t) & = R_{\rm p,+} - \widehat{r_{\rm p}}(t), ~~~ l(t)  = r  - \widehat{r_{\rm p}}(t), \\
z(t) & = \sqrt{\overline{\lambda} \left[s(t)^2 - l(t)^2 \right]} . 
\end{align}
The parameter $\lambda >0$ is designed to achieve faster convergence of the estimated concentration to true concentration. Moreover, the estimator for the moving interface position is given by the following dynamics: 
\begin{align} 
(c_{\mathrm{s}, \beta} -c_{\mathrm{s}, \alpha} ) \frac{d \widehat{r_{\rm p}}(t)}{dt} 
= - \kappa \left[c_{ \mathrm{ss,+} }(t) - \widehat{ c }_{ \mathrm{s,+} } ( R_{\mathrm{p},+}, t ) \right] \notag \\ 
- D_{\rm s,+} \frac{ \partial \widehat{ c } _{\mathrm{s},+} }{ \partial r } (\hat{r}_{\rm p}(t),t), 
\label{Ch8:ODE-observer}
\end{align} 
where the parameter $\kappa>0$ is designed to achieve fast convergence of the estimated interface position to the true value. 

The stability of the estimation error system is theoretically proven for the PDE observer \eqref{Ch8:eqn:cs:obs:PDE}--\eqref{Ch8:eqn:cs:obs:BC2} with gains \eqref{Ch8:Pgain}, \eqref{Ch8:Qgain} under the assumption $\widehat{r_{\rm p}}(t) \equiv r_{\rm p}(t)$ for all $t \geq 0$ in the next section. As the moving interface position $r_{\rm p}(t)$ is unknown in practice, we construct the estimator \eqref{Ch8:ODE-observer}, and use the estimated interface position $\widehat{r_{\rm p}}(t)$ in the gains \eqref{Ch8:Pgain}, \eqref{Ch8:Qgain} of PDE observer.

The sign of the observer gain in \eqref{Ch8:ODE-observer} (first term in the right hand side) is determined based on the monotonic relation, namely, as the surface concentration $c_{ \mathrm{ss,+} }(t)$ is increased the moving interface position $r_{\rm p}(t)$ is decreased. Physically, as the battery is discharged, the domain of the lithium rich $\beta$-phase in the positive electrode is expanded from the outer region. Hence, the observer \eqref{Ch8:ODE-observer} is designed so that if the measured surface concentration is larger than the estimated surface concentration, the battery is discharged more than estimated, and the domain of $\beta$-phase for the estimator is driven to be expanded.

\subsubsection*{Stability Analysis of the Estimation Error System with Known Interface Position}\label{Ch8:sec:stability}
%A rigorous proof of the convergence of the estimated lithium-ion concentration governed by PDE observer \eqref{Ch8:eqn:cs:obs:PDE}--\eqref{Ch8:eqn:cs:obs:BC2} to the true concentration governed by \eqref{Ch8:eqn:cs:pos:PDE}--\eqref{Ch8:eqn:cs:pos:BC2} is provided in this section under the assumption $\widehat{r_{\rm p}}(t) \equiv r_{\rm p}(t)$ for all $t \geq 0$. 
%\begin{assum}
%The interface position $r_{\rm p}(t)$ is known and used in PDE observer \eqref{Ch8:eqn:cs:obs:PDE}--\eqref{Ch8:eqn:cs:obs:BC2} by setting $\widehat{r_{\rm p}}(t) \equiv r_{\rm p}(t)$ for all $t \geq 0$. 
%\end{assum}

Let $\widetilde{ c }_{\mathrm{s},+} (r,t) $ be an estimation error defined by 
$ \widetilde{c}_{\mathrm{s},+} (r,t) := c_{\mathrm{s},+}(r,t) - \widehat{c}_{ \mathrm{s},+ }(r,t) $. 
The stability analysis of the estimation error system is presented in the following theorem. 
\begin{thm}\label{Ch8:thm:main}
Consider the plant PDE \eqref{Ch8:eqn:cs:pos:PDE}--\eqref{Ch8:eqn:cs:pos:BC2} and the PDE observer \eqref{Ch8:eqn:cs:obs:PDE}--\eqref{Ch8:eqn:cs:obs:BC2} with observer gains \eqref{Ch8:Pgain} and \eqref{Ch8:Qgain} under the properties of \eqref{Ch8:property1}, \eqref{Ch8:property2}, and the assumption $\widehat{r_{\rm p}}(t) \equiv r_{\rm p}(t)$ for all $t \geq 0$. Then, for any initial estimation error $\widetilde{ c_{\mathrm{s},+} }(r,0)$, the estimation error is exponentially stable at the origin in the sense of the norm 
\begin{align}
\int_{r_{\rm p}(t)}^{R_{\rm p,+}} r^2 \widetilde{ c }_{\mathrm{s},+} (r,t)^2 {\rm d}r. 
\end{align}
\end{thm}
Note that subtracting \eqref{Ch8:eqn:cs:obs:PDE}-\eqref{Ch8:eqn:cs:obs:BC2} from \eqref{Ch8:eqn:cs:pos:PDE}-\eqref{Ch8:eqn:cs:pos:BC2}
under $\widehat{r_{\rm p}}(t) \equiv r_{\rm p}(t)$ yields the estimation error dynamics 

\begin{align} \label{Ch8:eqn:cs:err:PDE} 
\frac{\partial \widetilde{c}_{ \mathrm{s},+ } }{ \partial t }(r,t) 
 = & 
\frac{ D_{\rm s,+} }{ r^{2} } \frac{ \partial }{ \partial r } 
\left[ r^{2} \frac{ \partial \widetilde{c}_{\mathrm{s},+} }{ \partial r }(r,t) \right] - P(r_{\rm p}(t),r) \widetilde{c}_{ \mathrm{s},+} (R_{\rm p,+},t), \\ 
\label{Ch8:eqn:cs:err:BC1} 
\widetilde{c}_{\mathrm{s},+} (r_{\rm p}(t),t)
 =& 0, \\
\label{Ch8:eqn:cs:err:BC2}
 D_{\mathrm{s},+} 
 \frac{ \partial \widetilde{c}_{\mathrm{s},+} }{ \partial r }(R_{\rm p,+},t)
 = & - Q \left( r_{ \mathrm{p} }(t) \right) \widetilde{ c }_{ \mathrm{s,+} } ( R_{\rm p,+},t ). 
\end{align}

%The proof of Theorem \ref{Ch8:thm:main} is established by showing the stability of the estimation error system \eqref{Ch8:eqn:cs:err:PDE}--\eqref{Ch8:eqn:cs:err:BC2} through the remainder of Section \ref{Ch8:sec:stability}. 

\textbf{Change of coordinate:} 
First, we introduce the following change of coordinate and state variable to simplify the structure of the estimation error dynamics in a cartesian coordinate:
\begin{align}\label{Ch8:ch-co}
x & = R_{\rm p,+} - r, \\
\label{Ch8:ch-state} \widetilde{u}(x,t) & =  r  \widetilde{c}_{\mathrm{s},+} (r,t) , \\
\label{Ch8:s-def} s(t) & =  R_{\rm p,+} - \widehat{r_{\rm p}}(t). 
\end{align}
%Then, 
%\begin{align} 
%\widetilde{u}_x(x,t) =& -  \widetilde{ c_{\mathrm{s},+} } (r,t)  -  r  \frac{\partial \widetilde{ c_{\mathrm{s},+} }}{\partial r} (r,t) 
%\end{align} 
%Thus, 
%\begin{align} 
%u_x(s(t),t) =&  -  r_{\rm p}(t)  \frac{\partial  c_{\mathrm{s},+} }{\partial r} (r_{\rm p}(t),t) 
%\end{align} 
%
%
The estimation error dynamics \eqref{Ch8:eqn:cs:err:PDE}-\eqref{Ch8:eqn:cs:err:BC2} is rewritten by the new coordinate and state as 
\begin{align}
\label{Ch8:eqn:cs:err2:PDE} 
\frac{\partial \widetilde{ u} }{ \partial t }(x,t) 
& = D_{\rm s,+}\frac{ \partial^2 \widetilde{ u} }{\partial x^2} (x,t) 
- \overline{P}(s(t),x) \widetilde{ u } (0,t)  , \\
 \label{Ch8:eqn:cs:err2:BC1} 
\widetilde{ u }(s(t),t)
& = 0, \\
\label{Ch8:eqn:cs:err2:BC2}
 \frac{\partial \widetilde{u}}{\partial x}(0,t)
& = - \overline{Q}(s(t)) \widetilde{ u }(0,t) ,
\end{align}
where 
\begin{align}
\overline{P} (s(t),x) 
& = \frac{r}{R_{\rm p, +}} P(r_{\rm p}(t),r) , \\
\overline{Q} (s(t)) 
& = \frac{1}{R_{\rm p,+} } - \frac{1}{D_{\mathrm{s},+}} Q(r_{\rm p}(t)) .
\end{align}
With respect to the variable \eqref{Ch8:s-def}, the properties \eqref{Ch8:property1} and \eqref{Ch8:property2} presented in Remark \ref{Ch8:rem:shrinking} are equivalent to 
\begin{align} \label{Ch8:s-property1}
\dot s(t) &>0, \\
0 < s(t) &\leq 	R_{\rm p,+}. \label{Ch8:s-property2}
\end{align}

\textbf{Derivation of observer gains:} 
Consider the following invertible transformation from the estimation error $\widetilde{u}(x,t)$ to the transformed state $\widetilde{w}(x,t)$:
\begin{align}\label{Ch8:eqn:bst1}
\widetilde{w} (x,t) 
& = \widetilde{u}(x,t) + 
\int_{0}^{x} q(\overline{x},\overline{y}) \widetilde{u}(y,t) {\rm d}y, \\
\label{Ch8:eqn:bst2}
\widetilde{u} (x,t) 
& = \widetilde{w}(x,t) + 
\int_{0}^{x} p(\overline{x},\overline{y} ) \widetilde{w}(y,t) {\rm d}y, 
\end{align}
where $\overline{x} = s(t)-x$, $\overline{y} = s(t)-y$. We can show that if the gain kernel functions and the observer gains satisfy the following conditions: 
\begin{align}\label{Ch8:gain_cond_1} 
\frac{\pa^2 p}{\pa \bar x^2}(\bar x,\bar y) -\frac{\pa^2 p}{\pa \bar y^2}(\bar x,\bar y)= &- \bar \lambda p(\bar x,\bar y), \\
p(\bar x,\bar x) =& \frac{\bar \lambda}{2} \bar x, \\
 p(0,\bar y) =& 0, \\
 \frac{\pa^2 q}{\pa \bar x^2}(\bar x,\bar y) -\frac{\pa^2 q}{\pa \bar y^2}(\bar x,\bar y)= & \bar \lambda q(\bar x,\bar y), \\
q(\bar x,\bar x) =& -\frac{\bar \lambda}{2} \bar x, \\
 q(0,\bar y) =& 0,  \label{Ch8:gain_cond_6}  \\
 \overline{P}(s(t),x) = & D_{\rm s,+} p_{\bar y}(\bar x, s(t)), \label{Ch8:Pcond}\\
 \overline{Q}(s(t)) = & - p(s(t),s(t)), \label{Ch8:Qcond} 
\end{align}
then, the following target $\widetilde{w}$-system is obtained:
\begin{align}\label{Ch8:eqn:cs:tar:PDE} 
\frac{\partial \widetilde{ w} }{ \partial t }(x,t) 
= &~ D_{\rm s,+} \frac{ \partial^2 \widetilde{ w} }{\partial x^2}(x,t)-\lambda\widetilde{w}(x,t) + \dot{s}(t) \int_{0}^{x} q'(\overline{x},\overline{y}) \notag\\
& \times \left(\widetilde{w}(y,t)+\int_{0}^{y} p(\overline{y},\overline{z}) \widetilde{w}(z,t) {\rm d}z\right) dy,\\
\label{Ch8:eqn:cs:tar:BC1} \widetilde{w}(s(t),t) = & ~ 0,\\
\label{Ch8:eqn:cs:tar:BC2} \frac{ \partial \widetilde{ w} }{\partial x}(0,t) = &~ 0,
\end{align}
where $q'(\overline{x},\overline{y}) = \fr{\pa q}{\pa \overline{x}} (\overline{x},\overline{y}) + \fr{\pa q}{\pa \overline{y}} (\overline{x},\overline{y})$.  The equations \eqref{Ch8:gain_cond_1}--\eqref{Ch8:gain_cond_6} lead to the following explicit solutions:  
\begin{align}\label{Ch8:p-sol} 
p(\overline x,\overline y) =& \overline{\lambda} \overline x \frac{I_1 \left(\sqrt{\overline{\lambda}\left[\overline{y}^2 - \overline{x}^2 \right]}\right)}{ \sqrt{\overline{\lambda}\left[\overline{y}^2 - \overline{x}^2 \right]}}, \\
q(\overline{x},\overline{y}) =& -\overline{\lambda} \overline{x} \frac{J_1\left(\sqrt{\overline{\lambda}\left[\overline{y}^2 - \overline{x}^2 \right]}\right)}{\sqrt{\overline{\lambda}\left[\overline{y}^2 - \overline{x}^2 \right]}}, \label{Ch8:q-sol} 
\end{align}
with a modified Bessel function $I_1(\cdot)$ and a Bessel function $J_1(\cdot)$ of the first kind, respectively. Substituting the solution \eqref{Ch8:p-sol} to the conditions \eqref{Ch8:Pcond}, \eqref{Ch8:Qcond} (note that $\frac{d I_1(z)}{{\rm d}z} = \frac{I_2(z)}{z}$ for all $z$), and taking back to the original coordinate and variables, the observer gains are derived as \eqref{Ch8:Pgain} and \eqref{Ch8:Qgain}. 

\textbf{Stability proof:} 
We consider the time evolution of the following Lyapunov function:
\begin{align}\label{Ch8:wlyap}
W(t)=\frac{1}{2}\int_0^{s(t)} \widetilde{w}(x,t)^2 {\rm d}x. 
\end{align}
Taking the time derivative of \eqref{Ch8:wlyap} along with \eqref{Ch8:eqn:cs:tar:PDE}-\eqref{Ch8:eqn:cs:tar:BC2} yields 
\begin{align}
\dot{W}(t)=& - D_{\mathrm{s},+} \int_0^{s(t)}\left( \frac{ \partial \widetilde{ w} }{\partial x}(x,t)\right)^2 {\rm d}x - \lambda  \int_0^{s(t)}\widetilde{w}(x,t)^2 {\rm d}x\notag\\
&+\dot{s}(t)\int_0^{s(t)} \widetilde{w}(x,t) \left[\int_{0}^{x} q'(\overline{x},\overline{y})\right. \notag\\
&\left. \left(\widetilde{w}(y,t)+\int_{0}^{y} P(\overline{y},\overline{z}) \widetilde{w}(z,t) {\rm d}z\right) dy\right]{\rm d}x. 
\end{align}
Applying Young's, Cauchy Schwartz, and Poincare's inequalities with the help of the properties \eqref{Ch8:s-property1} and \eqref{Ch8:s-property2}, one can show that there exists a constant $a>0$ such that the following inequality holds:
\begin{align}
\dot{W}(t) \leq & - b W(t) + a \dot{s}(t) W(t), 
\end{align}
where $b = \frac{D_{\mathrm{s},+}}{4 R_{\rm p,+}^2} + \lambda$. With the help of \eqref{Ch8:s-property1} and \eqref{Ch8:s-property2}, it yields the exponential decay of $W(t)$ as 
\begin{align}
W(t) \leq e^{a R_{\rm p,+}} W(0) e^{-b t}.   
\end{align}
Hence, the origin of $\widetilde{w}$-system is shown to be exponentially stable, from which we conclude Theorem \ref{Ch8:thm:main}.  
\subsubsection*{Observer for Negative Electrode}
%The observer design for negative electrode presented in this section imposes the following assumption on the known variables through measurements. 
%\begin{assum} 
%The position and velocity of the interface $ r_{\rm p}(t), \dot{r}_{\rm p}(t)$ along with the surface concentration in the positive electrode $c_{\rm ss}(t)$ are known for any $ t \geq 0 $.
%\end{assum}

The observer design for lithium ion concentration in the negative electrode is constructed by the copy of the dynamics \eqref{Ch8:eqn:cs:neg:PDE}-\eqref{Ch8:eqn:cs:neg:BC2}  together with the output injection of the \emph{positive} electrode
\begin{align} \label{Ch8:eqn:cs:neg:obs} 
\frac{ \partial \widehat{c}_{ \mathrm{s},-} }{ \partial t }(r,t) 
= & \frac{D_{\rm s,-}}{r^{2}} \frac{\partial}{\partial r} \left[ r^{2} \frac{ \partial \widehat{ c }_{ \mathrm{s}, - } }{\partial r} (r,t) \right]   + P_{-}(r_{\rm{p}}(t))  \widetilde{ c }_{ \mathrm{s,+} } (R_{\rm p,+},t) ,
\end{align}
	for $ r \in (0,R_{\rm p,-}) $, $ t>0$ with boundary conditions
\begin{align}
\frac{ \partial \widehat{c}_{\mathrm{s},-} }{\partial r}(0,t) 
= & 0, \\
D_{ \mathrm{s},- } \frac{ \partial \widehat{c}_{\mathrm{s},-} }{ \partial r }( R_{\mathrm{p},-}, t ) 
= & - j_{ \mathrm{n},- }(t) + Q_{-} \left( r_{\rm p}(t) \right) \widetilde{ c }_{ \mathrm{s,+} } (R_{\rm p,+},t). \label{Ch8:eqn:cs:neg:obs:bc2}
\end{align}
Observer gains in the negative electrode are computed to conserve the total amount of lithium ions in the state observer defined as 
\begin{align} 
\widehat{ n }_{{\rm Li}} (t) 
= & 
\frac{ 3 \epsilon_{ \rm s,+} L_{+} }{ R_{ \rm{p},+ }^3 } 
\int_{0 }^{ R_{ \rm{p},+} } \widehat{c}_{ \rm{s},+} (r,t) r^{2} {\rm d}r  + \frac{ 3 \epsilon_{ \rm s,-} L_{-} }{ R_{ \rm{p},- }^3 } \int_{0 }^{ R_{ \rm{p},-} } 
\widehat{c}_{ \rm{s},-} (r,t) r^{2} {\rm d}r. \label{Ch8:eqn:lithium-mass-estimate}
\end{align} 
Taking the time derivative of \eqref{Ch8:eqn:lithium-mass-estimate} along with the dynamics \eqref{Ch8:eqn:cs:obs:PDE}--\eqref{Ch8:ODE-observer} 
and \eqref{Ch8:eqn:cs:neg:obs}--\eqref{Ch8:eqn:cs:neg:obs:bc2} leads to 
\begin{align} 
\frac{ d \widehat{n}_{{\rm Li},+} }{ dt } 
= & 
- a_{ \rm s,+} L_{+} j_{\mathrm{n},+}(t) - a_{ \rm s,-} L_{-} j_{\mathrm{n},-}(t)  + F \widetilde{c} _{ \mathrm{s,+} } (R_{\rm p,+},t), \label{Ch8:eqn:n-der}
\end{align} 
where $F$ is defined by 
\begin{align}
 F = &~ 
 a_{ \rm s,+} L_{+} \left( \frac{\kappa}{R_{ \rm{p},+ }^2} \widehat{r}_{\rm p}(t)^2 
 + Q(\widehat{r}_{\rm p}(t)) \right) 
 + a_{ \rm s,-} L_{-} Q_{-}(\widehat{r}_{\rm p}(t)) \notag \\
 & + 
 \frac{3 \epsilon_{ \rm s,+} L_{+}}{ R_{ \rm{p},+ }^3 } \int_{ \widehat{r}_{\rm p}(t) }^{R_{ \rm{p},+ }} r^2 P( \widehat{r}_{\rm p}(t),r )dr 
 + \epsilon_{ \rm{s},-} L_{-} P_{-}( \widehat{r}_{\rm p}(t)). 
\end{align} 
By the balance of the ionic molar fluxes given in \eqref{Ch8:eqn:jn:neg}, the first line in the right hand side of \eqref{Ch8:eqn:n-der} is canceled. 
Therefore, by designing the observer gains as
\begin{align}
Q_{-}( r_{\rm{p}}(t) ) 
& = - 
\frac{ a_{s,+}L_{+} }{ a_{s,-}L_{-} } \left( Q ( r_{ \rm{p} }(t)  ) + 
\frac{ \kappa }{ R_{ \rm{p},+ }^2 } \widehat{r}_{ \rm p}(t)^2 \right), \\
P_{-}(r_{\rm p}(t)) 
& =   - \frac{ \epsilon_{\rm{s},+}L_{+}}{\epsilon_{\rm{s},-}L_{-}} \frac{3}{R^{3}_{\rm{p},+}} 
\left[ \int_{ \widehat{r}_{\rm{p}}(t)}^{R_{\rm{p,+}}} P(r_{\rm{p}}(t)) r^{2}dr  \right] , 
\end{align}
one can show that $\frac{ d{ \widehat{n}_{ {\rm Li},+} } }{dt} =0$ from \eqref{Ch8:eqn:n-der}. Hence, the observer error in the negative electrode approaches to zero uniformly in space with the help of Theorem \ref{Ch8:thm:main}.

\subsection{Numerical Simulation}

\begin{figure}
	\begin{center} 
	\includegraphics[width=0.9\linewidth]{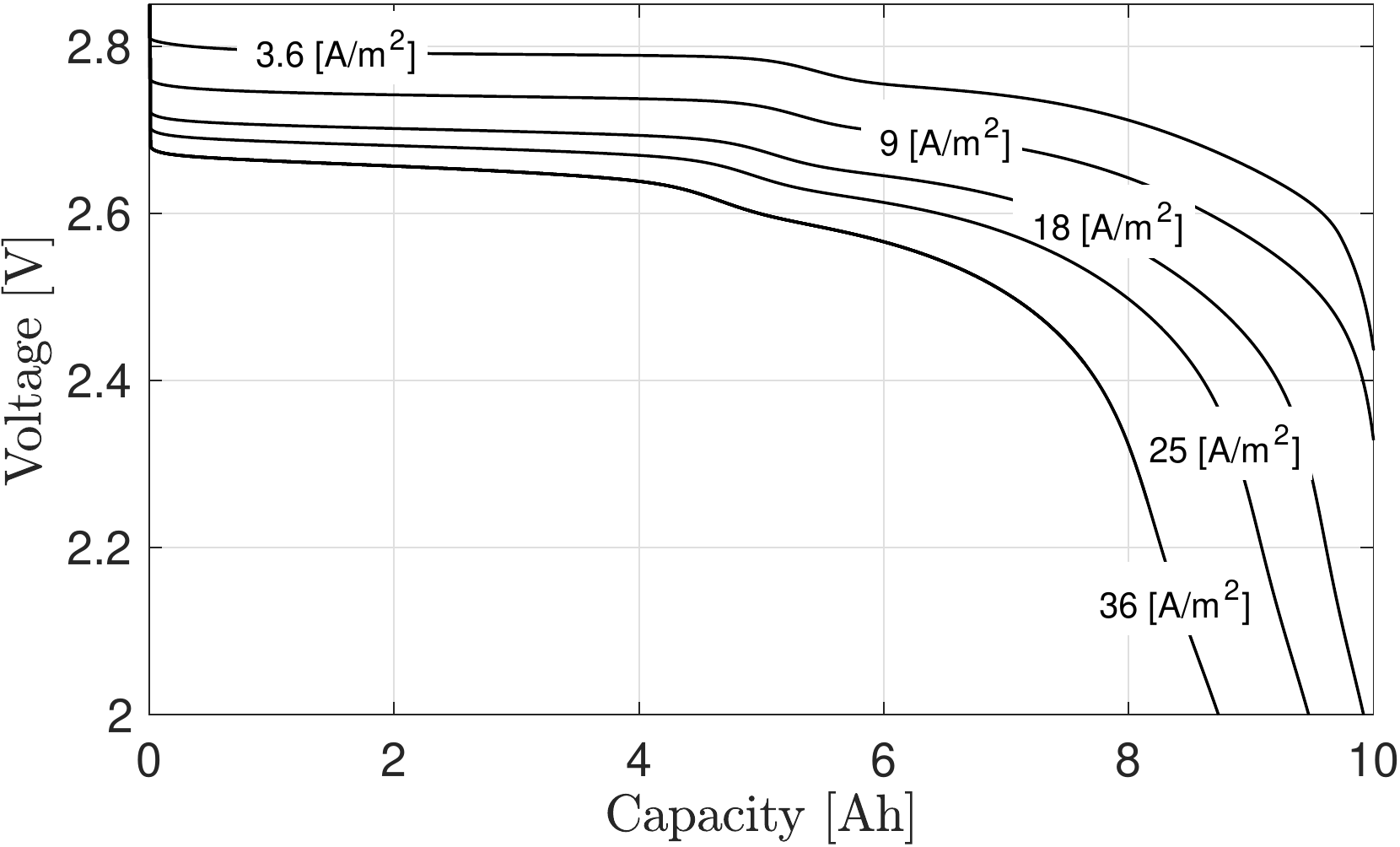}
	\caption{ Voltage plot for different (constant) current discharge inputs, which shows the analogous behavior to \cite{srinivasan04}.}
	\label{Ch8:discharge}
	\end{center}
\end{figure}

\begin{figure}
	\begin{center}
	\includegraphics[width=0.9\linewidth]{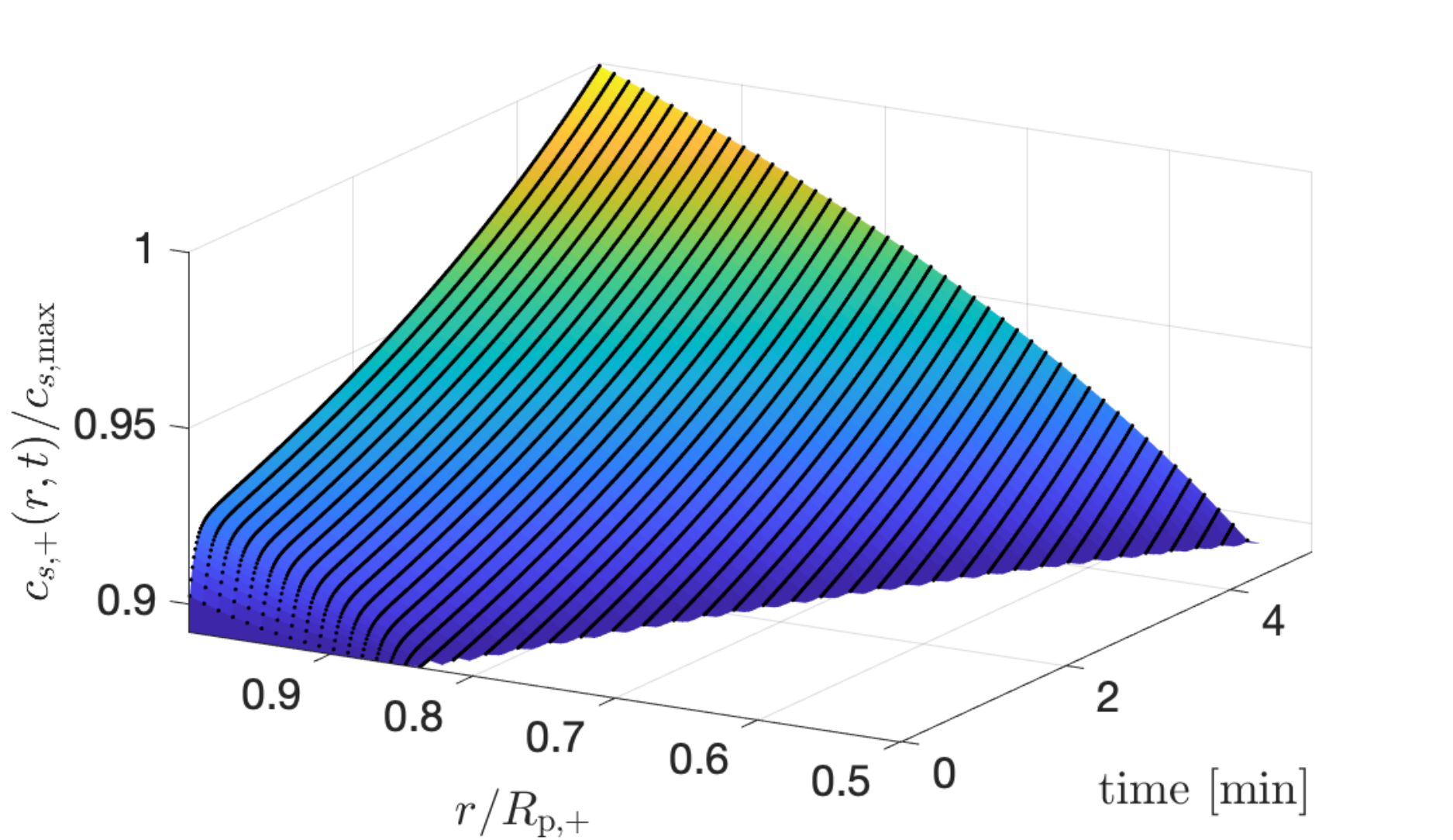}
	\caption{ Normalized concentration of lithium ions in a growing $\beta$-phase region.
		The plot corresponds to a $5[ \rm{min}  ]$ simulation of constant $5\rm[C-rate] $ discharge. The plot does not show 
		the $\alpha$-phase portion of the concentration since it is assumed to be constant. }
	\label{Ch8:betaphase}
	\end{center}
\end{figure}

\begin{figure}
	\begin{center} 
	\includegraphics[width=0.9\linewidth]{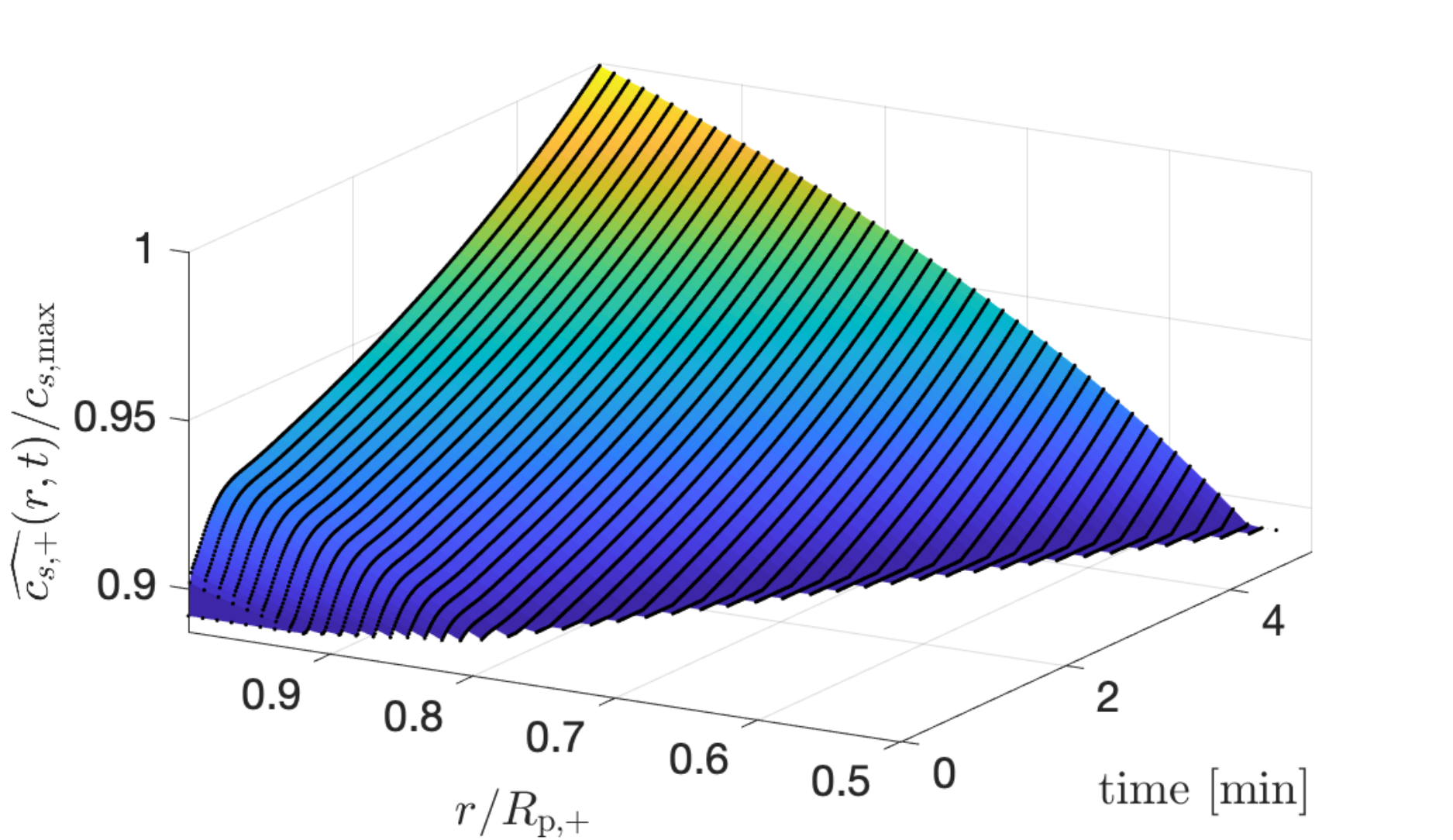}
	\caption{Estimate of the concentration of lithium ions in the positive particle. Starting from the initial error, the estimated profile converges to the true profile in Fig. \ref{Ch8:betaphase}.  }
	\label{Ch8:estimation:state}
	\end{center} 
\end{figure}

\begin{figure}
	\begin{center} 
	\includegraphics[ width=0.9\linewidth]{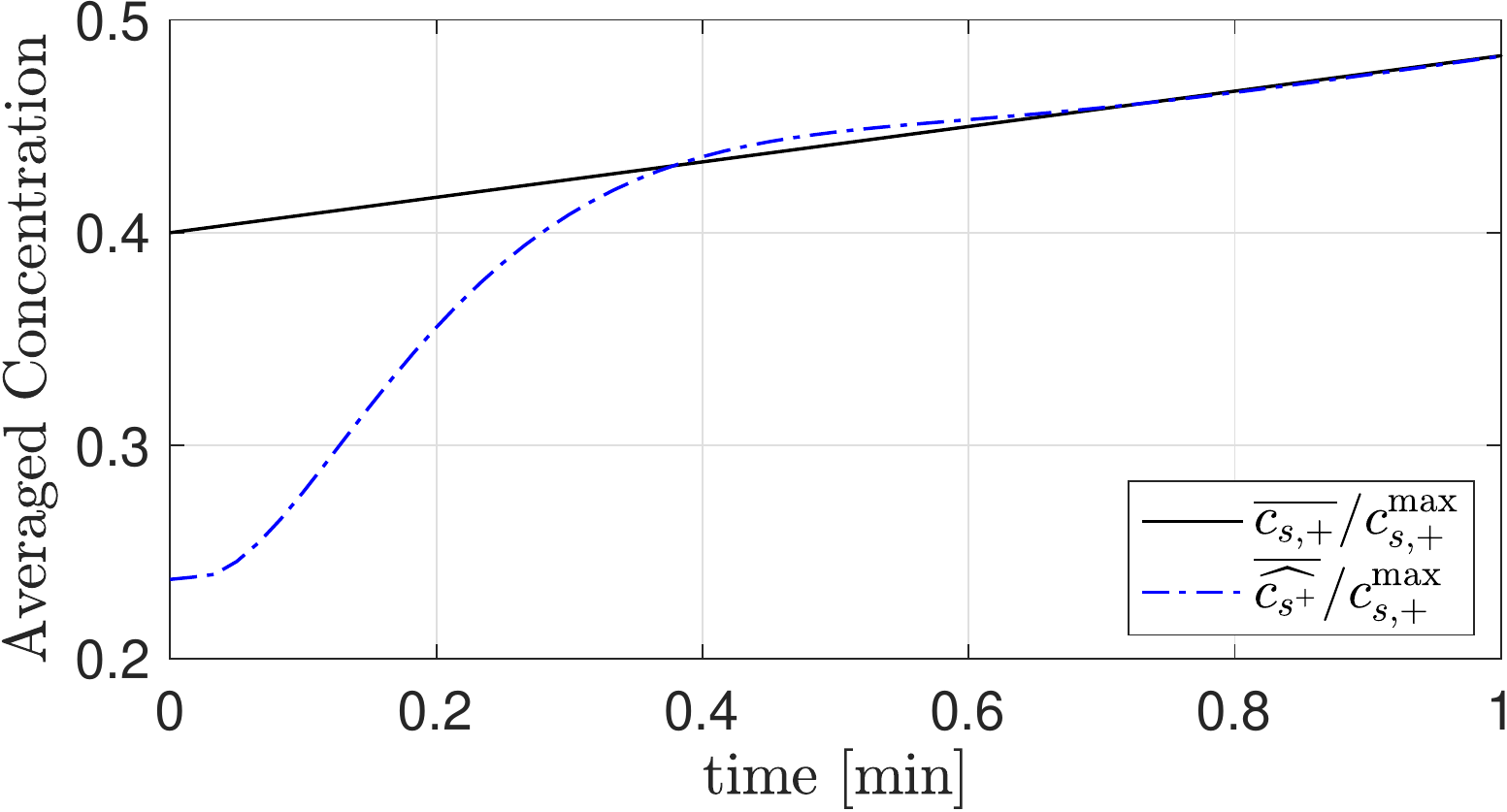}
	\caption{  Averaged concentration of true value (black solid) and estimated averaged concentration (blue dashed) in the positive particle normalized by the maximum concentration. }
	\label{Ch8:estimation:average}
	\end{center}
\end{figure}

\begin{figure}
	\begin{center}
	\includegraphics[ width=0.9\linewidth]{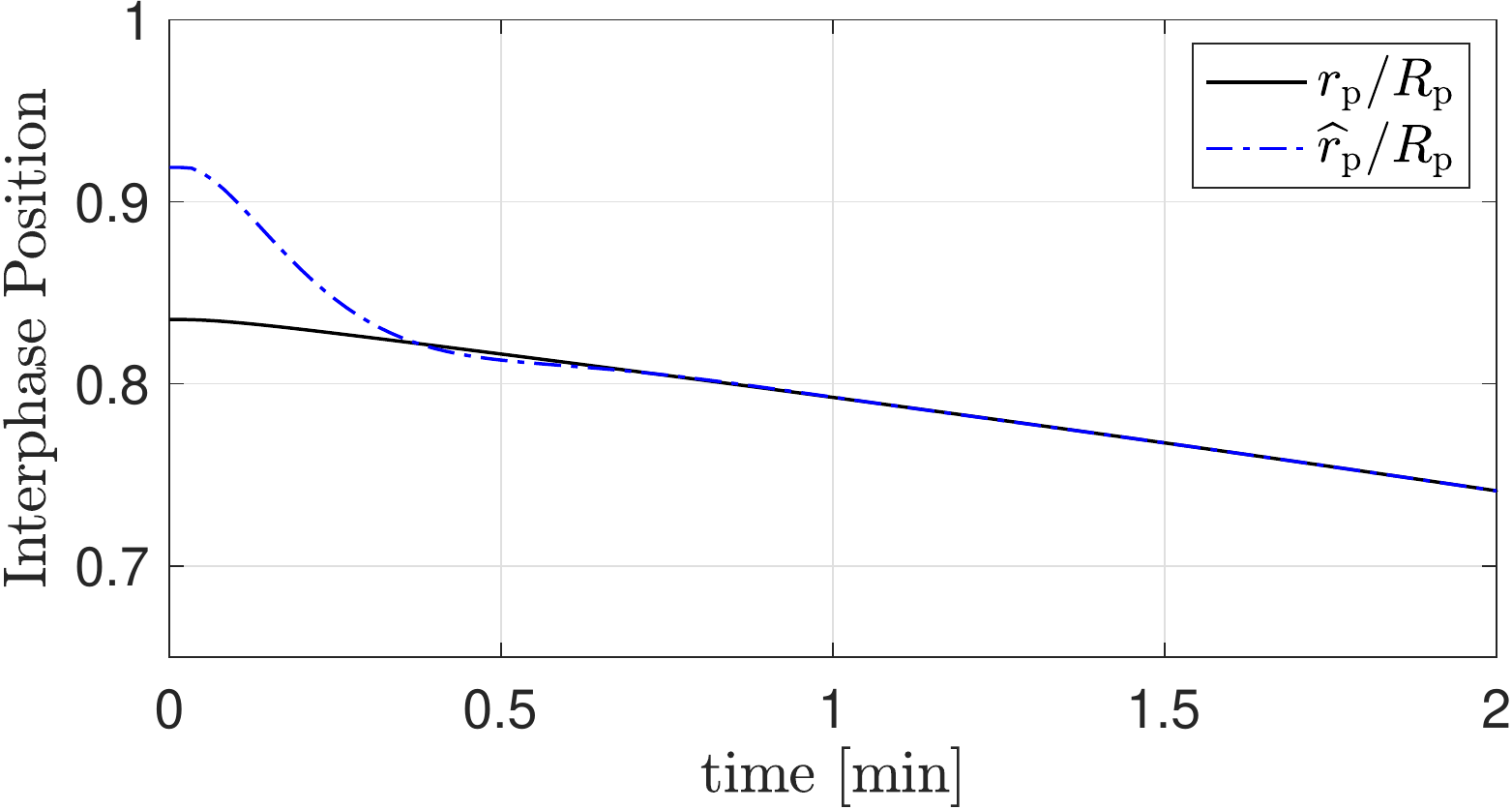}
	\caption{The estimated interface position becomes the same value as the true interface position after 0.5 [min]. }
	\label{Ch8:estimation:interface}
	\end{center} 
\end{figure}

\subsubsection*{Test with Constant Discharge Input}

To test the observer we run a numerical example with a constant discharge current of $5$ [C-rate].
We are assuming $ c_{\mathrm{ss},+} $ is available directly from measurements to be used as output error injection
in the observer. In practice, this quantity could be estimated from measurements. Figure \ref{Ch8:estimation:state} shows the estimated concentration of lithium ions in $\beta$-phase in the positive particle; 
one can compare this to the true concentration in Figure \ref{Ch8:betaphase}. Figure \ref{Ch8:estimation:average} shows the averaged concentration in the positive particle, both true value (black) and estimated value (blue). Convergence of the estimate to the true value is achieved within 0.8 [min], a relatively short time. Furthermore, Fig. \ref{Ch8:estimation:interface} shows the time evolution of the moving interface of the both true value (black) and estimated value (blue), which also illustrates the convergence of the estimate to the true value. Note that SoC is directly proportional to the averaged concentrations; then the importance to evaluate the estimation of this quantity.

\subsubsection*{Comparison with the Extended Kalman Filter} 
\begin{figure}[t]
\begin{center}
\subfloat[Estimation without noise. ]{\includegraphics[width=0.5\linewidth]{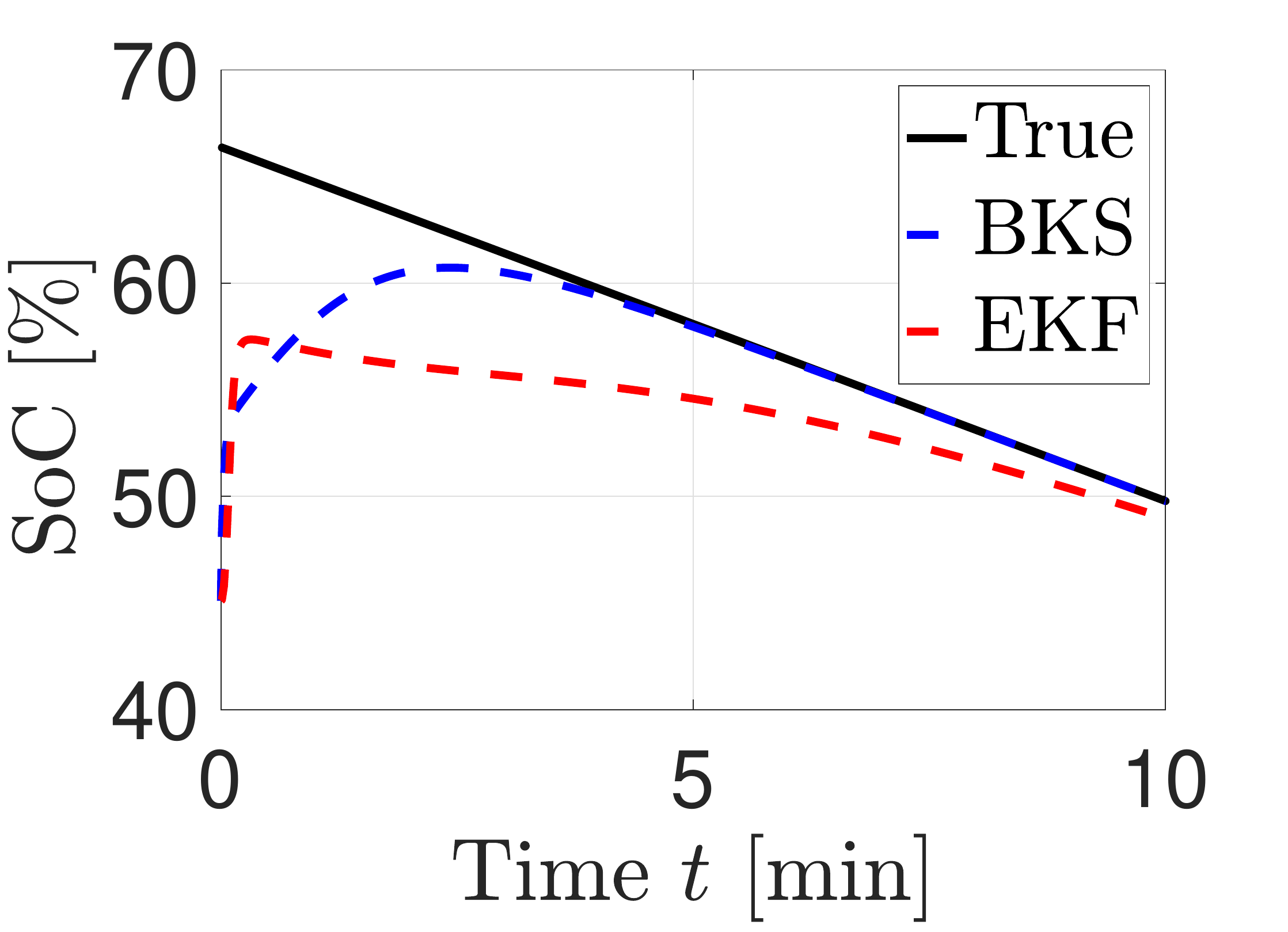}}
\subfloat[Estimation with measurement noise. ]{\includegraphics[width=0.5\linewidth]{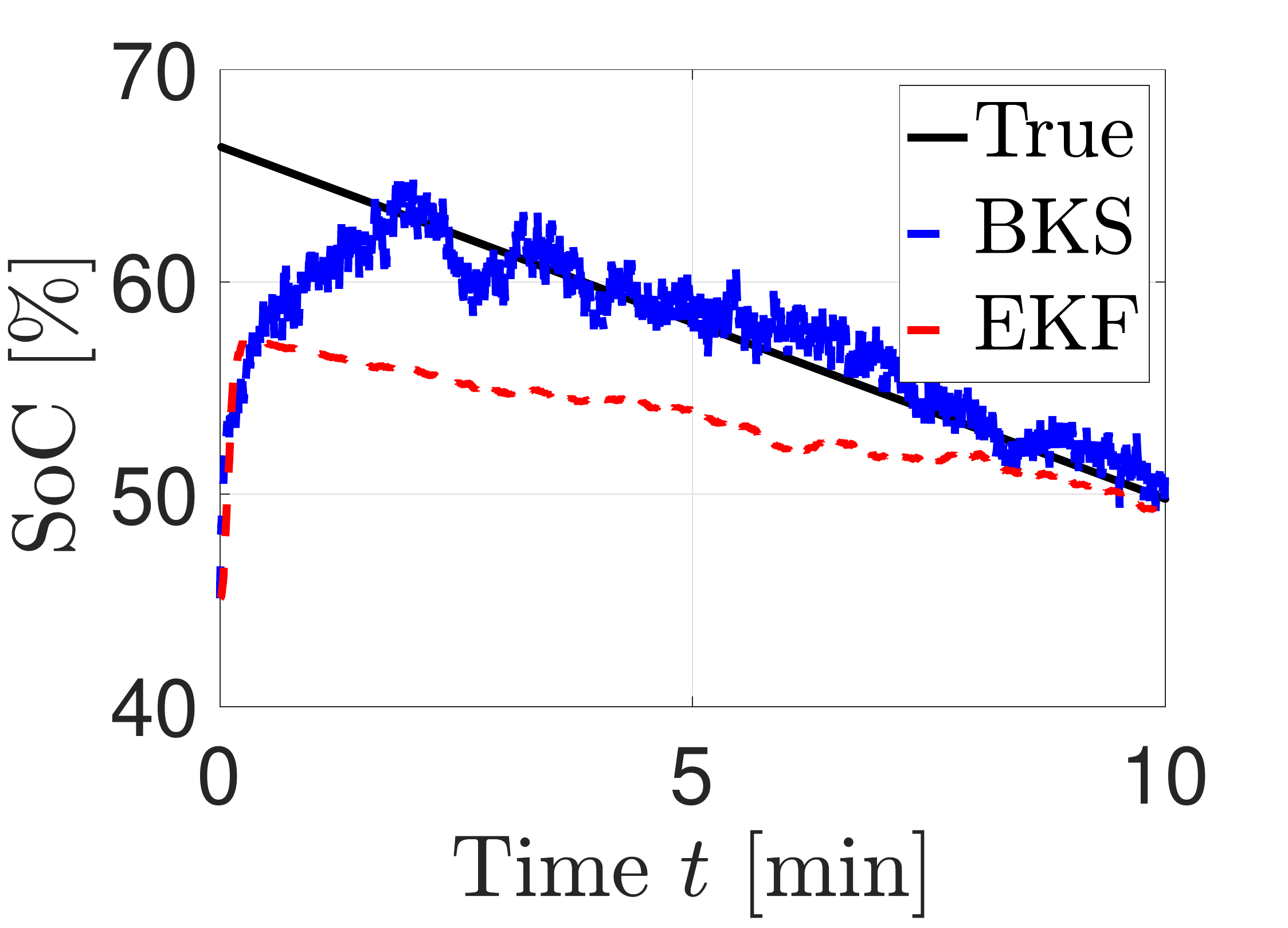}}
\caption{ Comparison of SoC estimation between the proposed BKS observer and the EKF. This sample simulation illustrates that the proposed BKS is superior in convergence speed, while the EKF is superior in noise attenuation. Note that each method can reduce the drawback by appropriately tuning the free parameters. }
\label{Ch8:fig:comparison}
\end{center}
\end{figure}

Since the spatial discretization of the diffusion equations is performed for computation of the electrochemical model, we can apply the Extended Kalman Filter\index{extended Kalman filter} (EKF) to the reduced-order model as another approach for SoC estimation. Here, we compare the performance of SoC estimation\index{SoC estimation} between the proposed backstepping (BKS) observer and the EKF with incorporating a measurement noise. 

Fig. \ref{Ch8:fig:comparison} shows a simulation result of SoC estimation via the BKS observer (blue) and the EKF (red). The initial SoC in the model is around 66 \%, while the initial SoC in the estimator is around 46 \% for both BKS and EKF. Fig. \ref{Ch8:fig:comparison} (a) depicts the result under the noise-free measurement, which illustrates that the BKS estimate converges and almost stays at the true value within 5 minutes, while the estimate by the EKF converges quickly first but does not stay at the true value even after 10 minues. Next, we incorporate the Gaussian noise in the measured value of the surface concentration, and compare the performance in  Fig. \ref{Ch8:fig:comparison} (b). The result illustrates that the BKS estimate still converges to the true value but it accompanies a noisy signal, while the EKF estimate has less noisy signal. Hence, in this one sample simulation, we observe that the BKS estimator is superior in convergence speed, while the EKF estimator is superior in noise attenuation. However, by lowing the gain parameters $(\lambda, \kappa)$ in the BKS estimator, we observe that the amplitude of the noisy estimate can be reduced in exchange for the convergence speed. Moreover, the EKF estimate also can be improved to accelerate the convergence speed by appropriately tuning the free parameters. Hence,  there is an essential tradeoff between the convergence speed and the noise attenuation for both methods, and it is not appropriate to address which method is superior in general. Nevertheless, one of the advantages of the proposed BKS estimator is to have only two free parameters to tune for any given discretization number, while the EKF algorithm requires a tuning of covariance matrices in which the number of free parameters is increased as the discretization number is increased.

% \subsection{Comments and Remarks} \label{sec:conclusion}

% This section develops the estimation algorithm for SoC via electrochemical model-based moving-boundary PDE observer, and provides the numerical study illustrating the desired performance of the proposed method. Towards a complete SoC estimation algorithm for lithium ion batteries with phase transition materials, an extension the existing SoC estimation algorithms from SPM
% to complex electrode settings will be considered, as was already achieved for electrodes with multiple active material \cite{camacho2016}. It was noted in \cite{SN2004} that two different
% particles sizes are needed to correctly model LFP electrodes, this correction can be added to our results following \cite{camacho2016}. One of the main assumptions for the model in this section is the restriction to only two coexisting phases in a single particle reduced further to a single phase problem by assuming a constant core phase. The relaxation of this assumption could be achieved through designing the state observer of concentration of lithium-ions in \emph{two phases} together with the estimation of the interface position. Furthermore, the robustness of the estimator's performance under some additive measurement noise can be studied in terms of input-to-state stability\index{input-to-state stability} (ISS) following \cite{Leobardo2018}. 

\begin{table}[t] \label{Ch8:table:param}
	\caption{ Parameters of LFP used in the simulation. }
	\centering
	\renewcommand{\arraystretch}{1.2}
	\begin{tabular}{rccc}
		\toprule[0.5mm]
		\multicolumn{4}{l}{\bf{ Parameters }} \\ & Negative & Separator &  Positive  \\ 
		\midrule
		$ L[\rm m] ^a$
		& $ 50 \times 10^{-6}$ 
		& $ 25 \times 10^{-6} $ 
		& $ 74 \times 10^{-6} $ \\
		
		$ c^{\rm max}_{\rm s} [\rm{mol}/m^3]^a$ 
		& $ 27760 $ 
		&
		& $ 20950 $ \\
		
		$ c_{\mathrm{s}, \alpha} [\rm{mol}/m^3]^b$ 
		&
		&
		& 0.0480$ \times c^{\rm max}_{\rm s,+}$ \\	
		 
		$ c_{\mathrm{s}, \beta  } [\rm{mol}/m^3]^b $  
		&
		&
		& 0.8920$ \times c^{\rm max}_{\rm s,+} $ \\		
																				  
		$ R_{\rm p} ~ [\rm m]^a $
		& $ 11 \times 10^{-6} $ 
		&
		& $ 52 \times 10^{-9} $ \\
		
		$ D_{\rm s} ~ [\rm m^2/s]^a$
		& $ 9 \times 10^{-14} $ 
		&
		& $ 8  \times 10^{-18} $ \\
		
		$ \epsilon_{\rm s}~[-]^{\rm{a}}$
		& $ 0.33 $
		&
		& $ 0.27 $ \\ 
		
		$ R_{\rm f} ~ [\rm \Omega m^2]^b $
		& $ 1 \times 10^{-5} $
		&
		& $ 0 $ \\ 
		
		$ R_{\rm c} ~ [\rm \Omega m^2]^b $
		& $ 0 $
		&
		& $ 6.5 \times 10^{-3} $ \\

		$ k ~ [\rm m^{2.5}/mol^{0.5}s ]^a $
		& $ 3 \times 10^{-5} $ 
		&
		& $ 3 \times 10^{-17 } $ \\  
		\midrule												     
	\end{tabular}\\
	\begin{tabular}{rccc}
		\multicolumn{4}{l}{\bf{Other Parameters and Physical Constants }}  \\
		\midrule
		$ A ~ [\rm m]^b $ & 1  \\ 
	    $ F ~ [\rm As/mol] $ & $ 96487  $  	\\ 
		$ R ~ [\rm J/Kmol ] $ & 8.314472 \\ 
		$ T ~ [\rm K]^b$ &  298 \\
		$ c_{\textrm{e}} ~ [\rm mol/m^3]^a $ & $ 1 \times 10^{3} $ \\ 
		$ \alpha_{\rm a}, \alpha_{\rm c} ~ [\rm -]^a  $ & 0.5 	\\														   
		\midrule 		
	\end{tabular} 
	\\
	$^{\rm{a}}$ borrowed from \cite{SN2004} \\
	$^{\rm{b}}$ assumed \\	
\end{table}

%\section{Notes and References}

%These further investigation will be addressed in our future work.

%%%%%%%%%%%%%%%%%%%%%%%%%%%%%%%%%%%%%%%%%%%%%%%%%%%%%%%%%%%%%%%%%%%%%%%%%%%%%%%%%%%%%%%

\section{Conclusion and Open Problems}

In this tutorial article, we have presented the state estimation of the Stefan PDE system modeling the dynamics of phase change phenomena of a material, and its applications to polar ice and lithium-ion batteries. The estimator is designed by an backstepping observer, which is given by a copy of the plant plus output error injection, where the observer gain is derived explicitly by solving a gain kernel equation of the state transformation. The convergence of the designed observer to the plant state from the Stefan system is ensured by Lyapunov analysis. The simulation results conducted for the thermodynamic model of Arctic sea ice illustrate 
%under suitable choice of the gain parameter, the temperature distribution in the Arctic sea ice converges to the true temperature profile within three days without overshooting, which is approximately ten times faster than the the 
the robust performance of the designed observer with respect to the neglected salinity effect and parameter uncertainty, where the convergence of the estimated temperature distribution to the true temperature is achieved within three days. The simulation performed for the electrochemical model of the lithium-ion batteries with phase transition materials has shown that the reduction of the error of more than 15 \% in the SoC estimate is achieved within five minutes even in the presence of sensor noise. 

There are various exciting open problems for the state estimation of the Stefan system, both from control-theoretic and application-driven perspectives. We summarize them in the following list for the control-theoretic problems: 
\begin{itemize}
\item sensor-delay compensation in the observer for the Stefan system (see \cite{ahmed2019observer}), 
    \item adaptive observer design for the Stefan system (see \cite{ahmed2015adaptive,ahmed2017adaptive}),
    \item observer for two-phase Stefan system under a single-boundary measurement (see \cite{liu2016backstepping}), 
    \item sampled-data observer for the Stefan system (see \cite{Ali17PDEODE,Karafyllis19sampled}), 
    \item prescribed-time observer design for the Stefan system (see \cite{Drew20,Drew19}), 
%\end{itemize}
and for the application-driven problems: 
%\begin{itemize} 
\item cancer treatment via cryosurgery (see  \cite{Rabin1998, Rabin03,Kumar17}), 
\item spreading of invasive species in ecology (see \cite{Du2010speading}),  
\item information diffusion on online social networks (see \cite{Wang20book}), 
\item domain walls in ferroelectric thin films (see \cite{mcgilly2015}), 
\item Black-Scholes model of American option pricing (see \cite{Chen2008}). 
\end{itemize} 

While the list of control-theoretic problems have all been considered for PDE and PDE-ODE systems on fixed domains, virtually all of them are unexplored for PDEs with moving boundaries. Moreover, the list of application-driven topics  owns a significant impact in the real-world problems, by the utilization of the PDE-based estimation method. This tutorial review provides supporting technical materials to tackle those challenging topics.

%\input{Section_1} %Introduction
%\input{Section_2} %Model
%\input{Section_3} %Estimator designs
%\input{Section_4} %Sea ice applications
%\input{Section_5} %Battery applications
%\input{Section_6} %Conclusion and Open Problems



\bibliographystyle{elsarticle-harv} 


\bibliography{BIB_ARC_21.bib}

\end{document}